\documentclass[sn-mathphys, Numbered]{sn-jnl}
\usepackage{dsfont}
\usepackage{graphicx}%
\usepackage{amsmath, amssymb, amsfonts, amsthm}
\usepackage[title]{appendix}%
\usepackage{xcolor}%
\usepackage{bm}
\usepackage{array, hhline}
\usepackage{enumerate}
\usepackage{ifthen}
\usepackage{nicefrac}
\usepackage{manyfoot, multirow, soul}%
\usepackage{fontenc, lmodern}
\definecolor{rose}{rgb}{1, 0, 0.6}
\newcommand{\ignore}[1]{}

\newcommand{\kav}[1]{{\color{red}{#1}}}

\newcommand{\indic}[1]{{\mathds{1}_{\{#1\}}}}

 \makeatletter
\newcommand\appendix@section[1]{%
  \refstepcounter{section}%
  \orig@section*{Appendix \@Alph\c@section: #1}%
  \addcontentsline{toc}{section}{Appendix \@Alph\c@section: #1}%
}
\let\orig@section\section
\g@addto@macro\appendix{\let\section\appendix@section}
\makeatother
 
\newtheorem{thm}{Theorem}

\newtheorem{defn}{Definition}
\newtheorem{lemma}{Lemma}

\newcommand{\Prob}{{\mathbb{P}}}
\newcommand{\E}{{\mathbb{E}}}
\newcommand{\R}{{\mathcal{R}}}
\newcommand{\WB}{{\mathcal{D}}}
\newcommand{\eL}{{\phi_L^*}}
\newcommand{\dpr}{$\mathsf{dpr}$}
\newcommand{\idp}{$\mathsf{IDP}$}
\newcommand{\BR}{{\mathbb B}}
\newcommand{\eU}{{\phi_R^*}}

\newcommand{\bphi}{{\bm \phi}}
\newcommand{\phimr}{\eU} 

\newcommand{\VPhi}{{\bm \phi}}

\newcommand{\PB}{{\cal B}}

\newcommand{\DA}{{\cal U}}

\newcommand{\MR}{{\cal M}}
\newcommand{\du}{{d_u}}
\newcommand{\db}{{d_b}}
\newcommand{\dm}{{d_m}}
\newcommand{\eiw}{{e}}

\newcommand{\eow}{{e}}

\newcommand{\niw}{{n_{i}}}
\newcommand{\nir}{{r_{i}}} 

\newcommand{\lamol}{{{\lambda}_{1}({\VPhi; 0})}}

\newcommand{\lamtl}{{{\lambda}_{2}({\VPhi; 0})}}
\newcommand{\N}{{\bar N}}

\newcommand{\lmo}[1]{{\lambda_{1}({#1})}}
\newcommand{\lmt}[1]{{\lambda_{2}({#1})}}
\newcommand{\lmi}[1]{{\lambda_{i}({#1})}}

\newcommand{\lamo}{\lmo{\VPhi;\beta}}
\newcommand{\lamt}{\lmt{\VPhi;\beta}}
\newcommand{\lami}{\lmi{\VPhi;\beta}}

\newboolean{showcomments}
\setboolean{showcomments}{false}

\usepackage[figuresright]{rotating}
\usepackage{caption}
\usepackage{subfig}
\usepackage{multirow}
\usepackage{booktabs}

\usepackage{ulem}

\newcommand{\extra}[1]{\ifthenelse{\boolean{showcomments}}{{\textbf{(Extra:  #1)}}} {}  }

\newcommand{\tw}[1]{\ifthenelse{\boolean{showcomments}}{{\textbf{(TW says:  #1)}}} {}  }
\newcommand{\rev}[1]{{\textcolor{black}{#1}}}

\begin{document}

\title{On the interplay between pricing, competition and QoS in ride-hailing}

\author*[1]{\fnm{Tushar Shankar} \sur{Walunj}}\email{tusharwalunj@iitb.ac.in}

\author[1]{\fnm{Shiksha} \sur{Singhal}}\email{shiksha.singhal@iitb.ac.in}

\author[2]{\fnm{Jayakrishnan} \sur{Nair}}\email{jayakrishnan.nair@ee.iitb.ac.in}

\author[1]{\fnm{Veeraruna} \sur{Kavitha}}\email{vkavitha@iitb.ac.in}

\affil[1]{\orgdiv{Industrial Engineering \& Operations Research}, \orgname{IIT Bombay, India}}

\affil[2]{\orgdiv{Electrical Engineering}, \orgname{IIT Bombay, India}}


\ignore{\author{Tushar Shankar Walunj}
\ead{tusharwalunj@iitb.ac.in}
\author{Shiksha Singhal}
\ead{shiksha.singhal@iitb.ac.in}

\author{Jayakrishnan Nair}
\ead{jayakrishnan.nair@ee.iitb.ac.in}

\author{Veeraruna Kavitha}
\ead{vkavitha@iitb.ac.in}

\address{IEOR, Indian Institute of Technology Bombay, India}

{\let\thefootnote\relax\footnote{{The second author's work is partially supported by Prime Minister's Research Fellowship (PMRF), India.}}}
}
\abstract{
We analyse a non-cooperative game between two competing ride-hailing platforms, each of which is modeled as a two-sided queueing system, where drivers (with a limited level of patience) are assumed to arrive according to a Poisson process at a fixed rate, while the arrival process of (price-sensitive) passengers is split across the two platforms based on Quality of Service (QoS) considerations. As a benchmark, we also consider a monopolistic scenario, where each platform gets half the market share irrespective of its pricing strategy. The key novelty of our formulation is that the total market share is fixed across the platforms. The game thus captures the competition between the platforms over market share, with pricing being the lever used by each platform to influence its share of the market. \rev{The market share split is modeled via two different QoS metrics: (i) probability that an arriving passenger obtains a ride, and (ii) the average passenger pick-up time.} The platform aims to maximize the rate of revenue generated from matching drivers and passengers.

In each of the above settings, we analyse the equilibria associated with the game in certain limiting regimes. We also show that these equilibria remain relevant in the more practically meaningful `pre-limit.' Interestingly, we show that for a certain range of system parameters, no pure Nash equilibrium exists. Instead, we demonstrate a novel solution concept called an \textit{equilibrium cycle}, which has interesting dynamic connotations. Our results highlight the interplay between competition, passenger-side price sensitivity, and passenger/driver arrival rates.
}

\keywords{
BCMP queueing network, ride-hailing platforms, two-sided queues, Wardrop equilibrium, Nash equilibrium, cooperation
}

\maketitle

\section{Introduction}

The ride-hailing industry, exemplified by platforms like Uber, Lyft, and Ola, has gained tremendous popularity and growth in recent years. These platforms operate as two-sided matching systems, where passengers in need of a ride and willing drivers are matched in real time. The scale of these matching platforms is remarkable, with OLA alone boasting over 1.5 million drivers across 250 cities in India~\cite{enwiki:1158169707}.

In the ride-hailing space, we see strategic behavior on part of passengers as well as platforms.\footnote{In certain settings, drivers are also strategic, though this aspect is not addressed in this paper.} On one hand, passengers choose which platform to use based on not just its pricing, but also its (history of) driver availability. (Indeed, today, technology has considerably cut the friction associated with switching between platforms.) On the other hand, platforms seek to maximize their revenues while competing with one another for market share; each platform must optimize its pricing to balance its market share and its revenue per-ride. It is this interplay between pricing, competition, market segmentation and QoS that we seek to analyse in this paper.

Although there has been some recent research on competing (two-sided) matching platforms, most studies overlook the crucial \textit{queueing} aspects that are inherent in practical ride-hailing systems, including random passenger/driver arrivals, the likelihood of driver availability, as well as random transit durations. Moreover, most existing studies rely on numerical computations of the system equilibria to draw their insights on impact of the strategic interaction between competing platforms (we provide a detailed review of the related literature later in this section). This paper aims to bridge this gap by formally analyzing a non-cooperative game between two ride-hailing platforms that compete for market share. The platforms compete via their pricing strategies, to serve a passenger base which is both impatient as well as price and delay sensitive.

Specifically, we assume symmetric platforms, a single geographical zone of operation, and static (not state-dependent) pricing. The market share of each platform is modeled as a Wardrop equilibrium~\cite{WE} (defined in terms of a certain QoS metric). We then characterize the equilibria associated with this game, approximating the payoff functions along certain scaling regimes. These equilibria shed light on the influence of passenger price sensitivity and passenger/driver arrival rates on pricing, platform surplus, and passenger QoS in the presence of inter-platform competition.

Our main contributions are as follows:\footnote{A preliminary version of this paper was presented at the Allerton conference in 2022~\cite{TR}.}

\begin{itemize}
    \item We model each platform as a BCMP network (\cite{BCMP}), admitting a product form stationary distribution. This model captures passenger-side price sensitivity, queueing (and impatience) of drivers, and transit times.
    (This BCMP modeling approach is quite powerful, and also allows for state-dependent pricing and probabilistic routing across multiple zones, though these features are not used in our analysis here.) 
    \item For benchmarking (i.e., to capture the impact of inter-platform competition), we first consider a \textit{monopolistic} scenario, where each platform gets half the market share, irrespective of its pricing strategy. We characterize the optimal pricing strategy of each platform under a certain limiting regime where driver patience grows unboundedly---we refer to this scaling regime as the Infinite Driver Patience (\idp) regime. 
    \item To capture competition between platforms, we model market share bifurcation between the platforms in the form of a Wardrop equilibrium, where the passenger base splits in manner that seeks to equalize a certain QoS metric across the platforms. We consider two reasonable QoS metrics: \rev{(i) the stationary probability that an arriving passenger is not served (which also accounts for those cases where a passenger declines a ride because the price quoted by the platform was unacceptable), and (ii) the average passenger pick-up time (where the delay/dis-utility in the event of not obtaining service is taken to be a suitably large constant).} \rev{Under each of the above passenger-side QoS metrics, we characterize the equilibria of the non-cooperative game between the platforms under certain limiting regimes; in fact it turns out the equilibria under the two QoS metrics coincide. In the case of the former QoS metric, we employ an \text{infinite driver patience} (\idp) scaling regime, where driver patience is scaled to infinity; in the case of the latter QoS metric, we additionally scale the delay/dis-utility in the event of not obtaining service to infinity.}
    
    Comparing these (symmetric) equilibria with the monopolistic scenario described above sheds light on the impact of inter-platform competition on pricing, QoS and platform revenues. Specifically, we find that competition drives platforms to deviate from monopoly pricing when \textit{passengers are scarce relative to drivers}. 
    When this happens, the equilibrium price deviates so as to enhance the very QoS metric that governs the market segmentation. 
    Another consequence of this deviation from monopoly pricing is that the utility (revenue rate) of each platform falls. In other words, competition between platforms benefits the passenger base and `hurts' the platforms. 
    
    \item From a game theoretic standpoint, we show that for a certain range of system parameters, no pure Nash equilibrium (NE) exists. Instead, we demonstrate a \textit{mixed} NE---the absence of a pure NE stems from certain discontinuities in the limiting payoff functions in the \idp\ regime. Interestingly, under the same set of system parameters, we also discover a novel solution concept that we refer to as an \textit{equilibrium cycle}. Specifically, in this setting, each platform has the incentive to set prices/actions within a certain interval, though for each such price/action pair, at least one platform has the incentive to deviate to a different price/action within the same interval. The equilibrium cycle thus suggests an \textit{indefinitely oscillating pricing dynamic}.
    
    \item While our explicit equilibrium characterizations are made under the certain limiting regimes, we also show that these equilibria are meaningful in the more practical `pre-limit,' where drivers are highly (though not infinitely) patient. Specifically, we show that all our equilibria (including the equilibrium cycle) are $\epsilon$-equilibria in the pre-limit.
\end{itemize}

The remainder of this paper is organized as follows: After briefly surveying the related literature below, we describe our system model in Section~\ref{sec:model}. We then analyze the monopolistic setting in Section~\ref{sec:monopoly}. In Section~\ref{sec:overall_bp}, we consider the service availability metric, which depends on both (i)~driver unavailability and (ii)~the price quoted by the platform. \rev{In Section~\ref{sec:delay}, we consider the mean pick-up delay QoS metric.} In Section~\ref{sec:comparison}, we compare the prices and platform utilities characterized under the three scenarios above, both theoretically and numerically. In Section~\ref{sec:comparison}, we also demonstrate the applicability of our limiting regime equilibria in the pre-limit via numerical experiments. Finally, we conclude in Section~\ref{sec:conclusion}. 

\subsection*{Literature Review}

There is a considerable literature on two-sided queues, where customers/jobs as well as servers (drivers, in the specific context of ride-hailing platforms) arrive into their respective queues and get `matched' over time. Papers that have analysed a single two-sided platform from various perspectives (like optimal pricing, matching algorithms, fleet sizing, heavy traffic scaling, etc.) include \cite{sood2018pricing, siva,uber, BCMP_routing, sun2019optimal, johari, vaze2022non, mahavir2022heavy, feng2021we}. However, in this survey, we restrict attention to the (relatively recent) literature on \textit{competition} between two-sided platforms, which is the primary focus of the present paper. 

The literature on competing two-sided platforms can be categorized into (i) papers that consider \textit{single-homing}, wherein drivers (or just cars, in the case of autonomous vehicles) are attached to a platform a priori, and (ii) papers that consider \textit{multi-homing}, wherein drivers (or the owners of autonomous vehicles) can choose which platform to work with. We discuss these categories separately; note that present paper falls in the former category.

Papers that analyse competing two-sided platforms under the single-homing assumption include \cite{siddiq2022ride, randall_berry2022competition, sen2023pricing}. In~\cite{siddiq2022ride}, the authors analyse a pricing game between two competing platforms, capturing autonomous vehicles (AVs), conventional cabs, and price-sensitive passengers. Under a linear demand model, a explicit equilibrium characterization is provided; the authors additionally examine the impact of who owns the AVs (the platforms, or autonomous and strategic individuals) on the equilibrium.
In~\cite{randall_berry2022competition}, the authors analyse a game between two competing platforms operating over multiple zones. The game proceeds over $T$ time steps; each platform optimizes its initial fleet placement, and its pricing at each time step. While the generalized Nash equilibria of this game are hard to compute in general, the authors show that under certain conditions, the game admits a potential function, allowing the equilibrium to be computed by solving a certain convex optimization. Finally, in \cite{sen2023pricing}, the authors analyse a simpler game with two platforms and two zones; each platform chooses its price at both locations, as well as the rate of fleet re-balancing. An iterative algorithm for computing the equilibrium associated with this game is provided. Importantly, none of the above papers considers an explicit queueing model with random arrivals and random transit times; such a (BCMP) queueing model lies at the heart of the game analysed in the present paper.

Papers that analyse competition between two-sided platforms allowing for multi-homing include~\cite{ahmadinejad2019competition, zhang2021inter, bernstein2021competition}.
Specifically, in~\cite{ahmadinejad2019competition}, the authors analyse a game where both drivers as well as passengers choose which platform to use based on prices announced by both platforms. The main contribution here is the characterization of conditions for \textit{market failure}, where a \textit{tragedy of the commons} drives prices so low that drivers have no incentive to participate. \cite{zhang2021inter} proposes a more general model that also allows for pooling of rides between passengers. Via extensive numerical computations of the equilibrium, the authors argue that multi-homing is bad for all entities (platforms, drivers, and passengers), and propose mechanisms to discourage multi-homing. Finally, \cite{bernstein2021competition} analyses a game where the split of drivers and passengers between platforms is determined by Hotelling models, which additionally capture the impact of congestion (via a quadratic congestion model). This paper provides an explicit analytical characterization of the pricing equilibrium between the platforms, along with a sensitivity analysis of the equilibrium with respect to passenger arrival rate. Interestingly, \cite{bernstein2021competition} also arrives at a conclusion that is consistent with \cite{ahmadinejad2019competition, zhang2021inter}---multi-homing is bad for all parties involved. As before, none of the above mentioned papers considers the queueing aspects of ride-hailing, as alluded to earlier.

To summarize, to the best of our knowledge, the preceding literature on competing two-sided platforms disregards the queueing effects associated with ride-hailing. Moreover, most of the prior work (with the exception of \cite{siddiq2022ride, bernstein2021competition}) does not provide an explicit analytical characterization of the equilibrium of the game between platforms. In contrast, the present paper treats each platform as a BCMP queueing system; this model captures the queueing dynamics of passengers as well as drivers in a stochastic setting. Moreover, we provide an explicit analytical characterization of the equilibrium between platforms under a certain natural limiting regime, highlighting formally the interplay between competition, price sensitivity, and arrival rates.

\section{Model and Preliminaries}
\label{sec:model}

In this section, we describe our model for the interaction between two competing ride-hailing platforms, and state some preliminary results.

\subsection{Passenger arrivals}

We consider a system with a set~$\mathcal{N} = \{1,2\}$ of independent ride-hailing platforms, each of which is modeled as a single-zone two-sided queuing system that matches drivers and passengers. \rev{These platforms will be assumed to be symmetric in all respects, except for their pricing strategies.} Passengers arrive into the system as per a Poisson process with rate~$\Lambda$. This aggregate arrival process gets split (details in Subsection~\ref{sec:WE_QoS} below), such that the passenger arrival process seen by platform~$i$ is a Poisson process with rate~$\lambda_i$, where~$\sum_{i=1}^2 \lambda_i = \Lambda$.


When a passenger arrives into platform~$i,$ if there are no waiting drivers available, the passenger immediately leaves the system (a.k.a., gets blocked). On the other hand, if there are one or more waiting drivers, the passenger is quoted price~\rev{$\phi_i \in [0,\phi_h],$} where~$\phi_h$ denotes the maximum price the platform can charge.
The passenger accepts this price (and immediately begins her ride) with probability~\rev{$f(\phi_i);$} with probability~\rev{$1-f(\phi_i),$} the passenger rejects the offer and leaves the system. Note that the function~$f$ captures the price sensitivity of the passenger base.\footnote{An alternate interpretation is that~$1-f(\phi_i)$ denotes the probability that the passengers find a more suitable option (like public transport) outside the platform.} We make the following assumptions on the function~$f.$
\begin{enumerate}[{\bf A}.1]
\item $f(\cdot)$ is a strictly concave, strictly decreasing and differentiable function.
\item $0 < f(\phi) \le 1$ for all~$\phi \in [0,\phi_h]$ and~$f(0) =1$.
\end{enumerate}
The following extension of the inverse of~$f$ will be used in the statements of our results:
\begin{equation}\label{eqn_f_inverse}
    f^{-1}(x) :=
    \begin{cases}
      \phi_h  &\text{ for } 0 \leq x < f(\phi_h), \\
      \phi & \text{ s.t. } f(\phi) = x \text{ for } f(\phi_h) \leq x \leq 1, \\
      0 & \text{ for } x > 1.
    \end{cases}
\end{equation}

\subsection{Driver behavior}

Each platform has a pool of dedicated drivers, which arrive into the system according to a Poisson process of rate~\rev{$\eta$}. \rev{(Note that in our model, drivers are not strategic, and do not choose which platform to join.)} The drivers wait in an FCFS queue to serve arriving passengers. The number of waiting drivers in this queue at platform~$i$ is denoted by~$n_i$.

When an arriving passenger is matched with the head of the line driver, a ride commences and is assumed to have a duration that is
exponentially distributed with rate~$\rev{\nu}.$ At the end of this ride,
the driver rejoins the queue of waiting drivers (and therefore becomes
available for another ride) with probability~\rev{$p$}; with
probability~\rev{$1-p$}, the driver leaves the system.

\begin{figure}[t]
     \centering 
     \includegraphics[scale=0.32]{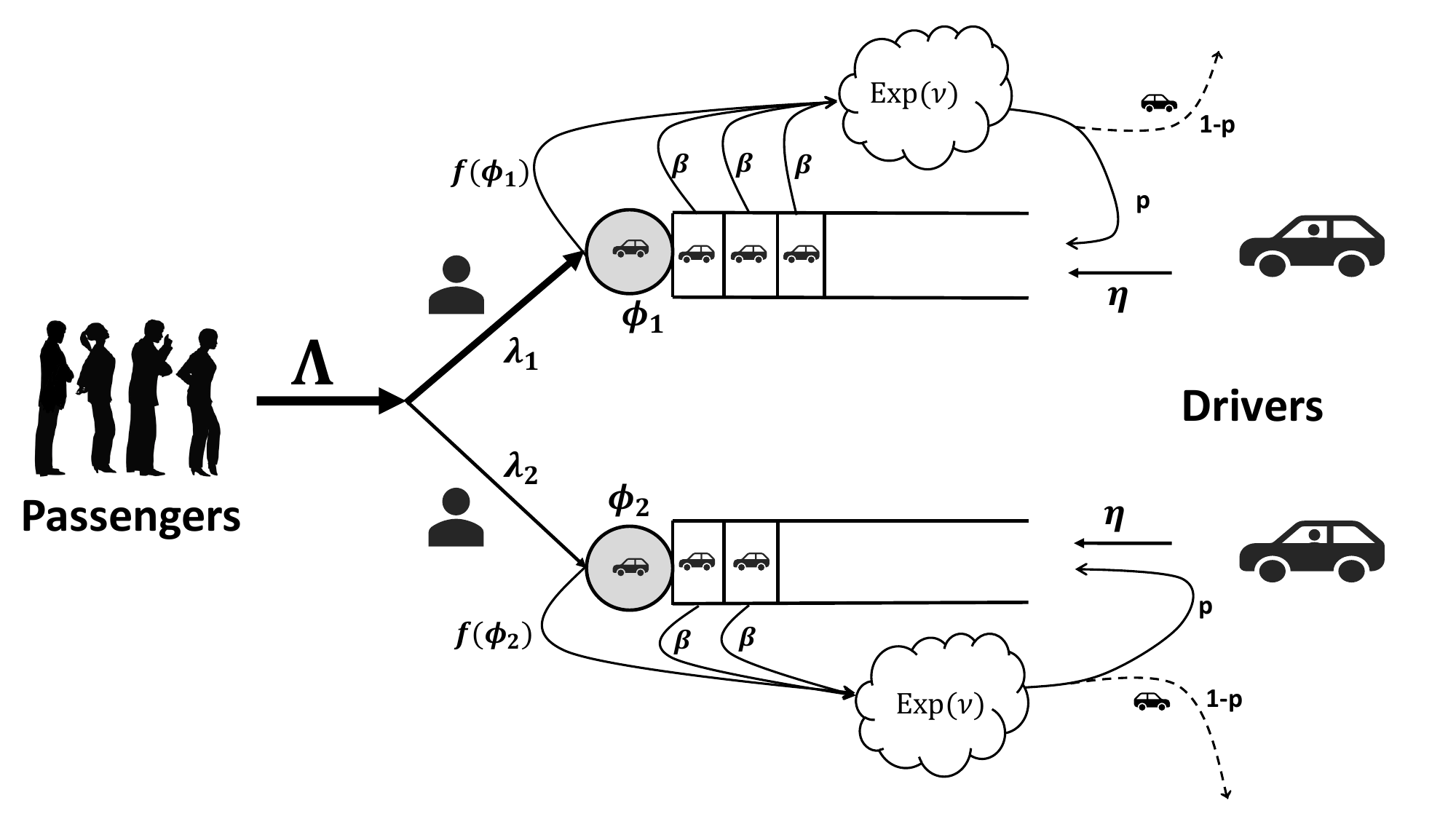}
\caption{Depiction of system model}
\label{fig:model}
 \end{figure}

Additionally, we also model driver abandonment from the waiting queue. Specifically, each waiting driver independently abandons the queue after an exponentially distributed duration of rate~\rev{$\beta$}. This
abandonment might capture a local (off platform) ride taken up by the
driver, or simply a break triggered by impatience. The duration of
this ride/break post abandonment is also assumed to be exponentially
distributed with rate~$\rev{\nu}$ (if drivers accept a ride off-platform,
it is reasonable to assume that the distribution of the duration of
this ride is the same as that of rides matched on the platform, given
that these rides are in the same zone); at the end of this duration, the
driver rejoins the queue with probability~\rev{$p$} and leaves the system
altogether with probability~\rev{$1-p$}. Figure~\ref{fig:model} presents a
pictorial depiction of our system model.

\subsection{BCMP modeling of each platform}

Under the aforementioned model, each platform can be described via a
continuous time (Markovian) BCMP network.\footnote{\rev{BCMP networks are a class of queueing networks that are known to admit a \textit{product-form stationary distribution}~\cite{BCMP}. They represent a generalization of the class of Jackson networks; see~\cite{balter}.}} Formally, the
state of platform~$i$ at time~$t$ is given by the tuple~$Z_{i,t} =
(N_{i,t},R_{i,t})$, where~$N_{i,t}$ is the number of waiting drivers,
and~$R_{i,t}$ is the number of drivers that are in the system but
unavailable to be matched with passengers (because of being in a ride,
or on a break). Realizations of the state of platform~$i$ are
represented as~$(\niw, \nir)$, and state space corresponding to
platform~$i$ is given by~$S_i
= \{(\niw, \nir)\colon \niw, \nir \in \mathbb{Z}_+\}$.

Each platform can be modelled using a BCMP network with two `service stations' as described in the following:
(i) service station 1 (SS1) is the queue of waiting drivers (the
occupancy of this queue is the first dimension of the state), and (ii)
service station~2 (SS2) is the `queue' of drivers on a ride/break (the
occupancy of this queue is the second dimension of the state). SS1 is
modeled as a single server queue having a state-dependent service
rate. Specifically, the (state dependent) departure rate
  from~SS1 is~\rev{$\lambda_i f(\phi_i) + n_i \beta$}. On the other
hand, SS2 is modeled as an infinite server queue, having exponential
service durations of rate~$\rev{\nu}.$ SS1 sees exogenous Poisson arrivals at
rate~$\eta,$ and departures from SS1 become arrivals into
SS2. Finally, departures from SS2 exit the system with
probability~\rev{$1-p$}, and join SS1 with probability~\rev{$p$}; see Figure~\ref{fig:model}.

Next, we describe the steady state distribution corresponding to the
above BCMP network. Define~$e :=
\nicefrac{\eta}{(1-p)},$ which is easily seen to be the effective driver
arrival rate seen by the service stations SS1 and SS2. The following
lemma follows from~\cite{BCMP} (see Sections~3.2 and~5 therein).
\begin{lemma} \label{lem_BCMP_Product_form}
The steady state probability of state~$s = (\niw, \nir)\in S_i$ is given by:
\begin{align*}
\pi_i(s) &= C_i
	\left[\frac{{e}^{\niw}}{\prod_{a=1}^{\niw} \left(\lambda_i f(\phi_i) + a {\beta} \right)}\right]
	\left[\left({\frac{{e}}{\nu}}\right)^{\nir} \left(\frac{1}{\nir!}\right) \right], \text{ with }  \\ 
C_i^{-1} 
&= \sum_{n=0}^{\infty} \left[\frac{{e}^{n}}{\prod_{a=1}^{n}(\lambda_i f(\phi_i)+ a {\beta})}\right]
\exp\left({\frac{{e}}{\nu}}\right).
\end{align*}
Here,~$C_i$ is the normalizing constant (we follow the convention that
$\prod_{a}^{k} (\cdot) = 1$ when~$ a> k$).
\end{lemma}

\subsection{Passenger split across platforms: Wardrop Equilibrium}
\label{sec:WE_QoS}

We model the split of the aggregate passenger arrival rate~$\Lambda$ into the system across the two platforms as a Wardrop equilibrium (WE) based on the Quality of Service (QoS); recall that platform~$i$ sees passenger arrivals as per a Poisson process with rate~$\lambda_i$. 
\rev{However, before describing our QoS metric formally, we first define the Wardrop split of the passenger arrival rate in terms of a generic QoS metric~$Q.$}

Let~$Q_i(\lambda)$ denote the QoS of platform~$i$ when passengers arrive at rate~$\lambda.$ Note that~$Q_i$ will, in general, also depend on the pricing policy~\rev{$\phi_i$}
employed by
platform~$i,$ though this dependence is suppressed for simplicity. The
Wardrop split~$(\lambda_1,\lambda_2)$ under pricing policy~$\VPhi = (\phi_1, \phi_2),$ and abandonment rate \rev{$\beta$} is defined as (see~\cite{WE}):
\begin{equation}\label{eqn_existence_uniq_WE_main}
\lamo \in \underset{\lambda \in [0,\Lambda] }{\arg\min} \left(Q_1  (\lambda) - Q_2 (\Lambda-\lambda) \right )^2, \quad \lamt = \Lambda - \lamo. 
\end{equation}
We address the uniqueness of the above Wardrop split in Lemma~\ref{lem_exist_unique_WE} stated at the end of this section. \rev{Note that under the Wardrop equilibrium, QoS is as \textit{balanced} as possible across platforms, and any infinitesimal passenger has no incentive to switch between platforms. Indeed, the Wardrop equilibrium may be interpreted as a Nash equilibrium among non-atomic passengers; see~\cite{WE}. Having defined the Wardrop
split in generic terms, we now describe the primary QoS metric considered in this paper.} 

Note that passengers can leave the system without taking a ride either because of driver unavailability, or because the price quoted was too high. Let~$\DA_i$ be the long-run fraction of passengers who leave the system on account of driver \textit{unavailability}. By the well-known PASTA property, this long-run fraction equals the stationary probability of zero waiting drivers on platform~$i$.
Thus from Lemma~\ref{lem_BCMP_Product_form}, \begin{equation}
\label{eqn_pi0_expression}
\rev{\DA_i = \DA_i (\phi_i, \lambda_i; \beta)} = \sum_{\nir = 0}^{\infty} \pi_i((0, \nir)) = \left(\sum_{n=0}^{\infty} \left[\frac{e^{n}}{\prod_{a=1}^{n}(\lambda_i f(\phi_i) + a \beta)}\right]\right)^{-1}.
\end{equation}
%
Let~$\PB_i$ denote the long-run fraction of passengers that get \textit{blocked} on platform~$i,$ i.e., those that leave platform~$i$ without taking a ride (due to driver unavailability \textit{or} a high price). Using PASTA again, 
\begin{align}\label{eqn_expression_of_PBi}
\rev{\PB_i  = \PB_i (\phi_i, \lambda_i; \beta) }&= \DA_i + \sum_{\niw=1}^{\infty} \sum_{\nir=0}^{\infty} (1 - f(\phi_i)) \ \pi_i(\niw, \nir) . 
\end{align}

\rev{Our main results assume $Q_i = \PB_i,$ i.e., the passenger Wardrop split seeks to equalize the stationary probability of passenger blocking across platforms (we emphasize once again that this blocking captures both driver unavailability as well as price sensitivity).
This analysis is presented in Section~\ref{sec:overall_bp}. A related QoS metric based on expected pick-up delay is subsequently analyzed in Section~\ref{sec:delay}.}

\rev{Finally, we establish existence and uniqueness of the Wardrop equilibrium under any QoS metric that satisfies a certain monotonicity assumption {\bf A}.3 (proof in Appendix~\ref{sec:appendix_A_model}); we also show that $\PB_i$ satisfies {\bf A}.3.}

\rev{
\begin{enumerate}[{\bf A}.3]
\item The QoS function~$Q_i$ for each~$i$ is continuous and strictly increasing in~$\lambda_i$.
\end{enumerate}}

 \rev{
\begin{lemma}
\label{lem_exist_unique_WE}
Consider any QoS metric satisfying Assumption~{\bf A}.3. Given any
price policy~$\VPhi = (\phi_1, \phi_2)$ and $\beta$, there exists a unique Wardrop
Equilibrium~$(\lamo,\lamt)$ solving \eqref{eqn_existence_uniq_WE_main}. Moreover, the QoS metric~$\{\PB_i\}$ satisfies {\bf A}.3 for~$\beta > 0$.
\end{lemma}
}
\ignore{
 Passengers either accept or reject service based on price. When they accept, the time before a driver arrives and the journey starts depends upon the number of drivers, $exp(n \nu)$ if $n$ drivers. If they don't accept, then $\alpha$ is the expected time to derive another service elsewhere.
 }
\rev{The above lemma establishes the existence and uniqueness of the WE for each $(\VPhi, \beta)$ with $\beta > 0,$ under QoS metric $Q_i = \PB_i,$ defined in~\eqref{eqn_expression_of_PBi}. 
}

\subsection{Non-cooperative (pricing) game between platforms}
\label{ssec:strategic_game}

We treat the action of each platform (say~$i$) to be its pricing policy~\rev{$\phi_i \in [0,\phi_h].$} We define the payoff of a platform as the (almost sure) rate at which it derives revenue from matching drivers with passengers, denoted by~$\MR_i.$ Note that~$\MR_i$ depends on the price (action) profile~$\VPhi$ and other parameters; the pricing policy of each platform influences the other's utility via the Wardrop split~\eqref{eqn_existence_uniq_WE_main} that determines the market shares of both platforms.

\rev{More formally, note that a platform~$i$ derives revenue $\phi_i$ every time it makes a successful match, i.e., when a passenger arrives, finds a waiting driver, and finds the offered price acceptable. 
Let $\R_i(t)$ represent the (random) cumulative matching revenue accumulated by platform~$i$ until time $t.$ Then the matching revenue rate is defined as $\MR_i = \lim_{t \to \infty} \nicefrac{\R_i(t)}{t}$ almost surely (if such a limit exists). 
}
In the following, we show that the above limit exists and obtain an expression of the  matching revenue rate (often referred to simply as the \textit{revenue rate}) of each platform in terms of the Wardrop split~$(\lambda_1(\VPhi;\beta),\lambda_2(\VPhi;\beta))$ (proof in Appendix~\ref{sec:appendix_A_model}).
\begin{lemma} \label{lem_mr_derivation}
The matching revenue rate of platform~$i$ 
is given by,
\begin{equation} \label{eqn_MR_expression} \MR_i(\VPhi,\lambda_i(\VPhi;\beta);\beta)  = \hspace{-2mm} \sum_{s_i \colon \niw \neq 0}  \hspace{-1mm} \lambda_i(\VPhi;\beta) \rev{f(\phi_i) \phi_i} \pi_i(s_i) = \rev{\lambda_i(\VPhi;\beta) f(\phi_i) \phi_i  (1- \DA_i (\phi_i, \lambda_i; \beta))}.
\end{equation}
\end{lemma}

\rev{Having defined the payoff function of each platform as a function of the actions (specifically, prices $(\phi_1,\phi_2)$, the remainder of this paper is devoted to analysing this non-cooperative game between the platforms. Since we assume symmetric platforms, we analyse only symmetric equilibria in this study.} As we will see, these equilibria will be parameterized by the \textit{driver-passenger ratio} (\dpr), defined as~$\rho := \nicefrac{e}{\Lambda};$ note that the \dpr\ captures the relative abundance of drivers and passengers in the system.

\subsection{Infinite Driver Patience (\idp) regime}


\rev{The non-cooperative game defined above is intractable even under the assumption of symmetric platforms; this is primarily because the WE, which dictates the arrival rate seen by each platform, does not admit an explicit closed form. However, we are able to analyse the game in a certain scaling regime, which is defined next.}

Specifically, we approximate the payoff functions of the platforms by
letting~$\beta \downarrow 0.$ This scaling regime, in which the patience times of waiting drivers are scaled to infinity, is referred to as the \textit{Infinite Driver Patience (\idp) regime}. While the \idp\ regime enables an
explicit characterization of the equilibria of the non-cooperative
game under consideration, it is also well motivated. Indeed,
$\beta \approx 0$ means that driver abandonment times are
stochastically much larger than the ride durations and the passenger
interarrival times. This is reasonable in many practical scenarios,
particularly when drivers are `tied' to their respective platforms,
and are available to take on rides for long durations (say an 8-10 hour work
shift) at a stretch.

We conclude this section with a result characterizing the limit of the
revenue rate as $\beta \to 0$ when the corresponding WE converge (proof in
Appendix~\ref{sec:appendix_A_model})\footnote{Throughout our notations, we emphasize
functional dependence on parameters of interest only as and when
required.}.

\rev{
\begin{lemma} 
\label{lem_approx_MR}
Suppose that for price policy~$\VPhi$, the passenger arrival rates (at WE)~$\lambda_i (\VPhi;\beta) $ converge as
$\beta \to 0$, with the limit represented by  ${\lambda}_i(\VPhi;0)$. Then the revenue rate (payoff) of platform~$i,$ $\MR_i
(\VPhi,\lambda_i(\VPhi;\beta); \beta),$ converges as~$\beta \to 0$.  Further, the limit  denoted by $\MR_i(\VPhi,\lambda_i(\VPhi;0); 0)$ equals: 
\begin{align}
\label{eq_approx_MR}
\MR_i(\VPhi,\lambda_i(\VPhi;0); 0) & =  \begin{cases}
        \eiw \phi_i & \text{ if }  \left(\frac{\eiw}{\lambda_i(\VPhi;0) f(\phi_i)}\right) < 1, \\
        \lambda_i(\VPhi;0) f(\phi_i) \phi_i & \text{ else.}
    \end{cases} 
\end{align}
\ignore{Along similar lines, we have
$\DA_i(\VPhi, \lambda_i(\VPhi,\beta); \beta) \to \DA_i (\VPhi, \lambda_i(\VPhi;0), 0)$,
where
\begin{align}
\label{eq_pi0_approx_lem}
\DA_i(\VPhi, \lambda_i(\VPhi;0), 0) &=
        \max \left\{1-\frac{\eiw}{\lambda_i(\VPhi;0) f(\phi_i)}, 0\right\}.
\end{align}
}
\end{lemma}
}

By virtue of the above result, we derive approximate equilibria of the actual system when~$\beta$ is small, by analysing the equlibria in the \idp\ regime, which is obtained by letting~$\beta \to 0.$ 
We also show that the equilibrium behavior in the \idp\ regime (which can be thought of as corresponding to~$\beta = 0$) is also meaningful in the more practical `pre-limit' where $\beta$ is small and positive. Formally, we show that the equilibria under the \idp\ regime are~$\epsilon$-equilibria when $\beta$ is small.

\ignore{
I understood the point, but the flow is back and forth and I feel the \idp\ regime is not defined properly.

\begin{itemize}
    \item  Can we cleanly define the \idp\ regime. Motivate why we study this (as is currently done in the paragraphs before the Lemma)

    \item We can then talk about the system with driver impatience and discuss its analysis via the \idp\ regime analysis. 

    \item should we name \idp\ regime with a more meaningful name, like 'permanent/loyal driver system'

\item  Loyal Drivers:  The drivers wait infinite amount of time and leave only after they get some count  of rides  (in our current model, which we call \idp\ regime) this count is a geometric random variable). 
    
\end{itemize}
}
\section{Monopoly}
\label{sec:monopoly}

In this section, we consider a monopolistic scenario, where a single platform optimizes its pricing in the absence of competition. This analysis serves as a benchmark for the game theoretic analysis in the following sections (where the competition between platforms is captured explicitly), shedding light on the impact of inter-platform competition on pricing and passenger/platform utility.

Specifically, we analyse the pricing of a `monopolistic' platform that sees a passenger arrival rate of~$\nicefrac{\Lambda}{2},$ irrespective of its pricing strategy. (Other model details, including  passenger price sensitivity and driver behavior, are as described in Section~\ref{sec:model}.) Note that this captures the behavior of each platform in our model from Section~\ref{sec:model}, if the passenger arrival rate seen by each platform is \textit{exogenously} set to half the total arrival rate (as opposed to being set \textit{endogenously} via the QoS based Wardrop split), to simulate the absence of competition.

In the monopolistic scenario, note that if the platform offers a price~$\phi$ to an incoming passenger, then the passenger accepts this price with probability~$f(\phi)$ (assuming a waiting driver is available). The revenue rate (utility) derived by the platform is
then given by Lemma~\ref{lem_mr_derivation} with~$\lambda = \nicefrac{\Lambda}{2}$ (see
\eqref{eqn_MR_expression}):
\begin{equation}
\label{eqn_mr_monopoly}
\MR(\phi; \beta) = \frac{\Lambda}{2} f(\phi) \phi (1 - \DA(\phi; \beta)). 
\end{equation}
This revenue rate can be approximated using Lemma~\ref{lem_approx_MR} when~$\beta$ is close to zero. Formally, we define the utility of the platform in the \text{\idp\ regime} ($\beta = 0$) as the pointwise limit of the utility as~$\beta \to 0,$ which under the monopoly setting equals, 
\begin{equation}
  \label{eq:MR_monopoly}
    \MR(\phi; 0) 
    = \begin{cases}
        e \phi & \text{ if }  \phi < f^{-1}\left(2\rho\right) , \\
        \frac{\Lambda}{2} f(\phi) \phi & \text{ else.}
    \end{cases}
\end{equation} 
The optimal pricing strategy for the monopolistic platform that seeks to maximize~$\MR(\phi; 0)$ is characterized as follows (proof in Appendix~\ref{sec:appendix_mono}). (Recall that the driver-passenger ratio (\dpr) is defined as~$\rho=\nicefrac{e}{\Lambda}.$)

\begin{thm} \label{thm_monopoly_opt}
Consider the \idp\ regime~$(\beta = 0)$. The optimal monopoly price is a non-increasing function of the~\dpr~$\rho.$
Specifically, define 
\begin{align*}
\dm(\phi):= f(\phi) + \phi f'(\phi),\ \phi_m := \max \{ \phi \in [0, \phi_h] : \dm(\phi) \geq 0\}.
\end{align*}
\begin{enumerate}[$(1)$]
    \item If~$\rho \leq f(\phi_m)/2$, then~$ f^{-1}\left(2\rho\right)$ is the unique optimizer;
    \item If~$\rho > f(\phi_m)/2$, then~$\phi_m$ is the unique optimizer. 
\end{enumerate}
Moreover, for any sequence of optimal prices corresponding to a sequence $\beta_n \to 0$, there exists a sub-sequence that converges to the unique optimal price of the \idp\ regime. 
\end{thm}

Theorem~\ref{thm_monopoly_opt} presents a complete characterization of the optimal pricing policy for a monopolistic platform in the \idp\ regime, while also establishing the validity of this approximation when~$\beta$ is positive and small, i.e., under the assumption of highly patient drivers. The optimal price is solely determined by the price sensitivity function~$f$ and the~\dpr~$\rho$.
Indeed, the optimal price quoted by the platform decreases (formally non-increases) as a function of~$\rho$, indicating that as the rate of driver arrivals increases (or equivalently, as the rate at which passengers enter the system decreases), the optimal price decreases, to induce a higher rate of matched rides. Formally, the optimal pricing may be interpreted as follows.
\begin{itemize}
\item It is instructive to first consider the special case~$\dm(\phi_h) \geq 0;$
this corresponds to the scenario where passengers are relatively price-insensitive.\footnote{Consider two price sensitivity functions~$f_1$ and $f_2$ satisfying Assumption~{\bf A}.1-2, such that~$f_2'(\cdot) \leq f_1'(\cdot) \leq 0.$ Note that $f_2$ corresponds to a more price-sensitive passenger base compared to $f_1.$ It is then easy to see that $\dm_1(\cdot) \geq \dm_2(\cdot),$ where $\dm_1(\phi) := f_1(\phi) + \phi f_1'(\phi)$ and $\dm_2(\phi) := f_2(\phi) + \phi f_2'(\phi).$ In this sense, the condition~$\dm(\phi_h) \geq 0$ captures the scenario where the passenger base is relatively price-insensitive.}
In this case, $\phi_m = \phi_h$ (as~$\dm$ is decreasing), which also implies  $f^{-1}(2\rho) = \phi_h$ (by definition) and hence  $\phi_h$ is the unique optimizer  under both parts~(1)-(2). 
Thus it is optimal for the platform to exploit the passenger price-insensitivity and set the maximum permissible price~$\phi_h,$ irrespective of the~\dpr.
\item Next, consider the case~$\dm(\phi_h) < 0$ and~$\rho \leq f(\phi_m)/2.$ This corresponds to the scenario where passengers are relatively price-sensitive, and the \dpr\ is small (i.e., drivers are scarce relative to passengers). In this case, the optimal price equals~$f^{-1}(2\rho),$ which is a decreasing function of \dpr, as expected---the scarcer the drivers, the higher the price. 
\item Finally, consider the case~$\dm(\phi_h) < 0$ and~$\rho > f(\phi_m)/2.$ This corresponds to the scenario where passengers are relatively price-sensitive, and the \dpr\ is large (i.e., passengers are scarce relative to drivers). This leads to a lower optimal price~$\phi_m,$ which does not vary further with the \dpr~$\rho$ (it is the maximizer of the function~$f(\phi)\phi$); basically the system is saturated with an abundance of drivers, and $\nicefrac{\Lambda}{2} f(\phi)\phi $ represents the revenue rate of the system, which is maximized at $\phi_m$. 
\end{itemize}
These observations are further corroborated by Figure~\ref{fig:Square} in Section~\ref{sec:comparison}.

In subsequent sections, we will draw insights by contrasting the monopolistic pricing strategy characterised here with the equilibrium strategy arising from the competition between the platforms. As noted before, this comparative analysis will shed light on how passenger-side churn affects the pricing policies of platforms.
\section{Duopoly  driven by price \& driver unavailability}
\label{sec:overall_bp}

In this section, we analyse the competition between platforms under a passenger QoS metric that captures the overall likelihood of finding a ride on arrival (note that a ride may not materialise either due to driver unavailability, or because the price is unacceptable). Formally, in the language of Section~\ref{sec:model}, this corresponds to~$Q_i = \PB_i$ for all~$i,$ where (see~\eqref{eqn_expression_of_PBi}) 
\begin{align}
\label{Eq_overall_exact}
    \PB_i(\VPhi; \beta) &= 
     \DA_i(\VPhi; \beta)  + (1-f(\phi_i))(1-\DA_i(\VPhi; \beta))   =  \DA_i(\VPhi; \beta) f(\phi_i) +(1-f(\phi_i)).
\end{align}
The goal of this section is to analyse and interpret the equilibria associated with the corresponding non-cooperative game between the platforms. Our main insight is that when passengers are scarce (relative to drivers), platforms are forced to deviate from monopoly pricing in order to `fight' for market share. 
In particular, we find that under certain conditions, no pure NE exists; instead, we demonstrate a mixed NE, and also a novel equilibrium concept having dynamic connotations, that we call an \textit{equilibrium cycle}.

We begin by defining the \idp\ regime formally. Recall that the payoff functions in the \idp\ regime are pointwise limits of our payoff functions as~$\beta \to 0.$ In light of Lemma~\ref{lem_approx_MR}, the limiting payoff functions can be deduced from the limit of the WE as~$\beta \to 0$ (with some abuse of notation, we refer the latter limit as the \textit{WE in the \idp~regime}). The following theorem characterizes these limits, and thereby establishes continuous extensions of the WE as well as the payoff functions over~$\beta \geq 0.$

\ignore{
use this limit to define the quantities at \idp~regime and subsequently discuss the continuity of the matching revenue with respect to $\beta$ including the \idp~regime.

Note that the passenger blocking probability on platform~$i$ when it uses static price policy~$\phi_i$ and when passengers arrive at rate~$\lambda_i$ is given by (see~\eqref{eqn_expression_of_PBi}),
\begin{align}
\label{Eq_overall_exact}
    \PB_i(\VPhi; \beta) &= 
     \DA_i(\VPhi; \beta)  + (1-f(\phi_i))(1-\DA_i(\VPhi; \beta))   =  \DA_i(\VPhi; \beta) f(\phi_i) +(1-f(\phi_i)),
\end{align}
\noindent the limiting value of which, as~$\beta \to 0,$ is given by (using Lemma~\ref{lem_approx_MR})
\begin{equation}
    \label{Eq_overall_aproxx}
   \PB_i(\VPhi; \lambda_i(\VPhi, 0))  = 
        \begin{cases}
            1 - \frac{\eiw}{\lambda_i} & \text{ if }  \left(\frac{\eiw}{\lambda_i f(\phi_i)}\right) < 1, \\
            1 - f(\phi_i) & \text{ else.}
        \end{cases}
\end{equation}
As before, the revenue rate of platform~$i$, when the two platforms operate with the static price vector~$\VPhi,$ equals~$\MR_i (\VPhi) = \MR_i (\VPhi; \lami)$, where the WE splits~$\lami$ are defined using the QoS metric~$\PB$.
In the following, \kav{we first characterize the limit of the WE as~$\beta \to 0$ (with slight abuse of notation, we refer the limit as WE at \idp~regime), use this limit to define the quantities at \idp~regime and subsequently discuss the continuity of the matching revenue with respect to $\beta$ including the \idp~regime.}
} 

\begin{thm}
\label{thm_WE_MR_PB}
Fix price vector~$\VPhi$. Extend the definition of WE and revenue rate~$\MR$ for $\beta = 0$ as in Table~\ref{tab:WE_MR_PB}.
Then for any~$i \in \{1,2\}$, the
two mappings~$\beta \mapsto \lambda_i(\VPhi; \beta)$ 
 and~$\beta \mapsto \MR_i(\VPhi;\beta)$ are continuous over~$[0,\infty)$.
\end{thm}

\begin{table}[ht]
\renewcommand{\arraystretch}{1.4}
\centering
\begin{tabular}{|c|ccc|cc|}
\hline
    \multirow{2}{*}{}    & \multicolumn{3}{c|}{Range of~$\phi_1$ when~$\phi_1 > \phi_2$}  & \multicolumn{2}{c|}{Range of~$\phi_1$ when~$\phi_1 = \phi_2$} \\
    \hhline{~-----} 
    & \multicolumn{1}{c|}{$ [0, f^{-1}(2\rho))$} & \multicolumn{1}{c|}{$ [f^{-1}(2\rho), f^{-1}(\rho))$} & $ [f^{-1}(\rho), \phi_h]$ & \multicolumn{1}{c|}{$ [0, f^{-1}(2\rho))$} & $ [f^{-1}(2\rho), \phi_h]$   \\ \hline
$\lambda_1(\VPhi; 0)$ & \multicolumn{1}{c|}{$\nicefrac{\Lambda}{2}$}         & \multicolumn{1}{c|}{$\Lambda - \nicefrac{e}{f(\phi_1)}$}        & 0                                   & \multicolumn{1}{c|}{$\nicefrac{\Lambda}{2}$}         & $\nicefrac{\Lambda}{2}$                \\ \hline
$\MR_1(\VPhi; 0)$     & \multicolumn{1}{c|}{$e\phi_1$}                       & \multicolumn{1}{c|}{$m(\phi_1)$}                                & 0                                   & \multicolumn{1}{c|}{$e\phi_1$}                       & $\nicefrac{\Lambda}{2}f(\phi_1)\phi_1$ \\ \hline
$\MR_2(\VPhi; 0)$     & \multicolumn{1}{c|}{$e\phi_2$}                       & \multicolumn{1}{c|}{$e\phi_2$}                                  &~$\Lambda f(\phi_2)\phi_2$           & \multicolumn{1}{c|}{$e\phi_2$}                       & $\nicefrac{\Lambda}{2}f(\phi_2)\phi_2$ \\ \hline
\end{tabular}
\caption{\idp\ regime under QoS~$\PB$. Here,~$m(\phi) = (\Lambda f(\phi)-\eow)\phi$,~$\lambda_2(\VPhi; 0) = \Lambda - \lambda_1(\VPhi; 0).$ When~$\phi_1 < \phi_2,$ reverse the role of~$\phi_1$ and~$\phi_2.$
\label{tab:WE_MR_PB}
}
\end{table}

In essence, Theorem~\ref{thm_WE_MR_PB} (proof in Appendix~\ref{sec:appendix_duopoly_PB}) defines the \idp\ regime ($\beta = 0$), by specifying the payoff of each player as a function of the price (action) profile.
In the remainder of this section, we analyse the equilibria in the \idp\ regime, and then show that these are also meaningful in the `pre-limit' (i.e., when~$\beta$ is small).
\subsection{Equilibria in the \idp\ regime}

The following theorem characterizes the Nash equilibria of the \idp\ regime (proof in Appendix~\ref{sec:appendix_duopoly_PB}).

\begin{thm}
\label{thm_B_sys_limit_sys_main}
Consider the \idp\ regime~$(\beta = 0)$. Define
\begin{align*}
\db(\phi):= f(\phi) + 2\phi f'(\phi), \ \phi_b = \max \{ \phi \in [0, \phi_h] : \db(\phi) \geq 0 \}.  
\end{align*}

\begin{enumerate}[$(1)$]
    \item If~$\rho \leq f(\phi_b)/2,$ then~$\left( f^{-1}\left(2\rho \right), f^{-1}\left(2 \rho \right)\right)~$ is the unique symmetric NE;
    \item If~$f(\phi_b)/2 < \rho < 1,$ then there does not exist a symmetric pure strategy NE. However, there exists a symmetric mixed strategy NE  ~$(\sigma^*,\sigma^*),$ supported over~$[\eL,\eU]$, where~$\eL, \eU$ are respectively the maximum value and the maximizer of the function $\phi \mapsto (f(\phi) - \rho) \phi /\rho$, and
   ~$$\sigma^*([\eL, \phi]) =  \frac{e(\phi + \eU) - \Lambda f(\eU)\eU}{ 2e\phi - \Lambda f(\phi)\phi} \text{ for } \phi \in [\eL, \eU];$$ 
    \item If~$\rho \ge 1$ then for any~$\epsilon > 0$, choose~$0 < \delta \leq \phi_h~$ such that  ~$\sup_{\phi \leq \delta} \Lambda f(\phi)\phi < \epsilon$. Then~$(\delta, \delta)$ is an~$\epsilon$-NE.
\end{enumerate}
\end{thm}

Theorem~\ref{thm_B_sys_limit_sys_main} characterizes the impact of competition between platforms on their equilibrium pricing strategies (with market shares determined by the QoS metric~$\PB_i$). Crucially, the nature of the equilibrium itself depends on the \dpr, which captures the relative abundance of drivers and passengers---a pure strategy NE arises over a certain range, a mixed strategy NE arises over another, and finally, an~$\epsilon$-NE arises the \dpr\ exceeds~1. Similar to the monopoly pricing analysed in Section~\ref{sec:monopoly}, the equilibrium price tends to decrease as a function of the \dpr~$\rho$. In other words, as the rate of driver arrivals increases (or conversely, as the rate at which passengers enter the system decreases), the competition between platforms intensifies, resulting in a lower equilibrium price.

\begin{figure}[ht]
\centering
\includegraphics[trim = {1.5cm 8cm 2.5cm 8cm}, clip, scale = 0.4]{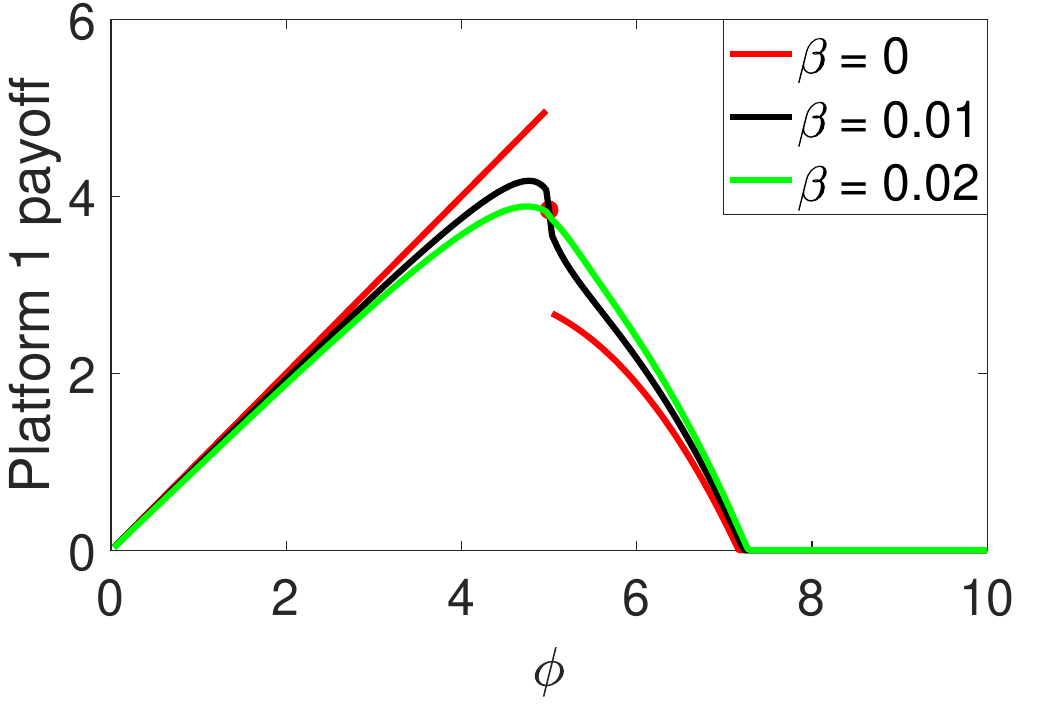}
\caption{Depiction of discontinuity of~$\MR_1(\phi,\phi_2)$ v/s~$\phi$ (fixing~$\phi_2 = 5)$ in the \idp\ regime, with~$\Lambda = 2, e = 1, \phi_h = 10, f(\phi) = 1 - (a\phi)^2$, where~$a = 0.1$.}
\label{fig:MR_discontinuity}
\end{figure}

Next, we formally interpret the equilibrium prices as established in Theorem~\ref{thm_B_sys_limit_sys_main}. Once again, it is important to emphasize that although there are clear structural similarities between the statements of Theorems~\ref{thm_monopoly_opt} and~\ref{thm_B_sys_limit_sys_main}, the quantities~$\db$ and~$\phi_b$ in the latter differ from their counterparts in the previous theorem.
\begin{itemize}
    \item 
    We begin by considering the case~$\rho \leq f(\phi_m)/2;$ this corresponds to the scenario where drivers are scarce relative to passengers. Under these favourable market conditions, the equilibrium is characterized by both platforms operating at the optimal monopoly price~$f^{-1}\left(2\rho \right)$. 
    \item Next, consider the case~$f(\phi_m)/2 < \rho \leq f(\phi_b)/2;$ this corresponds to the scenario where the \dpr~$\rho$ is somewhat larger (i.e., drivers are less scarce relative to passengers). In this case, competition between the platforms induces a deviation from the optimal monopoly price; each platform sets an equilibrium price~$f^{-1}\left(2\rho \right)$ that is strictly smaller than the optimal monopoly price~$\phi_m$.
    \item Next, consider the case~$f(\phi_b)/2 < \rho < 1;$ this corresponds to the scenario where the~\dpr~$\rho$ is quite large, indicating a scarcity of passengers (relative to drivers). In this case,  competition between the platforms intensifies, leading to a further deviation of the equilibrium behavior from the monopoly setting. In particular, a pure strategy NE does not even exist in this case. Formally, this is because of certain discontinuities in the payoff functions, which are in turn induced by discontinuities in the Wardrop split. Indeed, when passengers are scarce, the Wardrop split is primarily dictated by price; as a consequence, a slight reduction in price by one platform (below the price posted by the other) can result in a significant jump in its market share. Figure~\ref{fig:MR_discontinuity} provides an illustration of this payoff function discontinuity in the \idp\ regime.\footnote{Such discontinuities do not arise when~$\beta > 0$.}
    
    While a pure NE does not exist if~$f(\phi_b)/2 < \rho < 1$ owing to the above mentioned payoff discontinuities, we show that there does exist a symmetric \textit{mixed} NE in this case, characterized by a continuous distribution~$\sigma^*$ supported over the closed interval~$[\eL, \eU]$. We note that both~$\eL$ and~$\eU$ are decreasing (formally, non-increasing) functions of~$\rho$. While the mixed NE allows one to interpret the equilibrium behavior of the platforms via \textit{randomised} prices, a \textit{dynamic} interpretation in terms of oscillating prices is also possible. We formalise this interpretation below by introducing a novel equilibrium called \text{equilibrium cycle} (EC), and show that in this case,~$[\eL, \eU]$ is also an EC. \rev{Interestingly, the absence of pure Nash pricing, and the occurrence of randomised/oscillatory prices arises only when the \dpr\ is large, i.e., when passengers are scarce relative to drivers.}
    
     \item Finally, consider the case~$\rho \geq 1;$ this corresponds to the scenario where passengers are very scarce, driving both platforms to operate at near zero prices. Formally however, it can be shown that~$(0,0)$ is not a NE. Instead, we show that~$(\delta,\delta)$ is an~$\epsilon$-NE for suitably small values of~$\delta$ and~$\epsilon$, suggesting that players choose a price `close to zero' in such a competitive environment.

\end{itemize}

These observations are further corroborated by Figure~\ref{fig:Square} in Section~\ref{sec:comparison}. In Section~\ref{sec:comparison}, we also conduct a comprehensive comparison of the equilibrium prices derived from Theorem~\ref{thm_B_sys_limit_sys_main} and the optimal monopolistic price as characterized in Section~\ref{sec:model}. This analysis sheds light on the welfare implications of platform competition in different market conditions.

We now revisit the case~$f(\phi_b)/2 < \rho < 1,$ with the aim of developing a dynamic interpretation of the strategic interaction between platforms. It is instructive to consider the incentives for unilateral deviation in price by each platform in this case. Specifically, it turns out that when the opponent, say platform~$-i,$ offers a price within the interval~$[\eL, \eU],$ platform~$i$ is incentivised to also set its price within the same interval. In particular, if the opponent's price lies within the range~$(\eL, \eU],$ platform~$i$ has an incentive (though not a best response, which does not exist) to set its price slightly below the opponent. 
On the other hand, if the opponent's price is~$\eL,$ then platform~$i$'s best response is to set its price at~$\eU.$ In essence, this suggests an indefinite oscillation of prices within the interval~$[\eL, \eU].$ This motivates the following notion of an equilibrium cycle.

\begin{defn}[Equilibrium Cycle]
\label{def_equi_cycle}
A closed interval~$[a, b]$ is called an equilibrium cycle (EC) if:
\begin{enumerate}[$(i)$]
    \item for any~$ i \in \mathcal{N}$, for~$\phi_{-i} \in [a,b],$ there exists~$\phi_i \in [a,b]$ such that,~$$\mathcal{M}_i(\phi_i,\phi_{-i}) > \mathcal{M}_i(\tilde{\phi}_i,\phi_{-i}) \ \forall \ \tilde{\phi}_i \in [0,\phi_h] \setminus [a,b],$$
    \item for any price vector~$\VPhi \in [a,b]^2$, there exists~$ i \in \mathcal{N}$ and~$ \phi_i' \in [a,b]$ such that~$\MR_i(\phi_i', \phi_{-i}) > \MR_i(\VPhi)$ and~$$\mathcal{M}_i(\phi_i', \phi_{-i}) > \mathcal{M}_i(\tilde{\phi}_i, \phi_{-i}) \ \forall \ \tilde{\phi}_i \in [0,\phi_h] \setminus [a,b],$$
    \item no subset of~$[a,b]$ satisfies the above two conditions.
\end{enumerate}
\end{defn}

The first condition above establishes the `stability' of the interval~$[a,b];$ if any player plays an action in this interval, the other player is also incentivized to play an action in the same interval (this choice dominating any action outside the interval). The second condition establishes the `cyclicity' of the same interval; if both players play any actions within the interval~$[a,b]$, at least one player has the incentive to deviate to a different action within the same interval (this deviation also dominating any action outside the interval). The last condition ensures the minimality of the interval, i.e., no strict subset of an EC is an EC.\footnote{To the best of our knowledge, EC is a novel equilibrium concept in game theory. Formalizing the conditions under which it arises in a general class of non-cooperative games represents a promising avenue for future research.}

\begin{thm} \label{thm_EC}
Consider the \idp\ regime~$(\beta = 0)$. If~$f(\phi_b)/2 < \rho < 1,$ then the interval~$[\eL, \eU]$ is an equilibrium cycle.
\end{thm}
Theorem~\ref{thm_EC} (proof in Appendix~\ref{sec:appendix_duopoly_PB}) establishes that if~$f(\phi_b)/2 < \rho < 1,$ the support of the symmetric mixed NE as established in Theorem~\ref{thm_B_sys_limit_sys_main} is also an equilibrium cycle. Intuitively, this EC can be interpreted (in dynamic terms) as the limit set of alternating `better' response dynamics between the platforms (not that \textit{best} response may not exist in this case owing to the payoff function discontinuities discussed before). We provide numerical illustrations of these oscillating dynamics in Figure~\ref{fig:BR_dynamics_EC}.

Finally, we comment of the payoff of each platform in the case $f(\phi_b)/2 < \rho < 1.$
\begin{lemma} \label{lem_mixed_NE_payoff_and_sec_value}
If~$f(\phi_b)/2 < \rho < 1,$ then the average payoff of each platform~$i$ under the mixed NE~$(\sigma^*, \sigma^*)$ 
is given by~$\MR_i(\sigma^*,\sigma^*; 0) = m(\eU).$ Moreover, the security value and strategy of each platform are~$m(\eU)$ and~$\eU,$ respectively.\footnote{ \label{footnote_Security_value_def} Given a strategic form game,~$G = \left< \mathcal{N}, (\Phi_i), (\MR_i)\right>$, the security value of a player~$i$ is given by:
$$
\underline{\MR_i} = \max_{\phi_i \in \Phi_i}\min_{\phi_{-i} \in \Phi_{-i}} \MR_i(\phi_i, \phi_{-i}).
$$
Any strategy~$\phi_i^s \in \Phi_i$ that guarantees this payoff to player~$i$ is called a security strategy of player~$i$; see \cite[Definition 6.3]{narahari}.}
\end{lemma}
Interestingly, the average payoff of each platform under the mixed NE agrees with its security value. Moreover, this is also the payoff that each platform receives when one platform plays the price~$\eL$ (the left endpoint of the EC) and the other plays the price~$\eU$ (the right endpoint of the EC), which is its best response to the opponent's action, and also the security strategy.

\subsection{Applicability of \idp\ regime equilibria in the pre-limit (for small~$\beta > 0$)}
We now establish a connection between the equilibria in the \idp\ regime (Theorems~\ref{thm_B_sys_limit_sys_main} and~\ref{thm_EC}) and equilibrium behavior in the pre-limit, where~$\beta$ is small and positive. Specifically, we show that the equilibria of the \idp\ regime are~$\epsilon$-equilibria in the pre-limit. We recall here the definition of an $\epsilon$-NE (see~\cite[Section 3.13]{myerson}).

\begin{defn}[$\epsilon$-NE] \label{def_epsilon_equilibria}
Consider a strategic form game~$G = \left< \mathcal{N}, (\Phi_i), (\MR_i)\right>$ where~$\mathcal{N}$ is a finite set, the~$\Phi_i$ are nonempty compact metric spaces, and the~$\MR_i : \Phi_i \to \mathbb{R}$ are utility functions. Given~$\epsilon \geq 0,$ a mixed strategy~$\sigma$ is called an~$\epsilon$-NE if for all~$i \in \mathcal{N}$ and~$\phi_i \in \Phi_i$,
$$
\MR_i(\phi_i, \sigma_{-i}) \leq \MR_i(\sigma_i, \sigma_{-i}) + \epsilon.
$$
Note that when~$\epsilon = 0$, an~$\epsilon$-NE is just an NE in the usual sense.
\end{defn}
The notion of~$\epsilon$-NE relaxes the notion of the NE in that no player stands to make a utility gain exceeding~$\epsilon$ from a unilateral deviation. Next, we define an~$\epsilon$-EC as follows (note that this definition, like the definition of the EC itself, is novel):
\begin{defn} [$\epsilon$-equilibrium cycle]
\label{def_epsilon_equi_cycle}
 A closed interval~$[a ,b]$ is called an~$\epsilon$-equilibrium cycle ($\epsilon$-EC) if:
\begin{enumerate}[$(i)$]
\item  for any~$ i \in \mathcal{N}$, for~$\phi_{-i} \in [a + \epsilon, b - \epsilon],$  there exists~$\phi_i \in [a,b]$ such that, 
$$\mathcal{M}_i(\phi_i,\phi_{-i}) > \mathcal{M}_i(\tilde{\phi}_i,\phi_{-i}) \ \forall \ \tilde{\phi}_i \in [0,\phi_h] \setminus [a - \epsilon, b + \epsilon],$$
\item for any price vector~$\VPhi \in [a + \epsilon, b - \epsilon]^2$, there exists~$ i \in \mathcal{N}$ and~$ \phi_i' \in [a,b]$ such that~$\MR_i(\phi_i', \phi_{-i}) > \MR_i(\VPhi)~$ and~$$\mathcal{M}_i(\phi_i', \phi_{-i}) > \mathcal{M}_i(\tilde{\phi}_i, \phi_{-i})  \ \forall \ \tilde{\phi}_i \in [0,\phi_h] \setminus [a - \epsilon, b+\epsilon],$$

\item  there exists no subset~$[c,d]$ of~$[a + \epsilon, b - \epsilon]$ that satisfies conditions ~$(i)-(ii)$.
\end{enumerate}
\end{defn}

The definition of an~$\epsilon$-EC relaxes the notion of the EC, permitting platforms to deviate from EC strategies by a margin of~$\epsilon$. This relaxation preserves the key properties of the EC while allowing for small deviations in the space of actions.

\begin{thm}
\label{thm_B_system_epsilon_equilibrium}
Consider~$\rho < 1$. Then for any~$\epsilon > 0$, there exists~${\bar \beta}_\epsilon > 0$ such that for all $0 < \beta \le {\bar \beta}_\epsilon$, 
\begin{enumerate}[$(1)$]
    \item the NE defined in Theorem~\ref{thm_B_sys_limit_sys_main} for the \idp\ regime (both pure and mixed) are~$\epsilon$-NE, and
    \item if~$f(\phi_b)/2 < \rho < 1,$ then the interval~$[\eL, \eU]$ defined in Theorem~\ref{thm_EC} is an~$\epsilon$-EC.
\end{enumerate}
\end{thm}

Theorem~\ref{thm_B_system_epsilon_equilibrium} validates our approach of characterizing equilibria in the \idp\ regime (Theorems~\ref{thm_B_sys_limit_sys_main} and~\ref{thm_EC}), by showing that these equilibria are also meaningful in the pre-limit (i.e., when~$\beta$ is small and positive). In other words, under the practical assumption that drivers are patient, our characterization of the \idp\ regime equilibria are indicative of equilibrium behavior of the platforms.

{
\section{Duopoly driven by expected pick-up time}
\label{sec:delay}

In this section, we analyse the competition between the platforms under a different passenger QoS metric which captures the average delay (a.k.a., pick-up time) experienced by passengers. To formalise this metric, we define the `delay' in the event of a ride not materializing (either due to driver unavailability or because the price is unacceptable) as $\nicefrac{1}{\alpha},$ where $\alpha \in (0,1).$ Note that $\nicefrac{1}{\alpha}$ may also be interpreted as the dis-utility experienced by a passenger on not obtaining a ride on the platform. Under this QoS metric, we analyse the equilibria associated with the corresponding non-cooperative game between the platforms, in the regime where $\alpha, \beta \rightarrow 0,$ i.e., where the delay/disutility associated with not obtaining a ride is large, and drivers are highly patient. Our main insight is the equilibria from Section~\ref{sec:overall_bp} (under QoS metric~$Q_i=\PB_i$) remain $\epsilon$-equilibria for small~$\alpha,\beta.$ This establishes that the equilibrium pricing characterized in Theorems~\ref{thm_B_sys_limit_sys_main} and~\ref{thm_EC} are robust to pick-up delay considerations on part of passengers.\footnote{Indeed, we show in Section~\ref{sec:comparison} via numerical case studies that the aforementioned equilibria are accurate even when $\alpha$ is large.} 

Formally, we define the expected pick-up time $\E[W_i]$ on platform~$i$ as follows:
\begin{align*}
&\E\left[W_i \ |\ \text{no ride} \right] := \frac{1}{\alpha}  \\
&\E\left[W_i \ |\ \text{ride \& } N_i = n\right] := \indic{n \leq \N} \frac{1}{\gamma n}
\end{align*}
Here, $N_i$ denotes the steady state number of waiting drivers on platform~$i,$ and $\N$ is a large (but arbitrary) positive integer. This may be interpreted as follows. The average waiting time when a ride is taken is assumed to be inversely proportional to the number of available drivers. (This would be the case if each driver's travel time to the passengers's location is assumed to be independent and exponentially distributed with mean~$1/\gamma,$ and the platform assigns the ride to the driver who can arrive at the passenger's location soonest.) However, when the number of available drivers is large (formally, greater than $\N$), the average wait time is approximated to be zero. (This assumption, while reasonable, is made purely for analytical convenience.) 
Taking~$\gamma = 1$ without loss of generality, we have  (see~\eqref{eqn_expression_of_PBi},~\eqref{Eq_overall_exact}),
\begin{align}
\label{eqn_wait_metric_original_form}
\E\left[ W_i\right] 
    &= \frac{1}{\alpha} \PB_i(\VPhi; \lambda_i, \beta )  +  \E\left[\left. \tau(N_i)\right| N_i > 0\right] f(\phi_i) \Prob(N_i > 0) \\ 
    \label{eqn_WB_scratch}
    &= (1 - f(\phi_i)) \frac{1}{\alpha} + f(\phi_i) \E[\tau(N_i)], \mbox{ where }  \\
    & \tau(n) :=  
    \begin{cases}
      \frac{1}{\alpha}   &  \mbox{if } n = 0 , \\
     \left(\frac{1}{n}\right) \indic{n \leq \N} &  \mbox{else.} 
    \end{cases}
    \nonumber
\end{align}

Finally, we define the QoS metric under consideration as $\WB_i = \alpha \E\left[ W_i\right].$ (Note that the multiplicative scaling by~$\alpha$ does not impact the Wardrop split or the platform payoffs.) Formally, 
\begin{align} \label{eqn_wait_metric_defn}
 \WB_i(\VPhi; \lambda_i, \beta , \alpha) 
&= \PB_i(\VPhi; \lambda_i, \beta) + \alpha f(\phi_i)\frac{\sum_{n=1}^\N  \frac{\mu_{i,n}}{n} }{\sum_{n=0}^\infty  \mu_{i,n} }, 
\end{align}
where $$\mu_{i, n} := \prod_{a=1}^n\left(\frac{ e} {  \lambda_i f(\phi_i) + a\beta}\right)$$
(see Lemma~\ref{lem_BCMP_Product_form}). From~\eqref{eqn_wait_metric_defn}, it is intuitive (though far from trivial to prove formally) that the equilibria associated with the game with limiting payoff functions as $\beta,\alpha \rightarrow 0$ coincide with those established in Section~\ref{sec:overall_bp}, and these equilibria are also meaningful in the practically applicable pre-limit. In the remainder of this section, these intuitions are formalised.

\ignore{governed by the expected delay in receiving the (riding) service or the expected pick-up time. Any passenger who visits one of the platforms avails the riding service from them if the quoted price is acceptable and if at least one driver is available in the queue. In this case, the expected delay or the expected pick-up time   is the time it takes for one of the available drivers of  the visited platform to arrive at the passenger's location.
If the passenger has to venture outside of the two platforms,  we assume the delay it experiences   in receiving the service is given by $\nicefrac{1}{\alpha}$. 

Passengers arrive randomly at any spot within the zone, and waiting drivers are also distributed across the zone.
 When a passenger requests a ride, one of the waiting drivers of platform~$i$ (the quickest possible) arrives to the location of the passenger, where the travel/pick-up time of any driver is again exponentially distributed with rate $\gamma$ for some $\gamma > 0$.
Let $  W_i$ denote the delay/waiting a typical passenger experiences  at platform $i$. Conditioned that the quoted price is acceptable, and $n$ waiting drivers are available, the conditional expectation $\E[  W_i| N_i = n_i ] = \nicefrac{1}{n_i\gamma}$ (the expected value of the minimum of $n$ independent exponential random variables). On the other hand, conditioned on the event that the passenger has to avail service outside the platform, $E[  W_i | \text{outside ride}] = \nicefrac{1}{\alpha}$. \textit{We assume $\alpha \gamma < 1$ and without loss of generality set $\gamma = 1$.} 

Therefore, the expected delay for any typical passenger arriving to platform $i$, when the platform uses a static pricing policy $\phi_i$ and receives passengers at a rate $\lambda_i$, is given by (see~\eqref{eqn_expression_of_PBi},~\eqref{Eq_overall_exact}),
\begin{align}
\E\left[ W_i\right] 
    &= \frac{1}{\alpha} \PB_i(\VPhi; \lambda_i, \beta )  +  \E\left[\left.W_i\right| N_i > 0\right] f(\phi_i) \Prob(N_i > 0) 
\label{eqn_wait_metric_original_form}
\end{align}
The above equation can be rewritten as (see~\eqref{eqn_expression_of_PBi},~\eqref{Eq_overall_exact}),
\begin{equation} \label{eqn_WB_scratch}
\E[ W_i] = (1 - f(\phi_i)) \frac{1}{\alpha} + f(\phi_i) \E[\tau(N_i)], 
\mbox{ with }  \tau(n) :=  
    \begin{cases}
      \frac{1}{\alpha}   &  \mbox{if } n = 0 , \\
     \left(\frac{1}{n}\right) \indic{n \leq \N} &  \mbox{else.} 
    \end{cases}
 \end{equation}
Here $\tau(n)$ represents the expected pick-up time at an instance with $n$ waiting passengers; the expected pick-up time is almost zero when the number of waiting drivers is large and hence we set $\tau(n) = 0$ for $n > \N$ for some large enough $\N < \infty$. 
}

\ignore{
By Lemma \ref{lem_BCMP_Product_form}, the stationary (marginal) distribution of number of waiting drivers in platform $i$ is given by 
$$
\pi_i (N_i = n) = \sum_{r_i} \pi_i (n, r_i)
= \frac{ \mu_{i,n} }{\sum_{n'=0}^\infty  \mu_{i, n'}} \mbox{ where } 
\mu_{i, n} := \prod_{a=1}^n\left(\frac{ e} {  \lambda_i f(\phi_i) + a\beta}\right)
$$} 

%
%
 \ignore{
 Without loss of generality, and for ease of explanation, we use $\alpha \E[ W_i]$ as a QoS metric in this section. Therefore, the QoS   proportional to the expected delay in receiving the service for a typical passenger arriving at platform $i$, for $\beta > 0$ is given by :

 \begin{align} \label{eqn_exp_delay_expression}
\WB_i(\VPhi; \lambda_i, \beta , \alpha) &= 
 \PB_i(\VPhi; \lambda_i, \beta ) + \alpha f(\phi_i) \Prob(N_i > 0)  \E\left[\left.\frac{1}{N_i}\right| N_i > 0\right].
\end{align}

 
We would consider that the expected pick-up time $\tau(n)$ is zero if number of waiting drivers $n$  are more  i.e., if $n >\N$  for some $ \N < \infty$ and that the conditional expected value of the pick-up time is $\tau(n) = \nicefrac{1}{n}$, conditioned on   $N_i = n$. In this case, the expected delay at platform $i$ is given by:
} 


\ignore{
As before we consider IDP regime   for analysis (i.e., let $\beta \to 0$ in the above and define appropriate limit system for each $\alpha$). 
The  expression   for QoS  $\WB_i$  is simplified for $\beta > 0$ and  is derived for IDP regime, i.e., for $\beta = 0$ in the following lemma:
\begin{lemma}\label{lem_wait_metric_QoS}
Suppose $\rho_i = \nicefrac{e}{(\lambda_i f(\phi_i))}$, and $\mu_{i,n} = \prod_{a = 1}^n \left( \nicefrac{e}{(\lambda_i f(\phi_i) + a \beta)} \right)$. The QoS metric~$\WB_i$ of the platform~$i$ when $\beta > 0$ is
\begin{align} \label{eqn_QoS_WB_positive_beta}
\WB_i(\VPhi; \lambda_i, \beta , \alpha)
 &= 
     \PB_i(\VPhi; \lambda_i, \beta ) + \frac{\alpha f(\phi_i)\sum_{n=1}^{\N}  \frac{\mu_{i,n}}{n}}{\sum_{n=0}^{\infty} \mu_{i,n}} ,
\end{align}
and when $\beta = 0$ is
\begin{align}  \label{eqn_QoS_WB_zero_beta}
\WB_i(\VPhi;\lambda_i, 0 , \alpha)
&= \begin{cases}
  1 -  f(\phi_i)  + \alpha    f(\phi_i) \left(1 - \frac{e}{\lambda_i f(\phi_i)} \right) \left( \frac{1}{ \alpha} + \sum_{n=1}^\N  \frac{1}{n}\left(\frac{e}{\lambda_i f(\phi_i)}\right)^n\right) & \mbox{ if }     \frac{e}{\lambda_i f(\phi_i)} < 1 ,  \\
 1 -  f(\phi_i)    & \mbox{ else. } 
 \end{cases}
\end{align}
\end{lemma}

\begin{proof}[Proof of Lemma~\ref{lem_wait_metric_QoS}]
Using \eqref{Eq_overall_exact} and~\eqref{eqn_wait_metric_original_form}, it suffices to simplify the conditional expectation given in \eqref{eqn_wait_metric_original_form}. Thus,
\begin{align*}
\Prob(N_i > 0) \E\left[\left.\frac{1}{N_i}\right| N_i > 0\right]
=
\Prob(N_i > 0) \sum_{n =1}^{\N} \frac{\Prob(N_i = n)}{n \Prob(N_i > 0)} 
=  
\frac{ \sum_{n=1}^{\N} \frac{\mu_{i,n}}{n}}{ \sum_{n=0}^{\infty} \mu_{i,n}} .
\end{align*}
Using Lemma \ref{lem_PB_joint_continuity_from_scratch}, $\nicefrac{1}{\sum_{n=0}^{\infty} \mu_{i,n}}$ is jointly continuous in $(\VPhi, \lambda_i, \beta)$ for each $i$, hence the RHS of the above equation is also jointly continuous as $\N$ is finite. Thus at~$\beta = 0$, it is easy to find the sum of the series $\sum_{n = 0}^{\infty} \mu_{i,n} $ for $\nicefrac{e}{\left(\lambda_i f(\phi_i)\right)} < 1$ and the series diverges for $\nicefrac{e}{\left(\lambda_i f(\phi_i)\right)} \geq 1$.
\ignore{\begin{align} \label{eqn_WA_part_proof_1}
 \frac{\sum_{n=1}^{\infty} \mu_{i,n}}{\sum_{n=0}^{\infty} \mu_{i,n}} = \frac{\sum_{n=1}^{\infty} \varrho^{n} \left(\frac{e}{\lambda_i f(\phi_i)}\right)^n }{\gamma_i \sum_{n=0}^{\infty} \left(\frac{e}{\lambda_i f(\phi_i)}\right)^n}
= \begin{cases}
\frac{ \varrho \left(\frac{e}{\lambda_i f(\phi_i)}\right) (1 - \left(\frac{e}{\lambda_i f(\phi_i)}\right)) }{(1 - \varrho \left(\frac{e}{\lambda_i f(\phi_i)}\right))}    & \text{if } \left(\frac{e}{\lambda_i f(\phi_i)}\right) < 1, \\
0  &\text{else.} 
\end{cases}
\end{align}}
\end{proof}
Observe that for $\nicefrac{e}{(\lambda_i f(\phi_i))} \geq 1,$ the driver unavailability probability $\DA_i(\VPhi; 0) = 0$; in fact the probability that number of drivers is less than $\N +1$ is also zero. Consequently, in \idp~regime, the system experiences an overflow of drivers, resulting in an insignificant expected delay for the arriving passengers. Therefore, in this case, the passengers prefer local rides or other options only when the quoted price is not acceptable. As before, the revenue rate of platform~$i$, when the two platforms operate with the pricing policy~$\VPhi,$ equals~$\MR_i (\VPhi) = \MR_i (\VPhi; \lambda_i(\VPhi; \beta, \alpha))$, where the WE splits~$\lambda_i(\VPhi; \beta, \alpha)$ are defined using the QoS metric~$\WB_i$. In the following Theorem, we first characterize the limiting WE as $(\beta, \alpha) \to (0, 0)$ and as $\beta \to 0$ for a fixed $\alpha$ and, subsequently, the limiting payoff functions (proof in Appendix~\ref{sec:appendix_wait_metric}).
} 

First, we characterize the limiting payoff functions as $\beta, \alpha \rightarrow 0$ (similar to how we considered the \idp\ scaling $\beta \rightarrow 0$ in the preceding section). This is achieved by establishing a continuous extension of the Wardrop split and the payoff functions over~$(\beta,\alpha) \in [0,\infty)\times[0,1).$ 
\begin{thm}
\label{thm_WE_MR_new}
Fix price vector~$\VPhi$.
For~$\alpha = 0,$ $\beta \geq 0,$ define the WE and the platform payoffs as per Table~\ref{tab:WE_MR_PB}. For~$\beta =0,$ $\alpha > 0,$ define
\begin{align}  
\label{eqn_QoS_WB_zero_beta}
\WB_i(\VPhi;\lambda_i, 0 , \alpha) := 1 -  f(\phi_i) + \indic{ \frac{e}{\lambda_i f(\phi_i)} < 1} f(\phi_i) \left(1 - \frac{e}{\lambda_i f(\phi_i)} \right) \left(1 + \alpha \sum_{n=1}^\N  \frac{1}{n}\left(\frac{e}{\lambda_i f(\phi_i)}\right)^n\right).
\end{align}
Now, define the WE for this case as per Table~\ref{tab:WE_MR_WB}, with
\begin{equation}\label{eqn_bar_alpha}
    \bar \alpha = \frac{\left( \frac{e}{\Lambda } - f(\phi_1) \right)}{\left(\frac{e}{\Lambda } - f(\phi_2)  \right) \sum_{n=1}^\N  \frac{1}{n}  \left( \frac{ e}{\Lambda  f(\phi_2)} \right)^n }, 
\end{equation}
and the platform payoffs (via Lemma~\ref{lem_approx_MR}) as 
$$\MR_i(\VPhi; 0, \alpha) = e\phi_i \indic{e<  f(\phi_i) \lambda_i (\VPhi; 0, \alpha)} + \lambda_i (\VPhi; 0, \alpha)  f(\phi_i) \phi_i \indic{e \ge  f(\phi_i) \lambda_i (\VPhi; 0, \alpha)}.$$
Under these definitions, the mappings $(\beta, \alpha) \mapsto \lambda_i(\VPhi; \beta, \alpha)$ 
 and~$(\beta, \alpha) \mapsto \MR_i(\VPhi;\beta, \alpha)$ are continuous over~$[0, \infty) \times [0, 1)$ for each $i.$
 \begin{table}[ht]
\renewcommand{\arraystretch}{1.4}
\centering
\begin{tabular}{|c|cc|c|}
\hline
    \multirow{2}{*}{}     & \multicolumn{2}{c|}{Range of~$\phi_1$ when~$\phi_1 > \phi_2$} & \multirow{2}{*}{$\phi_1 = \phi_2$} \\  \hhline{~--~}
    %
    & \multicolumn{1}{c|}{$ [0, f^{-1}(\rho))$} & $ [f^{-1}(\rho), \phi_h]$ & \\  \hline 
$\lambda_1(\VPhi; 0, \alpha)$ & \multicolumn{1}{c|}{$\tilde \lambda^*(\alpha)$}                 & ${\tilde \lambda}^*(\alpha) \indic{\alpha > \bar \alpha}$                                  & \multicolumn{1}{c|}{$\nicefrac{\Lambda}{2}$}            \\ \hline
%
%
    %
    %
    \end{tabular}
    \caption{Here, ${\tilde \lambda}^*(\alpha)$ is the unique minimiser of $\left(\WB_1(\VPhi;\lambda_1, 0 , \alpha) - \WB_2(\VPhi;\Lambda - \lambda_1, 0 , \alpha)\right)^2,$ where $\WB_i(\VPhi;\lambda_i, 0 , \alpha)$ is as defined in~\eqref{eqn_QoS_WB_zero_beta}.
     \label{tab:WE_MR_WB}
    }
    \end{table}
\end{thm}
It follows from~Theorem~\ref{thm_WE_MR_new} (proof in Appendix~\ref{sec:appendix_wait_metric}) that the limiting payoff functions as $\beta,\alpha \rightarrow 0$ (note that this is a two-dimensional limit) coincide with those obtained in Section~\ref{sec:overall_bp} under the $\PB$ QoS metric as $\beta \rightarrow 0.$ Thus, the equilibrium characterizations in Theorems~\ref{thm_B_sys_limit_sys_main} and~\ref{thm_EC} are also applicable to the (limiting) game considered here. Crucially, we are also able to establish that these equilibria are $\epsilon$-equilibria in the (two-dimensional) pre-limit.
%
%
\begin{thm}
\label{thm_new_system_epsilon_equilibrium}
Consider~$\rho < 1$. Then for any~$\epsilon > 0$, there exists~${\bar \alpha}_\epsilon, {\bar \beta}_\epsilon > 0$ such that for all $0 < \alpha \le {\bar \alpha}_\epsilon$, $0 < \beta \le {\bar \beta}_\epsilon$, 
\begin{enumerate}[$(1)$]
    \item the NE defined in Theorem~\ref{thm_B_sys_limit_sys_main} (both pure and mixed) are~$\epsilon$-NE, and
    \item if~$f(\phi_b)/2 < \rho < 1,$ then the interval~$[\eL, \eU]$ defined in Theorem~\ref{thm_EC} is an~$\epsilon$-EC.
\end{enumerate}
\end{thm}
Theorem~\ref{thm_new_system_epsilon_equilibrium} (proof in Appendix~\ref{sec:appendix_wait_metric}) further reinforces the applicability of the equilibrium characterizations in Theorems~\ref{thm_B_sys_limit_sys_main} and~\ref{thm_EC}, establishing that these are robust to driver impatience as well as the impact of passenger pick-up delays. This robustness is further demonstrated via numerical case studies in the following section.
}
 \section{Comparisons and numerical results}
\label{sec:comparison}

In this section, we compare the equilibria characterized in the previous sections (under different assumptions on how the passenger arrival rate gets split between platforms), both analytically as well as via numerical experiments. We also demonstrate how our results under the analysed limiting regimes are applicable in the pre-limit.

\subsection{Comparisons}

The following lemma compares the monopoly optimum prices derived in Section~\ref{sec:monopoly} and the equilibrium prices derived in Section~\ref{sec:overall_bp}. \rev{Recall that the same equilibria are also meaningful under the pick-up time QoS metric analysed in Section~\ref{sec:delay}.}

\begin{lemma}
\label{lem_mono_price_less_than_D}
Consider the \idp\ regime ($\beta = 0$). The equilibrium price under QoS metric~$\PB$ is less than or equal to the optimal monopoly price.\footnote{Under the conditions where the equilibrium is characterized as a mixed NE or EC, the above dominance holds for all points in the support of the NE/EC.}
\ignore{\begin{enumerate}[$(1)$]
    \item The equilibrium price under QoS metric~$\DA$ is greater than or equal to the optimal monopoly price.
    \item The equilibrium price under QoS metric~$\PB$ is less than or equal to the optimal monopoly price.\footnote{Under the conditions where the equilibrium is characterized as a mixed NE or EC, the above dominance holds for all points in the support of the NE/EC.}
\end{enumerate}}

Moreover, the optimal platform payoff in the monopolistic scenario is always greater than or equal to the equilibrium payoff. In this comparison, we treat the equilibrium price and platform payoff under QoS metric~$\PB$ to be zero for $\rho \geq 1,$ in line with Statement~(3) of Theorem~\ref{thm_B_sys_limit_sys_main}.
\end{lemma}

Lemma~\ref{lem_mono_price_less_than_D} (proof in Appendix~\ref{sec:appendix_comparison}) provides valuable insights into the interplay between competition, QoS and platform utility. In particular, we find that competition drives platforms to set prices that deviate from the optimal monopoly price, so as to enhance the QoS metric that governs the market share split. 
Indeed, when the market split is governed by QoS metric~$\PB,$ equilibrium prices are lower than the monopoly price, given that stationary probability that an arriving passenger obtains a ride gets boosted by a price reduction. Moreover, we see that inter-platform competition enhances QoS on the passenger side, while diminishing the equilibrium payoff of each platform.

Next, we numerically evaluate and compare the optimal/equilibrium prices and corresponding platform utilities across the monopoly and competition models.
 
We consider a non-linear price sensitivity function,~$f(\phi)= 1 - (a \phi)^2$ for some~$a < 1/ \phi_h$ (this choice satisfies Assumptions~{\bf A}.1-2). By performing simple algebraic calculations, we compute the following quantities defined in Theorems~\ref{thm_monopoly_opt}--\ref{thm_B_system_epsilon_equilibrium}, which characterize the different equilibria.
\begin{table}[ht]
\renewcommand{\arraystretch}{2.1}
\begin{tabular}{lll}
$\phi_m = \frac{1 }{\sqrt{3}a}$, & $\phi_b = \frac{1 }{\sqrt{5}a}$, & $f^{-1} (2\rho) = \left(\frac{1}{a}\right) \sqrt{ 1 - 2\rho }$, \\
$\phi_u = \frac{1 }{\sqrt{2}a}$, & $\eL = \left(\frac{2}{3a}\right)\left(\frac{1}{\rho} - 1 \right) \sqrt{\frac{1 - \rho}{3}}$, & $\phimr = \left(\frac{1}{a}\right) \sqrt{\frac{1 - \rho}{3}}$. 
\end{tabular}
\end{table}

\begin{figure}[ht]
    \centering
    \subfloat[\centering Equilibrium/optimal price versus~$\rho$ ]{{\includegraphics[trim = {1.5cm 7.2cm 2.6cm 8cm}, clip, scale = 0.33]{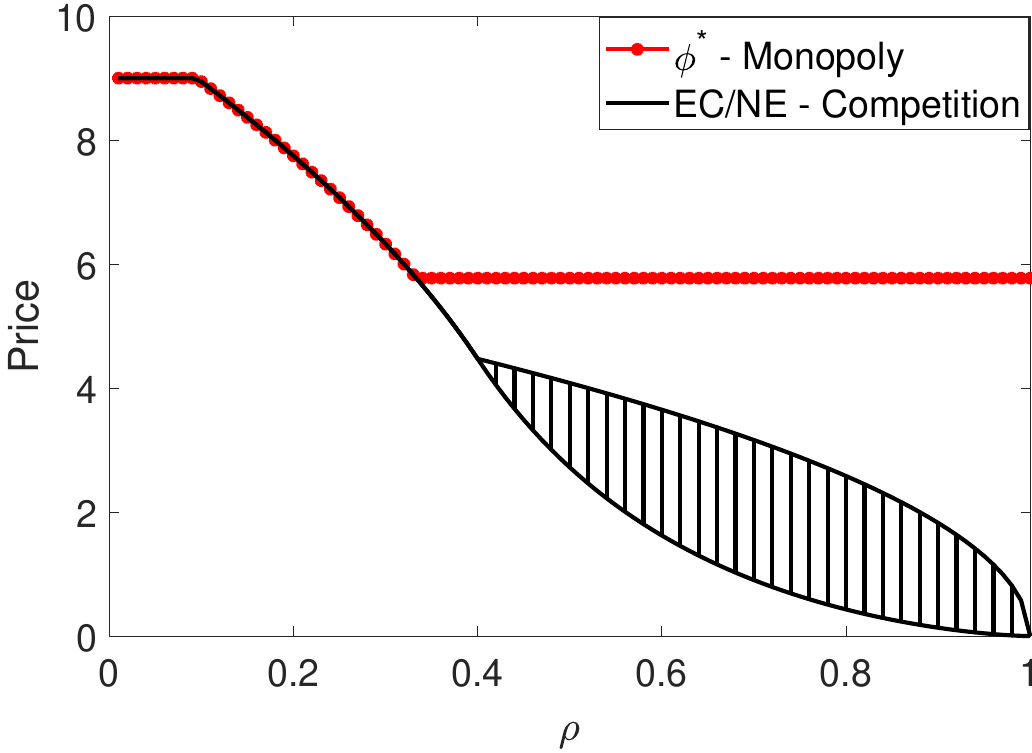} }}%
    \hspace{0.8cm}
    \subfloat[\centering Equilibrium/optimal platform payoff versus~$\rho$]{{ \includegraphics[trim = {1.5cm 7.2cm 2.6cm 8cm}, clip, scale = 0.33]{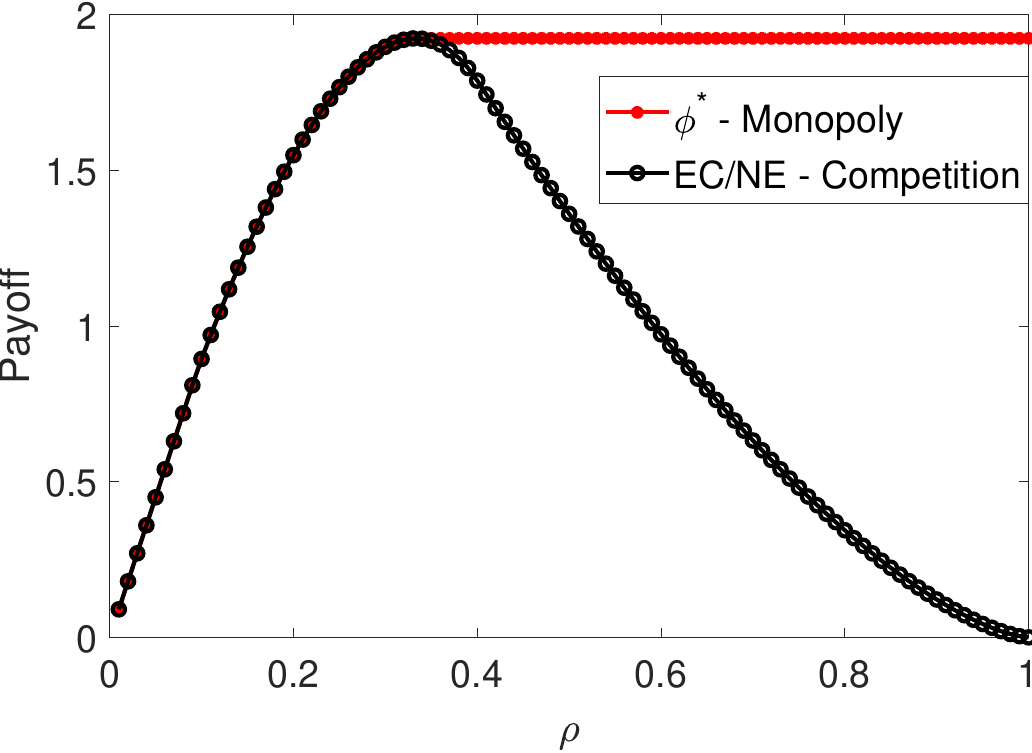} }}%
    \caption{Price sensitivity function~$f(\phi) = 1 - (a\phi)^2,$ where $a = 0.1$ and~$\phi_h = 9,$ $\Lambda = 1.$}%
    \label{fig:Square}%
\end{figure}

\noindent In Figure~\ref{fig:Square}, we plot the optimal/equilibrium price (left panel) and platform utility (right panel) as a function of the~\dpr~$\rho$ (specifically, we hold $\Lambda$ fixed and vary~$e$). 
The support of the mixed NE (or the EC) is marked via vertical black lines in the left panel; in the right panel, we plot the average payoff corresponding to the mixed NE (or equivalently, the security value; see Lemma~\ref{lem_mixed_NE_payoff_and_sec_value}). Note that the pricing is in line with the conclusions of Sections~\ref{sec:monopoly}--\ref{sec:overall_bp} and Lemma~\ref{lem_mono_price_less_than_D}. It is instructive to note the variation of the optimal/equilibrium platform utility as $\rho$ increases (specifically, as $e$ increases with $\Lambda$ being fixed): 
\begin{itemize}
\item Platform utility first increases and then saturates, in the monopoly setting.
This is because as driver availability grows, platforms are able to match more passengers with rides, until driver availability is no longer a constraint. 
\item In contrast, under inter-platform competition, 
platform utility first increases and then decreases, approaching zero as $\rho$ approaches~1. This is because competition drives the equilibrium price close to zero as $\rho$ approaches~1.
\item 
The \textit{price of anarchy} (from the standpoint of the platforms) exceeds one when~$\rho$ is large, i.e., when passengers are scarce relative to drivers.
\end{itemize}


\subsection{Cooperation versus competition}

We now study an alternate \textit{cooperative} scenario in which the platforms operate together, by pooling their driver resources to jointly serve the entire passenger demand. Mathematically, the analysis of this scenario is identical to that corresponding to the monopoly scenario (see Section~\ref{sec:monopoly}), but with passenger arrival rate $\nicefrac{\Lambda}{2}$ being replaced by~$\Lambda,$ and effective driver arrival~$e$ being replaced by~$2e.$ 
 \begin{figure}[ht]
    \centering
    \subfloat[\centering Optimal price versus~$\beta$]{{\includegraphics[trim = {1.5cm 7.2cm 2.5cm 8cm}, clip, scale = 0.33]{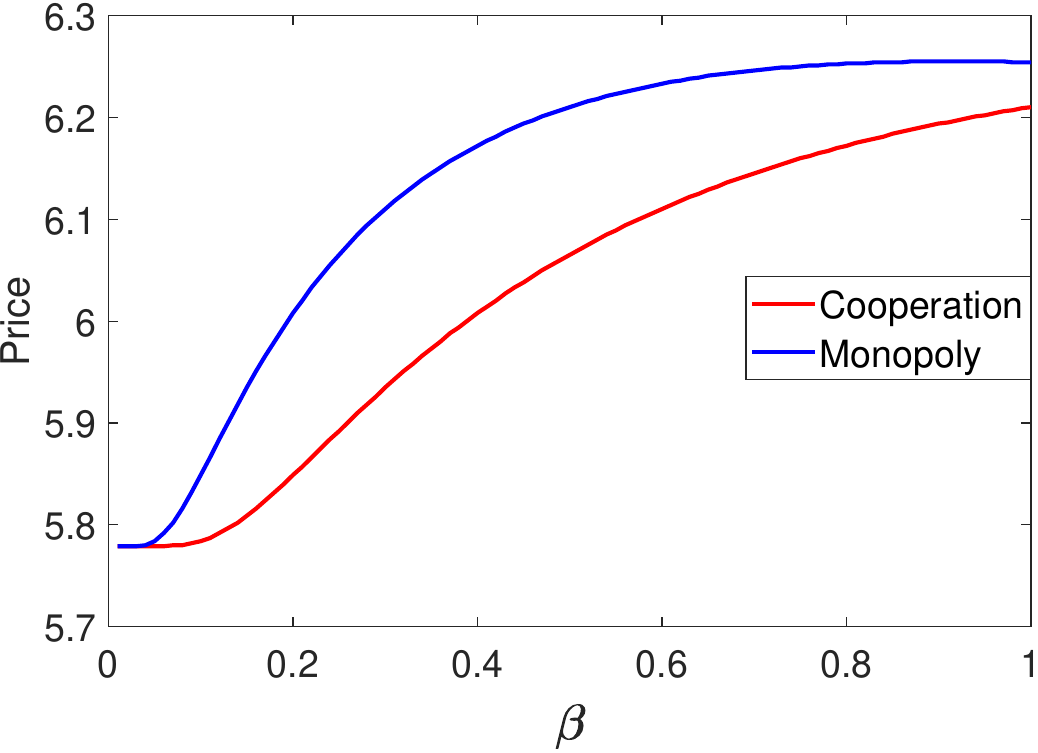} }}%
    \hspace{0.8cm}
    \subfloat[\centering Per-platform payoff versus $\beta$]{{\includegraphics[trim = {1.5cm 7.2cm 2.5cm 8cm}, clip, scale = 0.33]{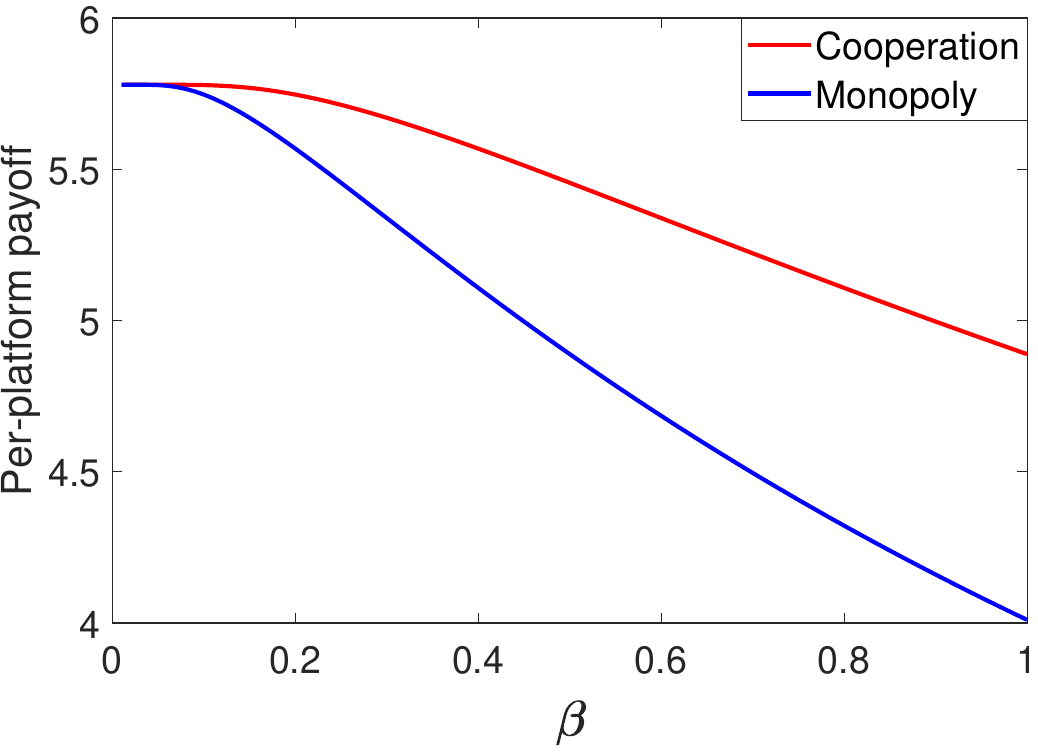} }}%
    \caption{Comparison between cooperative setting and monopoly setting; here, $e  = 2$,~$\Lambda = 3$,~$\phi_h = 10,$~$f(\phi) = 1 - (a \phi)^2$, with~$a = \nicefrac{1}{10.01}$}%
    \label{fig:monopoly_coop}%
\end{figure}
Since the \dpr~$\rho$ remains unchanged under the above scaling, the optimal pricing in the \idp\ regime ($\beta = 0)$ for this cooperative setting is exactly as stated in Theorem~\ref{thm_monopoly_opt}. This also means that the payoff of each platform (assuming an equitable revenue split) is also identical to that in the monopoly setting. In light of Lemma~\ref{lem_mono_price_less_than_D}, we conclude that \textit{cooperation is beneficial to the platforms} (while being detrimental to the passenger base).

Interestingly, the above observation suggests that in the \idp\ regime (and via continuity, also when~$\beta$ is small, i.e.,  assuming drivers are patient), the economies of scale one typically expects from resource pooling in queueing systems do not arise. However, when~$\beta$ is large, i.e., drivers are relatively impatient, we find that the platforms do indeed benefit from the economies of scale resulting from resource pooling---see Figure~\ref{fig:monopoly_coop}, which compares the optimal price and (per-platform) payoff in the monopoly setting (where each platform sees an exogenous passenger arrival rate of $\nicefrac{\Lambda}{2}$) and the cooperative setting (where both platforms jointly serve passengers at rate~$\Lambda$). Note that the per-platform payoff is higher in the cooperative setting (and this is achieved by operating at a lower price) compared to the monopoly  setting when $\beta$ is large. Given our previous observation that competition results in an equilibrium payoff less than the monopoly payoff, this suggests that cooperation is even more beneficial to the platforms if drivers are impatient.

\subsection{Numerical illustrations of the accuracy of limiting regime equilibria in the pre-limit}

\rev{We now provide some numerical evidence to complement the conclusions of Theorems~\ref{thm_B_system_epsilon_equilibrium} and~\ref{thm_new_system_epsilon_equilibrium}, which show that the equilibria characterized under certain limiting regimes (the \idp\ regime $\beta \rightarrow 0,$ and $(\beta,\alpha) \rightarrow (0,0),$ respectively) are also meaningful in the pre-limit.}

\begin{figure}[ht]
    \centering
    \subfloat[\centering {QoS metric $\PB$} ]{{\includegraphics[trim = {1.4cm 7.4cm 2.2cm 8.1cm}, clip, scale = 0.33]{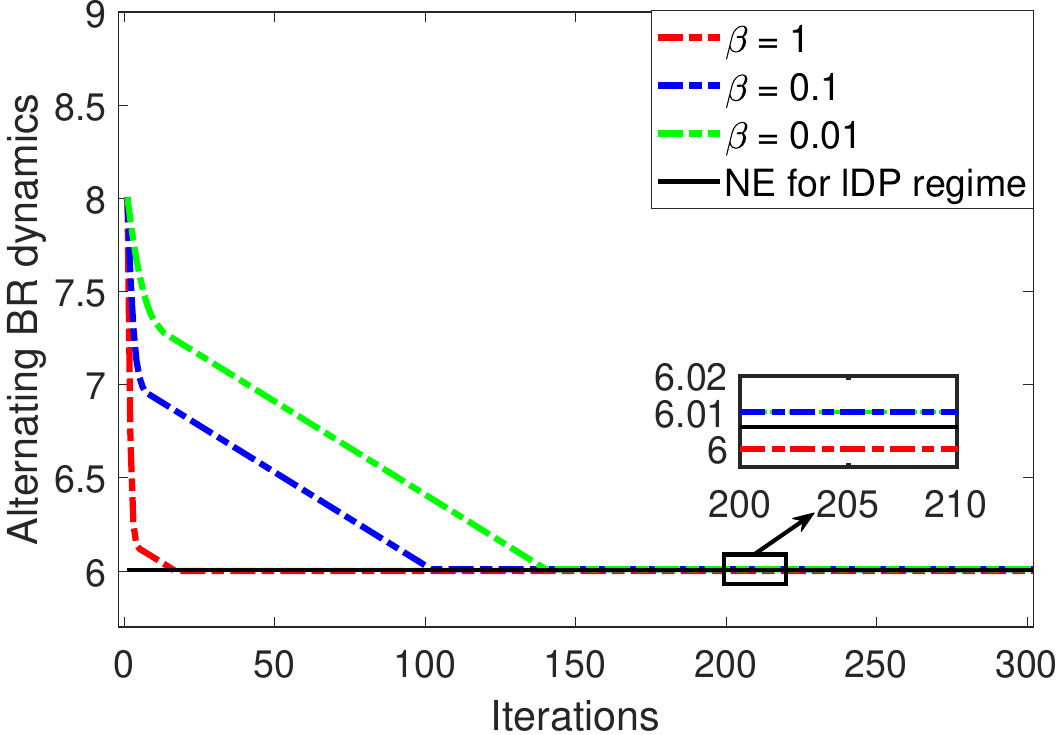} \label{fig:BR_dynamics_NE_PB} }}%
    \hspace{0.8cm}
    \subfloat[\centering {QoS metric $\WB$}]{{\includegraphics[trim = {1.2cm 8.2cm 2cm 9cm}, clip, scale = 0.33]{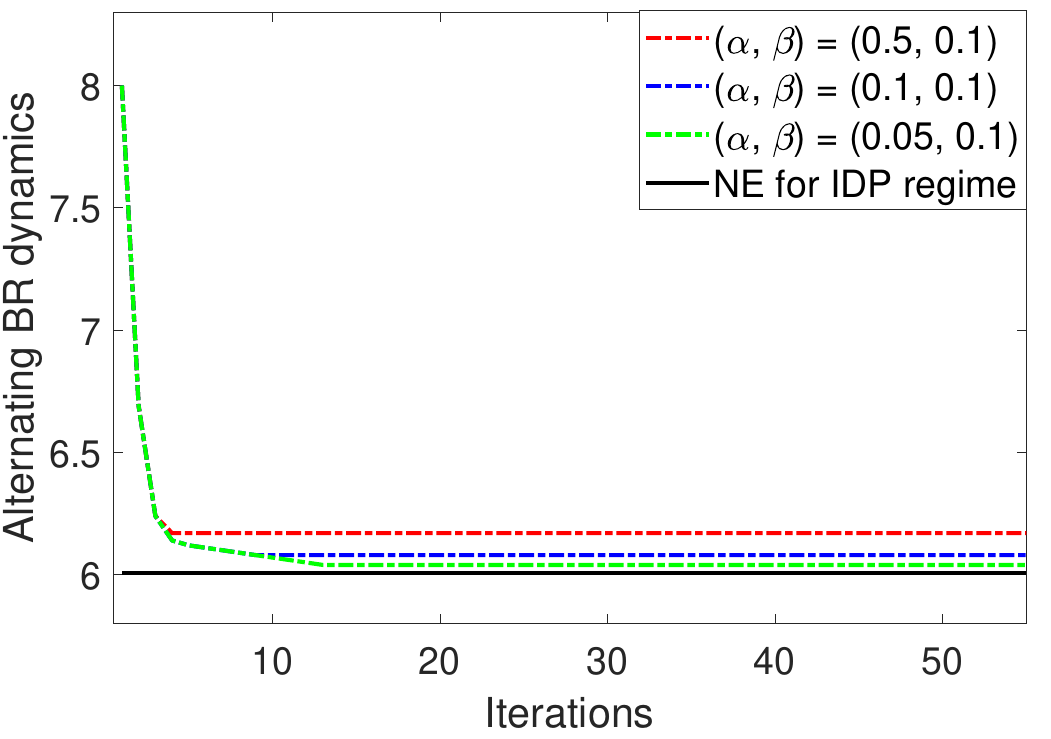} \label{fig:BR_dynamics_NE_WB} }}%
    \caption{Alternating best response dynamics; $ e = 1, \Lambda = 5, \phi_h = 10, f(\phi) = 1- (a \phi)^2$ with~$a = \nicefrac{1}{10.01}$.}%
    \label{fig:BR_dynamics_NE}%
\end{figure}

\rev{
In Figure~\ref{fig:BR_dynamics_NE_PB}, we plot the actions (prices) corresponding to alternating best response (BR) dynamics between the platforms for different choices of~$\beta.$\footnote{At each odd iteration, the price/action of Platform~1 is optimized, holding the price/action of Platform~2 fixed to its value in the previous iteration. Similarly, at each even iteration, the price/action of Platform~2 is optimized, holding the price/action of Platform~1 fixed to its value in the previous iteration.}  For reference, the figure also shows the NE corresponding to the \idp\ regime as characterized in Theorem~\ref{thm_B_sys_limit_sys_main}.
\rev{Similarly, in Figure~\ref{fig:BR_dynamics_NE_WB}, we plot the actions corresponding to alternating BR dynamics between the platforms for different choices of~$(\alpha, \beta).$} There are two key takeaways from these figures:
\begin{itemize}
\item Alternating BR dynamics appear to converge when~$\beta$ and $\alpha$ are small and positive. This suggests the existence of a pure NE even in the pre-limit.
\item The above limit is close to the NE corresponding to the limiting regime. This in turn suggests that the NE in the pre-limit is `close' to the NE in the limiting regime.
\end{itemize}
}

\begin{figure}[ht]
    \centering
    \subfloat[\centering Dynamics initialized above the \idp\ regime EC]{{\includegraphics[trim = {1.4cm 7.2cm 2.1cm 8cm}, clip, scale = 0.33]{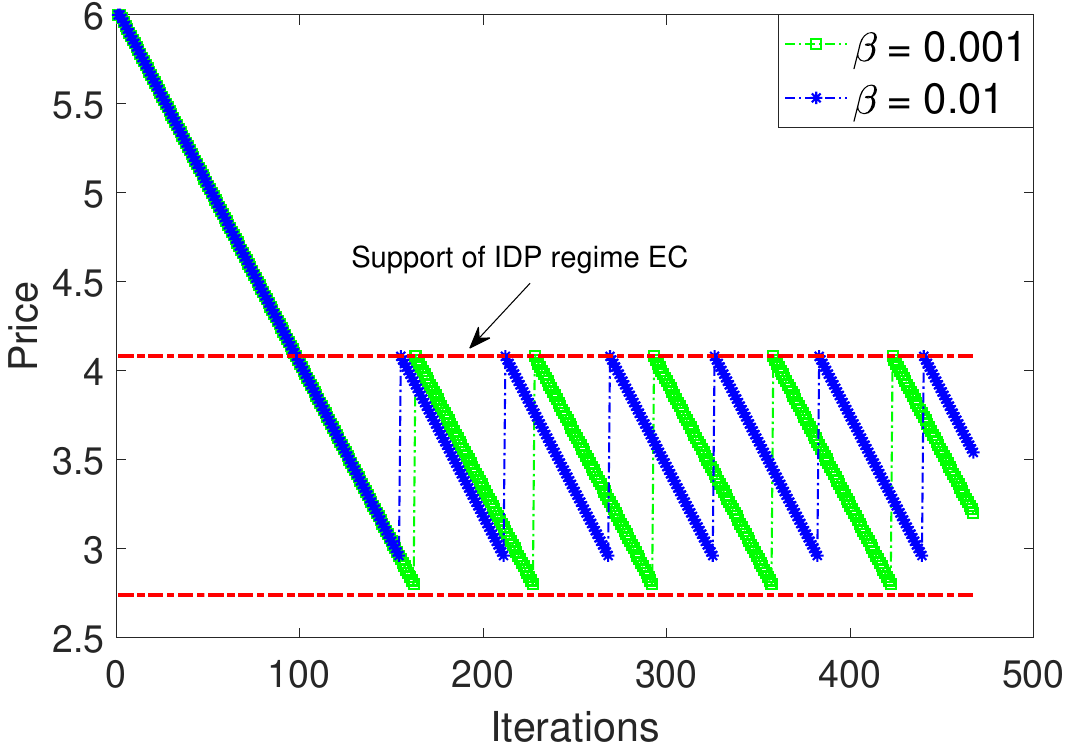} \label{fig:BR_dynamics_EC_PB_1} }}%
    \hspace{0.8cm}
    \subfloat[\centering Dynamics initialized below the \idp\ regime EC]{{\includegraphics[trim = {1.4cm 7.2cm 2.1cm 8cm}, clip, scale = 0.33]{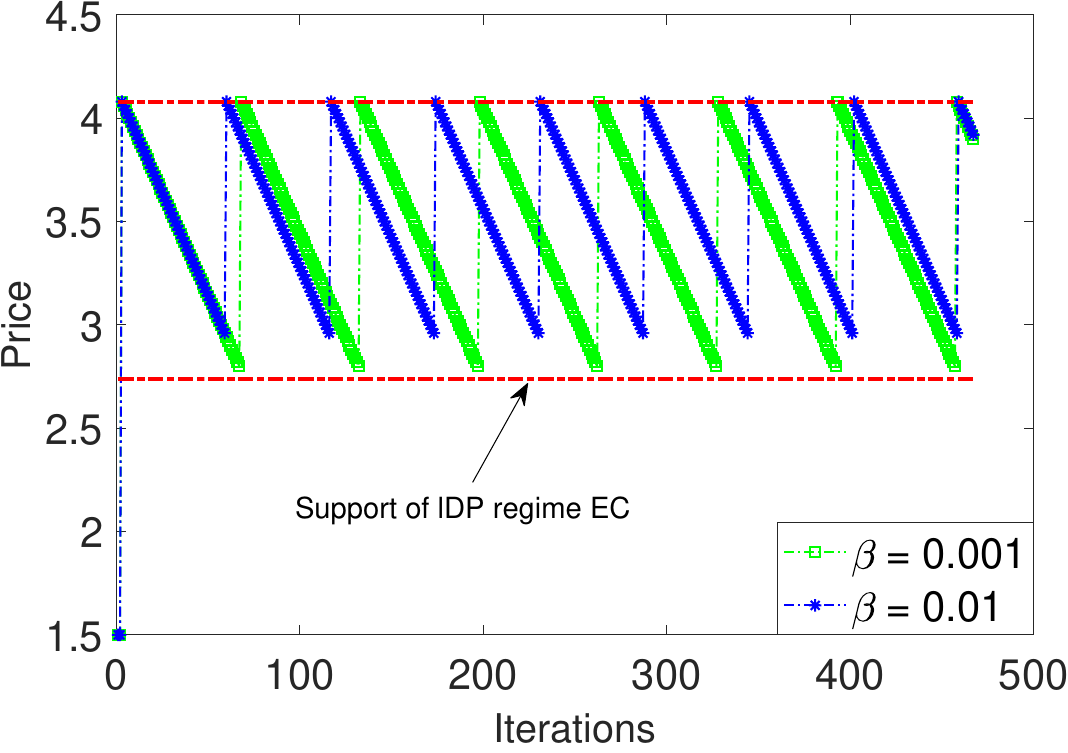} \label{fig:BR_dynamics_EC_PB_2} }}%
    \hspace{0.5cm}
    \subfloat[\centering Dynamics initialized above the \idp\ regime EC]{{\includegraphics[trim = {1.5cm 8.2cm 1.9cm 8.8cm}, clip, scale = 0.34]{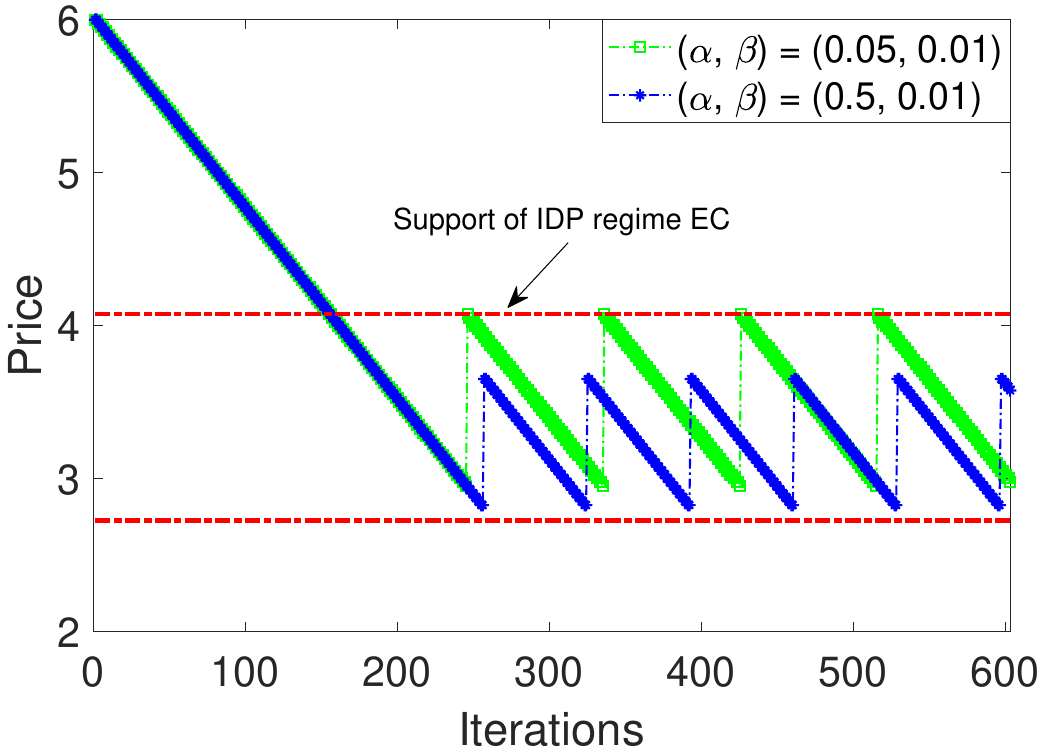} \label{fig:BR_dynamics_EC_WB_1} }}%
    \hspace{0.5cm}
    \subfloat[\centering Dynamics initialized below the \idp\ regime EC]{{\includegraphics[trim = {1cm 8.2cm 1.8cm 8.8cm}, clip, scale = 0.34]{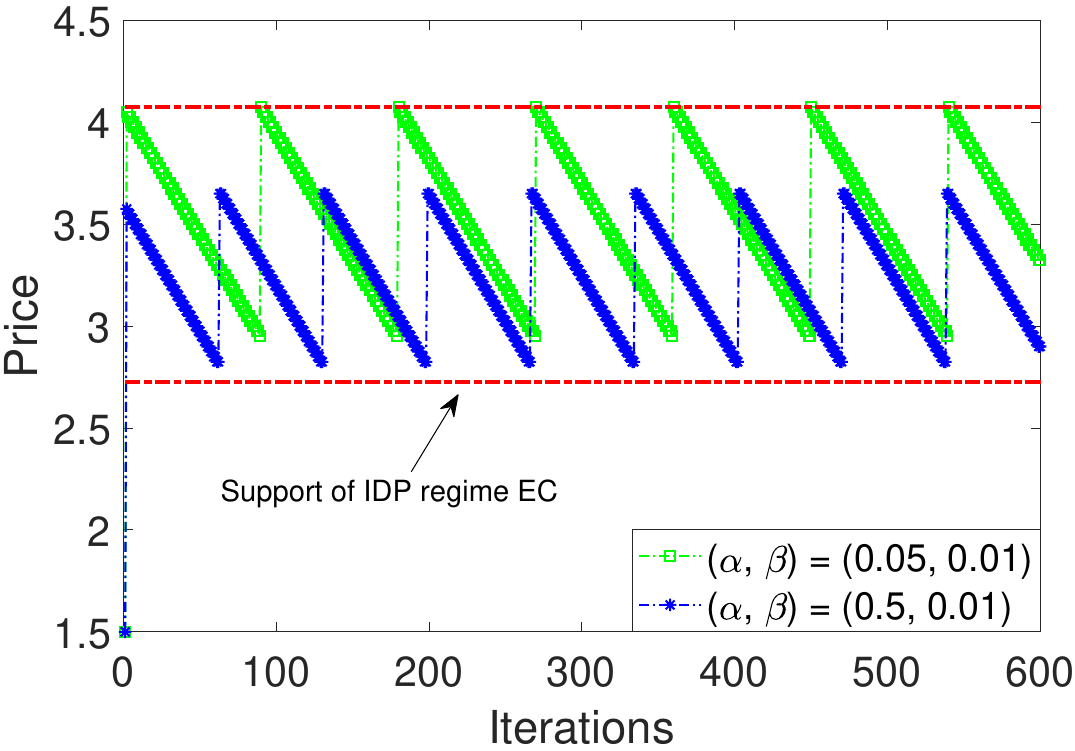} \label{fig:BR_dynamics_EC_WB_2} }}%
    \caption{Alternating BR dynamics under QoS metric $\PB$ (Figures~\ref{fig:BR_dynamics_EC_PB_1},~\ref{fig:BR_dynamics_EC_PB_2}) and QoS metric  $\WB$ (Figures~\ref{fig:BR_dynamics_EC_WB_1},~\ref{fig:BR_dynamics_EC_WB_2}) starting from two different points outside the \idp\ regime EC;~$e = 1$,~$\Lambda = 2$,~$f(\phi) = 1- (a \phi)^2$ with~$a = \nicefrac{1}{10.01}$.}%
    \label{fig:BR_dynamics_EC}%
\end{figure}

\rev{
In Figure~\ref{fig:BR_dynamics_EC}, we plot the actions (prices) corresponding to alternating best response (BR) dynamics between the platforms for different choices of~$(\alpha, \beta),$ for a choice of system parameters that induce an equilibrium cycle in the limiting regime. (Note that even though best responses do not exist in the limiting regime, they do exist in the pre-limit.) Panels~(a) and~(b) show the results for QoS metric~$\PB,$ and panels~(c) and~(d) show the results for QoS metric~$\WB.$ Interestingly, we find that the best response dynamics oscillate in the pre-limit; moreover, the support of this oscillation aligns closely with the equilibrium cycle in the limiting regime as characterized in Theorem~\ref{thm_EC} (the limiting regime EC is depicted using dotted red lines in Figure~\ref{fig:BR_dynamics_EC}). These observations are consistent with our dynamic interpretation of the equilibrium cycle in Section~\ref{sec:overall_bp}.}


\section{Concluding Remarks}
\label{sec:conclusion}

We consider competing ride-hailing platforms with impatient, price-sensitive passengers and impatient and revisiting drivers. On the one hand, each platform strives to meet its passengers' QoS goals so as to capture a larger market share (the total market share across platforms is conserved). On the other hand, they also seek to increase their prices so as to maximize their long-run revenue rate. Our analysis thus sheds light on the intricate interplay between competition, price-sensitivity, and the relative scarcity/surplus of drivers relative to passengers.

A key feature of our analysis is that we provide \textit{closed form} characterizations of the equilibrium prices in a certain limiting regime, which we refer to as the \textit{infinite driver patience} regime. Interestingly, we find that the equilibrium prices depend on (i) the price sensitivity of the passenger base, and (ii) the ratio of driver and passenger arrival rates. This latter quantity appears to play a role analogous to the utilization/load in classical queueing models. In particular, we find that competition drives platforms to set prices that deviate from the `optimal' values when passengers are price-sensitive and scarce (relative to drivers); moreover, the deviation seeks to enhance the QoS experienced by the passenger base. \rev{Crucially, we show that the above equilibria are robust to driver impatience, as well as pick-up delay considerations on part of passengers.}

From a game theoretic standpoint, our analysis also demonstrates some interesting equilibrium behavior. Specifically, in certain cases, we find that discontinuities in payoff functions can result in there being no pure Nash equilibrium. Instead, we demonstrate a mixed equilibrium, as well as a novel equilibrium concept with interesting dynamic connotations that we refer to as an \textit{equilibrium cycle}. 

\rev{Key limitations of the present study include the assumptions of symmetric platforms (made to simplify the equilibrium characterization), and non-strategic drivers. Addressing the latter limitation by analysing the impact of driver retention on inter-platform competition represents an important and immediate avenue for future work. This paper also motivates future work in other directions, including multiple zones and dynamic/surge pricing.} We believe the BCMP-style modeling approach we have adopted is amenable to these extensions. Additionally, we also believe the concept of equilibrium cycle can be demonstrated in a broad class of practically motivated game formulations (for example, capturing competition between firms providing price-sensitive, substitutable services), where small price perturbations can result in significant market share variations.

\section*{Compliance with Ethical Standards}

The authors did not receive support from any organization for the submitted work. The authors have no relevant financial or non-financial interests to disclose. All authors have agreed with the content and have given explicit consent to submit this manuscript. The authors do not require consent from IIT Bombay to submit this manuscript for publication. This work does not involve experimentation with human/animal subjects. The authors have no conflicts of interest to declare.

\bibliography{reference.bib}
\appendix

\section{Proofs related to Section~\ref{sec:model}}\label{sec:appendix_A_model}

\begin{proof}[Proof of Lemma~\ref{lem_exist_unique_WE}]
We begin with the proof of the second part (monotonicity property of QoS functions), with $Q_i = \PB_i$.
W.l.o.g. consider $Q_1$, and any $\lambda_1 < \lambda_1'$. Observe,
$$\prod_{a=1}^{n_1}(\lambda_1' f(\phi_1) + a \beta) - \prod_{a=1}^{n_1}(\lambda_1 f(\phi_1) + a \beta) > 0, \text{ for any } n_1>0,$$
and hence (since both the series are convergent for any~$\beta > 0$),
\begin{eqnarray}
\sum_{n_1=0}^{\infty} \frac{e_1^{n_1}}{\prod_{a=1}^{n_1}(\lambda_1 f(\phi_1) + a \beta)} -
\sum_{n_1=0}^{\infty} \frac{e_1^{n_1}}{\prod_{a=1}^{n_1}(\lambda_1' f(\phi_1) + a \beta)} > 0. \nonumber
\end{eqnarray}  
Thus, from~\eqref{eqn_pi0_expression},~$\DA_1$  is strictly increasing and continuous in $\lambda_1$. From~\eqref{eqn_expression_of_PBi},
\begin{align*}
\PB_1 &= \DA_1 + (1 - f(\phi_1)) (1 - \DA_1)
= (1 - f(\phi_1)) + f(\phi_1)\DA_1.
\end{align*}
Clearly,~$\PB_1$ is an affine transformation of~$\DA_1$ with positive coefficients. Therefore,~$\PB_1$ is also strictly increasing and continuous, i.e., the QoS satisfy {\bf A}.3, when~$\beta > 0$.

Next we show the existence of unique WE for any QoS which satisfy {\bf A}.3. Define~$$g(\lambda) := Q_1(\lambda) - Q_2(\Lambda-\lambda) \text{ where } \lambda \in [0,\Lambda].$$ 
By \textbf{A.}3,~$g$ is a continuous function. 
Further:
\begin{enumerate}[$(i)$]
    \item If~$g(0) < 0$ and~$g(\Lambda) \geq 0$   then using Intermediate value theorem, there exists a~$\lambda_1 = \lambda^* \in (0, \Lambda)$ such that~$g(\lambda^*) = 0$.  Then~$(\lambda^*, \Lambda-\lambda^*)$ is WE. Further, by strict monotonicity, we have uniqueness. 
    \item If~$g(0) < 0~$ and~$g(\Lambda) < 0$. By {\bf A}.3, and definition of~$g$, we have~$g(\lambda) < 0$ for all~$\lambda$. Hence~$( \Lambda, 0)$ is the unique WE. 
    In  a similar way, when~$g(0) > 0$ and~$g(\Lambda) > 0$, we have that~$(0, \Lambda)$ is the unique WE.
\end{enumerate}
 Rest of the cases (example~$g(0) > 0$ and~$g(\Lambda) \leq 0$ etc.,) follow using similar arguments.
 \end{proof}

\begin{proof}[Proof of Lemma~\ref{lem_mr_derivation}]
In this proof, we omit the subscripts~$i$ related to platform~$i$ for simplicity. Consider a renewal process with renewal epochs being the points where the Markov process~$\{Z_{t}\}_{t \ge 0}$ with~$Z_t = (N_t, R_t)$, visits some state~$s = (n, r)$ s.t.~$n > 0$. Let the overall transition rate\footnote{Observe the passenger rejecting the offered price does not alter the state of the system.} from state~$s$ be given by~$q(s) = \lambda f(\phi) + r \nu + n \beta + \eta$. Then by the stability of the Markov process, the expected length of the corresponding renewal cycle~$\E(\tau(s))$ is (see~\cite{balter}, $\pi(s)$ is the stationary probability),
 \begin{equation}
 \label{eq_tau}
 \E[\tau(s)] = \frac{1}{q(s) \pi(s)}.
 \end{equation}
 In any renewal cycle, the platform obtains the revenue if a driver is available and the arriving passenger accepts the ride (i.e., offered price) when the system is in state~$s$. Let~$\R(s,t)$ denote the revenue generated till time~$t$, obtained by matching while in the state~$s$. This component of the reward (for each $s$ with $n>0$) can be obtained using the renewal process mentioned above. Towards this, let~$\R(s)$ represent the revenue generated in one renewal cycle of the corresponding renewal process.
 Thus, by well-known Renewal Reward Theorem (see~\cite{balter}), the long-run revenue rate while in state~$s$ is,
\begin{equation}
\label{eq_long_run_rev}
\lim\limits_{t \rightarrow \infty} \frac{\R(s, t)}{t} \xrightarrow{\text{RRT}} \frac{\E[\R(s)]}{\E[\tau(s)]} \ \text{a.s.}
\end{equation}
Observe that the expected reward generated in one renewal cycle, when in state~$s$,  is the price  offered by the platform ($\phi$) multiplied with the probability  $\nicefrac{f(\phi) \lambda}{q(s)}$ of a passenger arriving to the system after accepting the price, 
i.e., 
$$\E[\R(s) ] = \phi \frac{f(\phi) \lambda}{q(s)}.$$
Thus,~\eqref{eq_long_run_rev} can be simplified using~\eqref{eq_tau}:
\begin{align} \label{eqn_revenue_using_rrt}
 \lim\limits_{t \rightarrow \infty} \frac{\R(s, t)}{t} 
&= \frac{\phi  \frac{\lambda f(\phi)}{q(s)}}{\frac{1}{q(s)\pi(s)}}
= \lambda f(\phi) \phi  \pi(s)  \text{ a.s.} 
\end{align}

Now, define~$\R(t)$ as the revenue generated by the platform till time~$t$. Then clearly, 
$$\R(t) = \sum_{s: n>0} \R(s,t).$$ 
Hence, the long-run revenue rate for the system is given by,
\begin{align*}
\lim_{t \rightarrow \infty} \frac{\R(t)}{t} &= \lim_{t \rightarrow \infty} \sum_{s: n>0}\left(\frac{ \R(s,t)}{t}\right)
\stackrel{(a)}{=}  \sum_{s: n>0} \lim_{t \rightarrow \infty} \left(\frac{ \R(s,t)}{t}\right)
\stackrel{(b)}{=} \sum_{s: n>0}  \lambda f(\phi) \phi  \pi(s) \text{ a.s.}
\end{align*}
where the equalities~$(a)$ and~$(b)$ follows from Lemma~\ref{lem_exchange_of_limits_mr} and~\eqref{eqn_revenue_using_rrt} respectively.
Observe 
  from \eqref{eqn_pi0_expression}
  that $\sum_{s: n>0} \pi(s)  = 1 - \DA_i $.
\end{proof}

\begin{lemma} \label{lem_exchange_of_limits_mr}
The long-run revenue rate of the system is equivalent to the long-run revenue rate of the system in state~$s$ when summed over all possible states, i.e.,
$$\lim_{t \rightarrow \infty} \sum_s \frac{\R(s,t)}{t} =  \sum_s \lim_{t \rightarrow \infty} \frac{\R(s,t)}{t} \ \ \text{a.s.}$$
\end{lemma}
\begin{proof}[Proof of Lemma~\ref{lem_exchange_of_limits_mr}]
Consider the following:
\begin{equation}
\label{eq_epsilon}
 \epsilon_t :=  \left| \sum_s \frac{\R(s,t)}{t} - \sum_s \R^*(s) \right|, \text{ where } \R^*(s) := \lim_{t \rightarrow \infty} \frac{\R(s,t)}{t}  \ \stackrel{\mbox{ \eqref{eqn_revenue_using_rrt} }}{=}  \lambda f(\phi) \phi  \pi(s)  \text{ a.s.} 
\end{equation}
To prove the required result, one needs to show that~$\epsilon_t \to 0$ as~$t \to \infty$ a.s.
Next, we define the following:
\begin{align*}
C_b &:= \{s: n+r \ge b\} \text{ for every }b \in \mathds{Z}^+ , \text{ and } \\
\mathcal{A}_b &:= \left \{\omega : \frac{V_t(C_b)(\omega)}{t} \rightarrow \pi(C_b) \right \} \cap \left \{\frac{\R(s,t)}{t} \rightarrow \R^*(s)  \ \forall \ s \right \},
\end{align*}
where,~$V_t(C_b)$ is the number of visits to the set~$C_b$ till time~$t$. By Law of Large Numbers applied to Markov chains as in standard textbooks (see \cite[Theorem 1.10.2]{norris}), we have that as $t \to \infty$, time-average ($\nicefrac{V_t(C_b)}{t}$) converges to stationary measure  ($\pi(C_b)$) almost surely and hence the first set has probability~$1$. Further, by RRT applied to each $s$ (as shown in~\eqref{eqn_revenue_using_rrt} of Lemma~\ref{lem_mr_derivation}), and then using the continuity of probability measure~\cite[Theorem 2.4]{protter} (as set of states $\{s\}$ is countable), we have $\Prob(\mathcal{A}_b) = 1$ for all~$b$.

Let~$\mathcal{A} := \cap_b \mathcal{A}_b$. In other words, we have
\begin{equation*}
\mathcal{A} = \left \{\omega : \frac{V_t(C_b)(\omega)}{t} \rightarrow \pi(C_b) \ \forall \ b \right \} \cap \left \{\frac{\R(s,t)}{t} \rightarrow \R^*(s)  \ \forall \ s \right \}.
\end{equation*}
Again
using continuity of probability, we get that~$\Prob(\mathcal{A}) = 1$.

Now, the final step is to show almost surely that for all~$\tilde{\delta} > 0$, there exists a~${T}_{\tilde \delta}$ (can depend upon sample path) such that~$\epsilon_t \le \tilde{\delta}$ when $t \ge {T}_{\tilde \delta}$. Towards this, 
 set $\delta = \nicefrac{\tilde \delta}{ (\lambda f(\phi) \phi)} $
and  first choose a~$b_\delta >0$ such that the stationary distribution, 
$$
\pi( C_{b_\delta}) =   \sum_{s\in  C_{b_\delta}} \pi(s)  = \sum_{s: n+r > b_\delta} \pi(s) < \frac{\delta}{4}.
$$
This is possible because the stationary distribution exists (see Lemma~\ref{lem_BCMP_Product_form}) and is summable. Now, for any~$\omega \in \mathcal{A}$, choose a~${T}_{\tilde \delta}$ such that for all~$t \ge {T}_{\tilde \delta}$,
\begin{equation}
 \frac{V_t(C_{b_\delta})}{t}(\omega) \le \pi( C_{b_\delta}) + \frac{\delta}{4}.  
 \label{eq_bound}
\end{equation}

From~\eqref{eq_epsilon}, we have the following,\footnote{\label{footnote_rrt_proof}Observe here that for any $t$, we have,  $\sum_{s \in C_{b_\delta}}\frac{\R(s,t)}{t} = \lambda f(\phi) \phi \frac{V_t(C_{b_\delta})}{t}$.}
\begin{align}
\epsilon_t & \leq \left| \sum_{s\in C_{b_\delta}}\left( \frac{\R(s,t)}{t} - \R^*(s) \right)\right| + \left| \sum_{s \in C_{b_\delta}^c}\left( \frac{\R(s,t)}{t} - \R^*(s) \right)\right| \nonumber \\ 
&
\leq \left| \sum_{s\in C_{b_\delta}}\left( \frac{\R(s,t)}{t} - \R^*(s) \right)\right| + \sum_{s \in C_{b_\delta}^c} \left| \left( \frac{\R(s,t)}{t} - \R^*(s) \right)\right| \nonumber \\ & \leq \lambda f(\phi)\phi  \left|  \frac{V_t(C_{b_\delta})}{t} -\pi ( C_{b_\delta} )  \right| + 
\sum_{s \in C_{b_\delta}^c} \left| \left( \frac{\R(s,t)}{t} - \R^*(s) \right)\right|;
\label{eq_epsilon_reduced}
\end{align}
note that the first inequality follows from the triangle's inequality; the second step is again a result of the application of the triangle's inequality and the fact that~$C_{b_\delta}^c$ is a finite set; the third step follows from~\eqref{eq_epsilon} and footnote~\footref{footnote_rrt_proof}.

 Next, for the chosen~$\omega$  pick~${T}_{\tilde \delta}$ further large (if required), such that (possible as $|C_{b_\delta}^c| < \infty$)
 $$\left| \frac{\R(s,t)}{t}(\omega) - \R^*(s) \right| < \frac{\delta}{2 \lvert C_{b_\delta}^c\rvert} \mbox{ for all } t \ge {T}_{\tilde \delta},   \text{  uniformly across } s \in C_{b_\delta}^c.$$ 
Thus, using~\eqref{eq_bound} and the above argument,~\eqref{eq_epsilon_reduced} can be further simplified to the following,
\begin{equation}
  \epsilon_t  \leq  \lambda \phi f(\phi)\phi  \left( \frac{\delta}{2} + \frac{\delta}{2} \right) = \lambda  f(\phi)\phi \delta = \tilde \delta,  \mbox{ for all } t \ge {T}_{\tilde \delta}. \nonumber 
\end{equation}
The above is true for all $\omega \in {\cal A}$, recall $\Prob(A) = 1$,
and thus the result. \end{proof}

\begin{proof}[Proof of Lemma~\ref{lem_approx_MR}]
For any static price policy~$\phi_i$, and $\beta > 0$ from Lemma~\ref{lem_mr_derivation} we have:
\begin{eqnarray}
\label{Eqn_MR_static}
\MR_i (\phi_i; \beta)  =    
  \lambda_i f(\phi_i) \phi_i \sum_{s_i: \niw \neq 0} {\pi}_i(s_i) 
 = 
\lambda_i f(\phi_i) \phi_i (1-\DA_i ),
\end{eqnarray}
where $\DA_i$ is defined in \eqref{eqn_pi0_expression}. 
The proof follows by Lemma~\ref{lem_PB_joint_continuity_from_scratch} 
and the convergence of WEs given in the hypothesis (observe substituting $\DA_i(\phi_i, \lambda_i (\phi_i; 0); 0)$ in \eqref{Eqn_MR_static} leads to~\eqref{eq_approx_MR}).
\end{proof}

\begin{lemma}\label{lem_PB_joint_continuity_from_scratch}
Define  $z(\phi, \lambda; \beta) = \nicefrac{1}{ \sum_{n=0}^{\infty} \mu_n(\phi, \lambda; \beta)}$, where $\mu_n(\phi, \lambda; \beta) := \prod_{a=1}^n \frac{e}{\lambda f(\phi) + a\beta}$, when the series converges, else define $z(\phi, \lambda; \beta) =  0$. 
Then the functions $z,$  $\DA_i$ of \eqref{eqn_pi0_expression}, $\PB_i$ of~\eqref{eqn_expression_of_PBi} and $\WB_i$ of \eqref{eqn_wait_metric_defn}  for any~$i\in \mathcal{N}$ are jointly continuous over $[0, \phi_h] \times [0, \Lambda] \times [0, \infty)$.
Further, at $\beta = 0$, we have:
$$\DA_i(\phi, \lambda; 0) = \begin{cases}
     1 - \frac{e}{\lambda f(\phi)}, & \mbox{if the   series  $\sum_n \mu_n (\phi, \lambda; 0)$ converges,} \\ 
     0, & \text{else.}
\end{cases}$$
\ignore{\begin{enumerate}[$(1)$]
    \item the function $z$ is continuous, i.e., $z(\phi_{i,m}, \lambda_{i,m}; \beta_m) \rightarrow z(\phi_i, \lambda_{i}; \beta)$, and 
    \item $\PB_i(\phi_{i,m}, \lambda_{i,m}; \beta_m) \rightarrow \PB_i(\phi_i, \lambda_{i}; \beta)$,
\end{enumerate}}
\end{lemma}

\begin{proof}
We begin to prove the continuity of function $z$ and consider a  sequence 
$(\phi_{m}, \lambda_{m}, \beta_m) \to (\phi, \lambda, \beta)$. We first consider the case with  $\beta > 0$, followed by the case with $\beta = 0$.
%
In the latter case, we will further have three sub-cases.

\begin{enumerate}[$(i)$]
    \item {${\beta > 0}$:}
As~$\beta_m \to \beta$, for any~$\epsilon \in( 0, \beta)$, there exists a~$M_\epsilon \in \mathbb{N}$ such that 
$$ 
\beta-\epsilon \le \beta_m \le \beta+ \epsilon \mbox{ and hence } \frac{e}{\beta_m}\leq \frac{e}{(\beta-\epsilon)}, \text{ for all } m \geq M_\epsilon.
$$ 
%
Further,  by {\bf A}.2,~$\lambda_{m}f(\phi_{m}) \ge 0$,  and hence  for all   $m \ge M_\epsilon$,
\begin{align*} 
\sum_n \mu_n (\phi_{m}, \lambda_{m}; \beta_m)  
\leq   \sum_{n=0}^{\infty} \prod_{a=1}^{n} \left( \frac{e}{a \beta_m} \right)
\leq \sum_{n=0}^{\infty}  \left( \frac{e}{ \beta - \epsilon} \right)^n\frac{1}{n!}
= \exp{\left( \frac{e}{\beta- \epsilon}\right)}.
\end{align*}
Clearly,  the term  for each $n$ in the above is upper bounded uniformly for all $m \ge M_\epsilon$ by,
$$ \mu_n (\phi_{m}, \lambda_{m}; \beta_m) \le \left( \frac{e}{ \beta - \epsilon} \right)^n\frac{1}{n!}, $$ and these upper bounds are summable (or integrable w.r.t. counting measure). 
Further, by continuity of the function $f$, as $m \to \infty$, we have  $\mu_n (\phi_{m}, \lambda_{m}; \beta_m)  \to \mu_n (\phi, \lambda; \beta) $
for each~$n$. Thus using the dominated convergence theorem, as applied to the counting measure on non-negative integers $\{n\}$, we have that 
\begin{equation}\label{eqn_MR_cont_B_proof_1}
\frac{1}{z(\phi_{m}, \lambda_{m};\beta_m)} = \sum_n \mu_n (\phi_{m}, \lambda_{m}; \beta_m)  \to \sum_n \mu_n (\phi, \lambda; \beta) = \frac{1}{z(\phi, \lambda;\beta)}.
\end{equation}

\item ${ \beta = 0}$: Define a function $h (\phi, \lambda) = \lambda f(\phi)$.  From {\bf A}.1 the function~$h$ is continuous    and thus~$h(\phi_{m}, \lambda_{m}) \to h(\phi, \lambda) $. Thus  for any~$\epsilon > 0$, there exists a~$M_\epsilon \in \mathbb{N}$ such that
\begin{align}
\label{eqn_bound_on_fun_h}
  h(\phi, \lambda) - \epsilon \leq  h(\phi_{m}, \lambda_{m})  \leq h(\phi, \lambda) +  \epsilon, \text{ for all } m \geq M_\epsilon.
\end{align}
Now, we consider the sub-cases depending on the driver and passenger load ratio.

\textbf{Case~(a)~$e < h(\phi, \lambda)$}:
  Choose~$\epsilon > 0~$ small enough such that~$e < h(\phi, \lambda) - \epsilon$. From~\eqref{eqn_bound_on_fun_h}, using similar arguments as before,
{\small \begin{align*}
\sum_{n=0}^{\infty} \prod_{a=1}^{n} \left( \frac{e}{h(\phi_{m}, \lambda_{m} ) + a\beta_{m}} \right) {\leq }
\sum_{n=0}^{\infty} \left( \frac{e}{h(\phi, \lambda ) - \epsilon} \right)^n = \frac{h(\phi, \lambda) - \epsilon}{h(\phi, \lambda) - \epsilon - e}, \text{ for all } m \geq M_{\epsilon}.
\end{align*}}
Applying the dominated convergence theorem as in part~(i), we obtain the desired continuity \eqref{eqn_MR_cont_B_proof_1}.

  {\bf Case (b)~$e > h(\phi, \lambda)$}:
 Choose~$\epsilon > 0$  small enough such that  $e  >  h(\phi, \lambda) + \epsilon$. From \eqref{eqn_bound_on_fun_h}, 
$  h(\phi_{m}, \lambda_{m})  \le h(\phi, \lambda) +\epsilon < e  $ for all $m \ge M_\epsilon$. 

 Observe that the following series diverges:
$$
\sum_{n=0}^{\infty} \left( \frac{e}{h(\phi, \lambda) + \epsilon} \right)^n  
$$
Consider any $\Delta > 0.$ Then there exists an $N_\Delta$ such that 
\begin{equation} \label{eqn_MR_Cont_PB_proof_2}
 \sum_{n=0}^{N_\Delta} \left( \frac{e}{h(\phi, \lambda) + \epsilon} \right)^n \geq \Delta.    
\end{equation}
Now, choose $\bar \beta$ small enough such that 
$$
 \sum_{n=0}^{N_\Delta} \left( \frac{e}{h(\phi, \lambda) + \epsilon + N_\Delta \bar \beta} \right)^n \geq \frac{\Delta}{2}.
$$
Choose (if required) $M_\epsilon$ further large (to ensure $\beta_m \leq \bar \beta$) such that for all $\beta_m \le \bar \beta $ and $m \ge M_\epsilon$,
\begin{align*}
 \sum_n \mu_n (\phi_{m}, \lambda_{m}; \beta_m) 
  &\geq \sum_{n=0}^{N_\Delta} \prod_{a=1}^{n} \left( \frac{e}{h(\phi_{m}, \lambda_{m} ) + a\beta_{m}} \right)\\
&\geq
\sum_{n=0}^{N_\Delta} \prod_{a=1}^{n} \left( \frac{e}{h(\phi, \lambda) + \epsilon + N_\Delta \bar \beta} \right) \\
&=
\sum_{n=0}^{N_\Delta} \left( \frac{e}{h(\phi, \lambda) + \epsilon + N_\Delta \bar \beta} \right)^n
\geq \frac{\Delta}{2}.    
\end{align*}
\ignore{
Now choose $\bar \beta$ and then choose  (if required) $M_\epsilon$ further large such that $\beta_m \le \bar \beta $ for all $m \ge M_\epsilon$, and such that (possible because the LHS of the term below converges to the LHS of \eqref{eqn_MR_Cont_PB_proof_2} as $\beta_m \to 0$, observe $N_\Delta < \infty$):
$$
\sum_{n=0}^{N_\Delta} \prod_{a=1}^{n} \left( \frac{e}{h(\phi, \lambda) + \epsilon + a\beta_m} \right) \geq \frac{\Delta}{2} \ \text{ for all } m \geq M_\epsilon.
$$
Thus in all, for all $m \geq M_\epsilon$,
$$
 \sum_n \mu_n (\phi_{m}, \lambda_{m}; \beta_m) 
  =  \sum_{n=0}^{\infty} \prod_{a=1}^{n} \left( \frac{e}{h(\phi_{m}, \lambda_{m} ) + a\beta_{m}} \right) {\geq }
\sum_{n=0}^{\infty} \prod_{a=1}^{n} \left( \frac{e}{h(\phi, \lambda) + \epsilon + a\beta_m} \right) \geq \frac{\Delta}{2}.
$$}
Therefore,  
$$
z  (\phi_{m}, \lambda_{m}; \beta_m) = \frac{1}{   \sum_n \mu_n (\phi_{m}, \lambda_{m}; \beta_m) } \le \frac{2}{\Delta} \ \text{ for every } m \geq M_\epsilon.
$$
Since $\Delta$ can be arbitrarily large, $
z  (\phi_{m}, \lambda_{m}; \beta_m) \to 0 = z  (\phi, \lambda; \beta)$ (by definition, $z  (\phi, \lambda; \beta) = 0,$  as the series $\sum_n \mu_n (\phi, \lambda; \beta = 0) $ is divergent when $e  \geq h(\phi, \lambda)$).   


\textbf{Case~(c)} $e = h(\phi, \lambda)$:
Let $\Delta >0$ be any given value,  pick $\epsilon = \epsilon_\Delta > 0$ and then $N_\Delta$ such that the following two are true (the second is possible as $e < h(\phi, \lambda) + \epsilon$):
  $$
   \frac{ h(\phi, \lambda) + \epsilon}{\epsilon}  \ge \Delta  \mbox{ and for that } \epsilon,  \mbox{ }
  \delta_N = \sum_{n = N + 1}^{\infty} \left( \frac{e}{h(\phi, \lambda) + \epsilon} \right)^n  \le  \frac{\Delta }{2}  \mbox{ for all } N \ge N_\Delta
.$$
Finally, pick $M_\Delta$ such that
~$h(\phi_{m}, \lambda_{m}) + N_\Delta \beta_m \leq h(\phi, \lambda) + \epsilon$ for all~$m \ge M_\Delta$  (possible as~$\beta_m \to 0$) . Then   for  any~$m \ge M_\Delta$,
\begin{eqnarray*}
\sum_n \mu_n (\phi_{m}, \lambda_{m}; \beta_m)
&>& 
\sum_{n=0}^{N_\Delta} \prod_{a=1}^{n} \left( \frac{e}{h(\phi_{m}, \lambda_{m}) + a\beta_m} \right)
\geq \sum_{n=0}^{N_\Delta} \left( \frac{e}{h(\phi, \lambda) + \epsilon} \right)^n  \nonumber
\\
&=&  \sum_{n=0}^{\infty} \left( \frac{e}{h(\phi, \lambda) + \epsilon} \right)^n -  \delta_{N_\Delta} 
= \frac{h(\phi, \lambda) + \epsilon}{h(\phi, \lambda) + \epsilon - e}  - \delta_{N_\Delta}
\ge \frac{\Delta}{2}.
%
\end{eqnarray*}
Since this is true for every $\Delta  > 0$, the proof goes through as in previous sub-case~(b).

\end{enumerate}
Therefore, in all the cases, $z(\phi_{m}, \lambda_{m};\beta_m) \to z(\phi, \lambda;\beta)$. The second part also follows as from~\eqref{eqn_pi0_expression}--\eqref{eqn_expression_of_PBi}, as we have:
\begin{equation*}
\DA_i(\phi_{i}, \lambda_{i}; \beta) = z(\phi_{i}, \lambda_{i}; \beta) \mbox{ and }
\PB_i(\phi_i, \lambda_i; \beta) = \DA_i  +(1-f(\phi_{i}))  (1-\DA_i).
\end{equation*}
Similarly, using \eqref{eqn_pi0_expression}, one can easily calculate~$\DA_i(\phi, \lambda; 0)$.

Observe that $\WB_i$ of \eqref{eqn_wait_metric_defn} can be rewritten a
$$\WB_i = \PB_i + \alpha f(\phi_i) \ \DA_i \sum_{n = 1}^{\bar N} \frac{\mu_{i,n}}{n}, $$
and thus, the joint continuity of $\WB_i$  follows (note that $\bar N$ is finite).
\end{proof}

\section{Proofs related to Section~\ref{sec:monopoly}}
\label{sec:appendix_mono}
\begin{proof}[Proof of Theorem~\ref{thm_monopoly_opt}]
We begin with the case when~$\dm(\phi_h) \ge 0.$
Differentiating the~$\MR(\phi;0)$ expression (see~\eqref{eq:MR_monopoly}) with respect to~$\phi$ (note that the derivative is not defined at~$\phi = f^{-1}(2\rho)$),
\begin{equation} \label{eqn_mr_deri_monopoly}
\frac{\partial}{\partial \phi} \left[\MR(\phi; 0)\right] = e \mathds{1}_{\left\{ \phi < f^{-1}\left(2\rho\right) \right\}} + \frac{\Lambda}{2} \dm(\phi) \mathds{1}_{\left\{ \phi > f^{-1}\left(2\rho\right) \right\}}.
\end{equation}
From the above, it is clear that the derivative of~$\MR$ is positive for all~$\phi \in [0, \phi_h)$ (since from Lemma~\ref{lem_monotonicity_of_h_1_function},~$\dm$ is strictly decreasing). Further,~$\MR$ is continuous for all~$\phi \in [0, \phi_h]$ and hence,~$\phi_h$ is the unique optimizer. Observe that in this case, for~$\rho \leq f(\phi_m)/2,$ we have~$f^{-1}(2\rho) = \phi_h;$ and for~$\rho > f(\phi_m)/2,$ we have~$\phi_m = \phi_h.$

Next, we consider the case when~$\dm(\phi_h) < 0.$
We begin by showing that zero of~$\dm$ exists in this case. Observe that~$\dm$ is a continuous and strictly decreasing function with $\dm(0)> 0$ and~$\dm(\phi_h)<0$. Using the Intermediate Value Theorem (IVT), there exists a unique zero of~$\dm$,~$\phi_m \in (0,\phi_h)$. The proof follows in two sub-cases:

\begin{enumerate}[$(1)$]
    \item If~$\rho \leq f(\phi_m)/2$

\noindent From the hypothesis of this case, we have~$\phi_m \le f^{-1}\left(2\rho\right)$ which implies~$\dm(\phi_m) \ge \dm(f^{-1}\left(2\rho\right)),$ and hence~$\dm(f^{-1}\left(2\rho\right)) \le 0$. From~\eqref{eqn_mr_deri_monopoly},~$\MR$ is strictly increasing for~$\phi \in [0, f^{-1}\left(2\rho\right))$ and strictly decreasing for~$\phi \in (f^{-1}\left(2\rho\right), \phi_h]$. Using continuity of~$\MR$,~$f^{-1}\left(2\rho\right)$ is the unique optimizer in this case.

\item If~$\rho > f(\phi_m)/2$

\noindent Again from the hypothesis of this case,~$\phi_m > f^{-1}\left(2\rho\right)$. Further using~\eqref{eqn_mr_deri_monopoly},~$\frac{\partial}{\partial \phi} \MR(\phi_m; 0) = 0$ since~$\phi_m$ is the zero of function~$\dm$. Using~\eqref{eqn_mr_deri_monopoly} and continuity of~$\MR$,~$\MR$ is strictly increasing for~$\phi \in [0, \phi_m),$ and strictly decreasing for~$\phi \in (\phi_m, \phi_h].$ Hence, the result follows from the continuity of~$\MR$.
\end{enumerate}
As~$f,$ and~$\dm$ are strictly decreasing functions, it is easy to observe that the optimal monopoly price is a non-increasing function of~$\rho.$

\textbf{Convergence of optimisers:} Consider a sequence~$\beta_n \to 0$. We first show that there exists an optimiser for every~$\beta_n$, say~$\phi_n^*$. The revenue rate of the system for any~$(\phi ;\beta)$, can be expressed as
$$
\MR(\phi;\beta)   \ = \   \frac{\Lambda}{2}f(\phi) \phi(1-\DA(\phi; \beta)).
$$
From the proof of Lemma~\ref{lem_PB_joint_continuity_from_scratch}, 
the above function is continuous in~$\phi$ over~$[0, \phi_h]$; thus there exists a~$\phi_n^*$ for every~$\beta_n$. 
An optimal price~$\phi^*$ for any~$\beta \in [0, \infty)$ maximizes the following,
$$
\phi^* \in 
 \Phi^*(\beta) := \text{Arg}\max_{\phi \in [0,\phi_h]} \MR(\phi;\beta).
$$
Clearly, the domain of optimization~$[0,\phi_h]$ is compact and is the same for all~$\beta$, and the function~$\MR(\phi; \beta)$ is jointly continuous in~$(\phi, \beta)$; hence, the hypothesis of the Maximum Theorem \cite[Theorem 9.14]{maximum} is satisfied. Thus the correspondence~$\beta \mapsto \Phi^*(\beta)$ is compact and upper semi-continuous. Consider one optimizer from~$\Phi^*(\beta)$ for each~$\beta_n$, that is~$\phi_n^*$. By \cite[Proposition 9.8]{maximum}, there exists  
 a sub-sequence of~$\{\phi^*_n\}$,  which converges to the unique optimal price of the \idp\ regime.
\end{proof}

\begin{lemma} \label{lem_monotonicity_of_h_1_function}
Under \textbf{A}.1, for any~$a,b > 0,$ the function~$\phi \mapsto a f(\phi) + b \phi f'(\phi)$ is strictly decreasing.
\end{lemma}

\begin{proof}[Proof of Lemma~\ref{lem_monotonicity_of_h_1_function}]
Under \textbf{A}.1,~$f(\phi)$ is a strictly concave and decreasing function. From strict concavity of~$f$, we have that~$f'(\phi)$ is negative and decreasing. Thus,~$f'(\phi)\phi$ is a decreasing function of~$\phi$. Since, sum of two monotonically decreasing functions is again a monotonically decreasing function, we have the result.
\end{proof}
 
\ignore{
\begin{proof}[Proof of Lemma~\ref{lem_pio_WE}]
From Lemma~\ref{lem_exist_unique_WE}, we know that there exist an unique WE with~$\DA$ as a WE metric. Thus, it suffices to check that at the proposed split~$(\lamo, \lamt),$~$\DA_1 = \DA_2;$ this is straightforward.
\end{proof}

\begin{proof}[Proof of Theorem~\ref{thm_uniqueness_symm_NE}]
Observe that~$\MR_i$ given in~\ref{eqn_MR_d_system_beta_zero} is continuous, thus by doing some algebra, the derivative of~$\MR_i$ in \idp\ regime, with respect to~$\phi_{i}$ for fixed~$\phi_{-i}$ is,
\begin{equation} \label{eqn_mr_der_lem}
{\Delta}_i (\phi_i,\phi_{-i}; 0) = \mathds{1}_{\{e < h(\VPhi)\}} e  + \mathds{1}_{\{e > h(\VPhi)\}} \frac{\Lambda f(\phi_{-i})}{(f(\phi_i) + f(\phi_{-i}))^2} \left[ f(\phi_i) (f(\phi_i) - f(\phi_{-i})) + f(\phi_{-i}) ( \du(\phi_{i}) )
\right].
\end{equation}
Note that the partial derivative given in~\eqref{eqn_mr_der_lem} is not defined when~$e = h(\VPhi).$
When~$\du(\phi_h) \geq 0,$ observe that~$\phi_u = \phi_h;$ and for~$\du(\phi_h) < 0,$ the unique~$\phi_u < \phi_h.$ Now we prove the existence and uniqueness of the NE with the two cases when~$\du(\phi_h) \geq 0,$ and~$\du(\phi_h) < 0$ using~\eqref{eqn_MR_d_system_beta_zero} and~\eqref{eqn_mr_der_lem}.

We begin with the case~$\du(\phi_h) \geq 0.$ Observe that when $\rho \leq f(\phi_h)/2,$ from~\eqref{eqn_f_inverse},~$f^{-1}\left(2\rho\right) = \phi_h$; and when~$\rho > f(\phi_h)/2,$~$\phi_u = \phi_h$. Thus, in this case, proving~$(\phi_h, \phi_h)$ is NE suffices. Fix~$\phi_{-i} = \phi_h$. Substituting this in~\eqref{eqn_mr_der_lem}, we have,
\begin{equation} 
 {\Delta}_i (\phi_i,\phi_h; 0) = \mathds{1}_{\{e < h(\VPhi)\}} e  + \mathds{1}_{\{e > h(\VPhi)\}} \frac{\Lambda f(\phi_h)}{(f(\phi_i) + f(\phi_h))^2} \left[ f(\phi_i) (f(\phi_i) - f(\phi_h)) + f(\phi_h) ( \du(\phi_{i}) )
\right] \nonumber
\end{equation}
Recall that~$f,$ and~$\du$ are strictly decreasing functions (Lemma~\ref{lem_monotonicity_of_h_1_function}), and~$\du(\phi_h) \ge 0$; thus from above,~${\Delta}_i (\phi_i,\phi_h; 0) > 0$ for all~$\phi_i \in [0, \phi_h)$ (excluding the single point of non-differentiability mentioned above). From the continuity of~$\MR_i$,~$\phi_h$ is the best response of platform~$i$ against~$\phi_{-i} = \phi_h$. Thus~$(\phi_h, \phi_h)$ is a NE.

To prove uniqueness, we use a contradiction based argument. Say there exists a~$\tilde \phi < \phi_h$ such that~$(\tilde \phi, \tilde \phi)$ is also a NE. Then using~\eqref{eqn_mr_der_lem}, the partial derivative at~${\VPhi = (\tilde \phi, \tilde \phi)}$ is,
\begin{equation}
{\Delta}_i (\tilde \phi, \tilde \phi; 0) = \mathds{1}_{\{e < h(\VPhi)\}} e  + \mathds{1}_{\{e > h(\VPhi)\}} \frac{\Lambda }{4} \du(\tilde \phi) > 0, \nonumber
\end{equation}
since~$\du$ is a strictly decreasing function and~$\du(\phi_h) \ge 0$. Thus, platform~$i$ has the incentive to increase the price from~$\phi_i = \tilde \phi$. This contradicts the assumption that~$\tilde \phi$ is also a NE.
Therefore~$(\phi_h, \phi_h)~$ is the unique NE.

We next consider the case when~$\du(\phi_h) < 0.$ The proof of existence of unique zero follows as in the proof of Theorem~\ref{thm_monopoly_opt}. The remaining proof follows in two cases:

\textbf{Case 1:} When~$\rho \leq f(\phi_u)/2$

Fix~$\phi_{-i} = f^{-1}\left(2\rho\right)$. Substituting this in~\eqref{eqn_mr_der_lem} we have,
\begin{equation} \label{eqn_mr_der_lem_phiunderbar}
 {\Delta}_i (\phi_i,f^{-1}\left(2\rho\right); 0) = \mathds{1}_{\{e < h(\VPhi)\}} e  + \mathds{1}_{\{e > h(\VPhi)\}} \frac{2\Lambda \rho}{(f(\phi_i) + 2\rho)^2} \left[ f(\phi_i) (f(\phi_i) - 2\rho) + 2\rho ( \du(\phi_{i}) )
\right].
\end{equation}
Using simple calculations, one can easily verify that if~$\phi_i < f^{-1}\left(2\rho\right),$ then~$e < h(\VPhi);$ and thus~${\Delta}_i (\phi_i,f^{-1}\left(2\rho\right); 0) > 0$. On the other hand, if~$\phi_i > f^{-1}\left(2\rho\right)$, then~$e > h(\VPhi);$ and hence~$f(\phi_i) < 2\rho \leq f(\phi_u)$ (since~$\rho \leq f(\phi_u)/2$).  This implies~$f(\phi_i)-2\rho <0$. Additionally,~$\phi_u < \phi_i$ and hence,~$\du(\phi_i) < 0$ (since~$\du$ is a strictly decreasing function). Thus from~\eqref{eqn_mr_der_lem_phiunderbar}, for~$\phi_i > f^{-1}\left(2\rho\right)$, we have~${\Delta}_i (\phi_i, f^{-1}\left(2\rho\right); 0)<0$. Hence, using continuity of~$\MR_i$,~$f^{-1}\left(2\rho\right)~$ is the best response of platform~$i$ against~$\phi_{-i} = f^{-1}\left(2\rho\right)$, and hence,~$(f^{-1}\left(2\rho\right), f^{-1}\left(2\rho\right))$ is a NE.

To prove the uniqueness, we use a contradiction based argument. Suppose there exists a $\tilde \phi \neq f^{-1}\left(2\rho\right)$ such that~$(\tilde \phi, \tilde \phi)$ is also a NE. Fix~$\phi_{-i} = \tilde \phi$. Now, we have two sub-cases.

(a)~If~$\tilde \phi >  f^{-1}\left(2\rho\right)$ and~$\phi_i = \tilde \phi$ then~$e > h(\VPhi)$ and thus using~\eqref{eqn_mr_der_lem},
${\Delta}_i (\tilde \phi,\tilde \phi; 0) = (\Lambda/4) \du(\tilde \phi).$
Since~$\du(\tilde \phi) < 0$ and the function~$\du$ is strictly decreasing, platform~$i$ has the incentive to reduce its price, which is a contradiction to our assumption.

(b) On the other hand, if~$\tilde \phi <  f^{-1}\left(2\rho\right)$ and~$\phi_i = \tilde \phi$ then~$e < h(\VPhi)$ and thus~${\Delta}_i (\tilde \phi,\tilde \phi; 0) = e$. Hence, platform~$i$ has an incentive to increase its price, which is again a contradiction.

Thus~$ ( f^{-1}\left(2\rho\right),  f^{-1}\left(2\rho\right))$ is the unique NE.

\textbf{Case 2:} When~$\rho > f(\phi_u)/2$

Fix~$\phi_{-i} = \phi_u$. Substituting this in~\eqref{eqn_mr_der_lem}, we obtain,
\begin{equation} \label{eqn_mr_der_phi_d}
{\Delta}_i (\phi_i,\phi_u; 0) = \mathds{1}_{\{e < h(\VPhi)\}} e  + \mathds{1}_{\{e > h(\VPhi)\}} \frac{\Lambda f(\phi_u)}{(f(\phi_i) + f(\phi_u))^2} \left[ f(\phi_i) (f(\phi_i) - f(\phi_u)) + f(\phi_u) ( \du(\phi_{i}) )
\right].
\end{equation}
Observe that if~$\phi_i < \phi_u,$ then~$f(\phi_i)> f(\phi_u),$ and~$\du(\phi_i) > 0;$ thus~${\Delta}_i (\phi_i,\phi_u; 0) > 0.$ If~$\phi_i > \phi_u,$ then~$e > h(\VPhi)$ (as~$\phi_u >  f^{-1}\left(2\rho\right)$), and~$\du(\phi_i) < 0$ (see~\eqref{eqn_mr_der_phi_d}); thus~${\Delta}_i (\phi_i,\phi_u; 0) < 0.$ Thus,~$\phi_u$ is the best response of platform~$i$ against~$\phi_{-i} = \phi_u$, and hence~$(\phi_u, \phi_u)$ is a NE.

To prove the uniqueness, we use a contradiction based argument.  Say there exists a~$\tilde \phi \neq \phi_u$ such that~$(\tilde \phi, \tilde \phi)$ is also a NE. Fix~$\phi_{-i} = \tilde \phi$. Now, we have two sub-cases:

(a) If~$\tilde \phi > \phi_u$ and~$\phi_i = \tilde \phi$ then~$e > h(\VPhi)$  and thus~\eqref{eqn_mr_der_lem} can be simplified to~${\Delta}_i (\tilde \phi,\tilde \phi; 0) = (\Lambda/4) \du(\tilde \phi)$. As~$\tilde \phi > \phi_u$ and~$\du$ is a strictly decreasing function,~$\du(\tilde \phi) < \du(\phi_u) = 0$ . Thus, platform~$i$ has the incentive to reduce its price, which is a contradiction.

(b) If~$\tilde \phi < \phi_u$ and~$\phi_i = \tilde \phi$ then using~\eqref{eqn_mr_der_lem},  
$$
{\Delta}_i (\tilde \phi,\tilde \phi, 0)  = \mathds{1}_{e<h(\VPhi)}e + \mathds{1}_{e>h(\VPhi)} \frac{\Lambda}{4} \du(\tilde \phi) > 0,
$$
as~$\du(\tilde \phi) > 0$ (since~$\du(\tilde \phi) > \du(\phi_u) = 0$). Therefore platform~$i$ has the incentive to increase its price, which is a contradiction.

Thus~$ (\phi_u, \phi_u)$ is the unique NE.
\end{proof}

\begin{proof}[Proof of Theorem~\ref{Thm_beta_nonzero_Dsystem}]
From Lemma~\ref{lem_MR_joint_cont}, the function~$\MR_i$ is continuous on~$[0, \phi_h]^2\times [0, \infty)$, and thus it is uniformly continuous on the compact interval~$[0, \phi_h]^2\times [0, \bar \beta]$ for some~$\bar \beta > 0$.
Therefore for any~$\epsilon > 0$, there exists~$ \ \bar \beta_{\epsilon} > 0$ such that for all~$ \beta \leq \bar \beta_{\epsilon}$, 
\begin{equation}\label{eqn_mr_D_sys_MR_close_by_epsilon}
\text{for all } (\VPhi, \beta) \in \mathcal{F}_\epsilon, \ | \MR_i(\phi_i, \phi_{-i}; 0) - \MR_i(\phi_i, \phi_{-i}; \beta) | < \epsilon.
\end{equation}
Using~\eqref{eqn_mr_D_sys_MR_close_by_epsilon}, and the definition~\ref{def_epsilon_equilibria}, we have the desired result.
\end{proof}}

\section{Proofs related to Section~\ref{sec:overall_bp}}
\label{sec:appendix_duopoly_PB}
\begin{proof}[Proof of Theorem~\ref{thm_WE_MR_PB}]
We begin by showing the joint continuity of the objective function given in~\eqref{eqn_existence_uniq_WE_main} with respect to~$(\lambda_1, \lambda_2, \beta)$ when~$\VPhi$ is fixed. To this end, it suffices to show that~$\PB_i(\lambda_i, \beta)$ is jointly continuous over~$[0, \Lambda] \times [0, \infty).$
For any~$\lambda_{i} \in [0, \Lambda]$ and~$\beta \in [0, \infty)$, consider~$(\lambda_{i,m}, \beta_m) \to (\lambda_i, \beta)$ as~$m \to \infty.$ From~\eqref{Eq_overall_exact} and~\eqref{eqn_pi0_expression}, we have
\begin{align*}
\PB_i(\lambda_{i,m}; \beta_{i,m})  &=  \DA_i(\lambda_{i,m}; \beta_{i,m}) f(\phi_i) +(1-f(\phi_i)) \\    
&= \left( \sum_{n=0}^{\infty} \left[\frac{e^{n}}{\prod_{a=1}^{n}(\lambda_{i,m} f(\phi_i) + a \beta_{i,m})}\right]\right)^{-1}  
f(\phi_i) +(1-f(\phi_i)) .
\end{align*}
Now the required joint continuity is given by  
Lemma~\ref{lem_PB_joint_continuity_from_scratch}.

The above continuity and the compactness of the domain $[0, \Lambda]$ establishes the existence of WE defined using \eqref{lem_exist_unique_WE}. Further using the above joint continuity,  one can apply  Maximum theorem~\cite{maximum} to the objective function \eqref{lem_exist_unique_WE} and obtain that the set of optimizers (or the WEs) is upper semi-continuous in $\beta \in [0, \infty)$ (here it is important to observe the domain of optimization is $[0, \Lambda]$ for all $\beta$
and hence is a continuous and compact correspondence w.r.t. the parameter $\beta$). 
By Lemma \ref{lem_exist_unique_WE}, we also have uniqueness of WE when $\beta > 0$ and hence the function $\beta \mapsto (\lambda_1(\VPhi; \beta), \lambda_2(\VPhi; \beta))$ (the unique WE as a function of $\beta$) is continuous at all $\beta > 0$ (see \cite[Theorem 9.12]{maximum}).

Thus $\lim_{\beta \to 0} (\lambda_1(\VPhi; \beta), \lambda_2(\VPhi; \beta))$ exists, and the rest of the proof computes these limits and shows them equal to the ones reported in Table \ref{tab:WE_MR_PB}. By upper semi-continuity of the set of WEs, these limits are WEs at $\beta = 0$  (proved via \cite[Proposition 9.8]{maximum} and because limit via a sub-sequence equals that via the entire sequence). Thus, we derive the limits of Table \ref{tab:WE_MR_PB}, either by establishing the uniqueness of WE at $\beta = 0$ (in some cases) or by finding the limits of  WEs as $\beta \to 0$ (in the rest of the cases).  
This part is established in different cases with $\beta = 0$ as below.
Once the limits in Table~\ref{tab:WE_MR_PB} are computed, we have continuity of (the chosen) WE $(\lambda_1(\VPhi; \beta), \lambda_2(\VPhi; \beta))$ at all $\beta \geq 0$. Towards the end, we also establish the continuity of the $\MR_i$ for all $\beta \ge 0$ for each $i.$



The limiting value of the blocking probability, as~$\beta \to 0$ is given by (using \eqref{Eq_overall_exact} and Lemma~\ref{lem_PB_joint_continuity_from_scratch}),
\begin{equation}
    \label{Eq_overall_aproxx}
   \PB_i(\VPhi; \lambda_i(\VPhi, 0))  = 
        \begin{cases}
            1 - \frac{\eiw}{\lambda_i} & \text{ if }  \left(\frac{\eiw}{\lambda_i f(\phi_i)}\right) < 1, \\
            1 - f(\phi_i) & \text{ else.}
        \end{cases}
\end{equation}

    \noindent \textit{Case~$(1)$}:~$\phi_1 = \phi_2$. From Lemma~\ref{lem_exist_unique_WE}, and~\eqref{Eq_overall_exact}, it is easy to see that the unique WE for~$\beta > 0$ system is~$(\lamo, \lamt) =(\nicefrac{\Lambda}{2}, \nicefrac{\Lambda}{2}).$  Thus, we set $\lamol = \lamtl = \nicefrac{\Lambda}{2}$ in Table \ref{tab:WE_MR_PB}.


\noindent \textit{Case~$(2)$}:~$\phi_1 > \phi_2$

\noindent \textit{Case $(2a)$}:~$\phi_1, \phi_2 \in [0, f^{-1}\left( 2\rho \right))$. From the definition of~$f^{-1}\left( 2\rho \right)$, we have~$\nicefrac{2e}{\Lambda f(\phi_i)} < 1$ for any~$i \in \{1,2\}$. Then, from~\eqref{Eq_overall_aproxx}, it is easy to see that~$(\nicefrac{\Lambda}{2}, \nicefrac{\Lambda}{2})$ is a WE at $\beta = 0$. We now argue that this is the unique WE at $\beta = 0$. Indeed, say there exists another WE with~$\epsilon>0$, such that
$$
(\lambda'_1(\VPhi; 0),\lambda'_2(\VPhi; 0)) = \left( \frac{\Lambda}{2}-\epsilon,\frac{\Lambda}{2}+\epsilon \right).
$$
\begin{enumerate}[]
\item If~$\PB_1'(\VPhi; 0) = 1-f(\phi_1)$ and~$\PB_2'(\VPhi; 0) = 1 -f(\phi_2)$:
 These blocking probabilities can never be equal  since~$f(\phi_1)\neq f(\phi_2)$.
 \item If~$\PB_1'(\VPhi; 0) = 1-\nicefrac{e}{\lambda'_1(\VPhi; 0)}$ and~$\PB_2'(\VPhi; 0) = 1-\nicefrac{e}{\lambda'_2(\VPhi; 0)}$: These cannot be equal since~$\lambda'_1(\VPhi; 0)\neq \lambda'_2(\VPhi; 0)$.
 \item If~$\PB_1'(\VPhi; 0) = 1 - f(\phi_1)$ and~$\PB_2'(\VPhi; 0) = 1 - \nicefrac{e}{\lambda'_2(\VPhi; 0)}$: This case is also not possible since,
$$ \frac{e}{\lambda'_2(\VPhi; 0)} = \frac{e}{\nicefrac{\Lambda}{2} + \epsilon} < \frac{e}{\nicefrac{\Lambda}{2}} = 2\rho  < f(\phi_1).$$
\item {If~$\PB_1'(\VPhi; 0) = 1 - \nicefrac{e}{\lambda'_2(\VPhi; 0)}$ and~$\PB_2'(\VPhi; 0) = 1 - f(\phi_2)$: This case is also not possible, as~$\PB_1'(\VPhi; 0) = 1 - \nicefrac{e}{\lambda'_2(\VPhi; 0)}$, which implies that~$\left(\frac{e}{\nicefrac{(\Lambda}{2} - \epsilon) f(\phi_1)}\right) \leq  1$
$$ \frac{e}{\lambda'_1(\VPhi; 0)} = \frac{e}{\nicefrac{\Lambda}{2} - \epsilon}  \leq f(\phi_1) < f(\phi_2).$$}
\end{enumerate}
Thus,~$(\nicefrac{\Lambda}{2}, \nicefrac{\Lambda}{2})$ is the unique WE at $\beta = 0$ in this sub-case and equals the one in the table.

\noindent \textit{Case $(2b)$}:~$\phi_1 \in [f^{-1}\left( 2\rho \right), f^{-1}\left( \rho \right))$. In this case a limit WE is given by,
\begin{equation} \label{eq_WE_b}
    \lamol   =  \left(\Lambda  - \frac{e}{f(\phi_1)}\right).
\end{equation} 
 One can verify the above by direct substitution in~\eqref{Eq_overall_aproxx}; the resulting~$\PB_1(\VPhi; 0),$ and~$\PB_2(\VPhi; 0)$ are equal and are given by~$1- f(\phi_1)$ and~$1 - \nicefrac{e}{\lamtl},$ respectively. 
Next, we show that the above WE is unique. Towards this we first study~$\PB_1(\lambda_1; 0)$ and~$\PB_2 (\Lambda-\lambda_1; 0)$ of~\eqref{Eq_overall_aproxx} as a function of~$\lambda_1$.  With~$x_1:=\nicefrac{e}{  f(\phi_1)}$ and~$x_2:=\Lambda- \nicefrac{e}{  f(\phi_2)},$ 
\begin{align*}
    \PB_1(\lambda_1; 0) &= 1- \left(f(\phi_1) \mathds{1}_{ \lambda_1 \le x_1}  + \frac{e}{\lambda_1}   \mathds{1}_{ \lambda_1 > x_1}\right),  \\
     \PB_2(\lambda_1; 0) &=1- \left(f(\phi_2) \mathds{1}_{ \lambda_1 \ge x_2}  + \frac{e}{(\Lambda-\lambda_1)}   \mathds{1}_{ \lambda_1 < x_2}\right).
\end{align*}
Thus from~\eqref{eq_WE_b}, and as~$\phi_1 \geq f^{-1}\left( 2\rho \right),$ we have~$\Lambda - \nicefrac{e}{f(\phi_1)} \leq \nicefrac{e}{f(\phi_1)}$; thus,~$\lamol \le x_1,$ and~$\lamol < x_2$.  
Since~$ \PB_1(\lambda_1; 0)$ is a non-decreasing function, and~$ \PB_2(\lambda_1; 0)$ is a non-increasing function, the later being strictly decreasing at the above identified WE, it follows that there is no other~$\lambda_1$ at which the two~$\PB_i$ become equal.

\noindent \textit{Case $(2c)$}:~$\phi_1 \in [f^{-1}\left( \rho \right), \phi_h]$. From~\eqref{Eq_overall_aproxx} ,~$ \PB_1(\VPhi; 0) = 1 - f(\phi_1)$ irrespective of~$\lamol$.
\begin{enumerate}
    \item If $\phi_2 \in [0, f^{-1}\left( \rho \right))$ then~$\PB_2(\VPhi; 0) = \max\left\{ 1 - \frac{e}{\lamtl}, 1-f(\phi_2) \right \},$ where~$\lamtl$ is unknown. Observe that for any~$\lamtl,$ since~$f(\phi_1) < f(\phi_2),$ $f(\phi_1) \leq \nicefrac{e}{\Lambda}$,
$$
1-f(\phi_1) >1-f(\phi_2) \text{ and } 1-f(\phi_1) \geq  1 - \frac{e}{\lamtl}, \text{ as }  1-f(\phi_1) \geq  1 - \frac{e}{\Lambda}.
$$
In this case, note that the blocking probability of platform~1 is constant; for platform~2, the blocking probability is initially constant value $1 - f(\phi_2)$ till $\nicefrac{e}{f(\phi_2)}$, and then it is an increasing function of~$\lambda_2$, i.e, $1 - \nicefrac{e}{\lambda_2}$. If $\phi_1 = f^{-1}\left( \rho \right),$ then by equating the respective blocking probabilities, we get $(\lamol, \lamtl) = \left(0, \Lambda \right),$ as limit WE, which is the unique optimizer of~\eqref{eqn_existence_uniq_WE_main}). Else, if $\phi_1 > f^{-1}\left( \rho \right)$, it is not possible to equate the blocking probabilities of the two platforms.
Thus, again, we have $(\lamol, \lamtl) = \left(0, \Lambda \right),$ as limit WE, which is the unique optimizer of~\eqref{eqn_existence_uniq_WE_main}.

\item $\phi_2 \in [f^{-1}\left( \rho \right), \phi_h]$, from~\eqref{Eq_overall_aproxx}, $ \PB_2(\VPhi; 0) = 1 - f(\phi_2)$ irrespective of~$\lamtl$. Observe that for any~$\lamtl,$ since~$f(\phi_1) < f(\phi_2),$ we have $1-f(\phi_1) >1-f(\phi_2).$
Thus, it is not possible to equate the blocking probabilities of the two platforms. 
Thus, we consider the following as limit WE, which also optimizes~\eqref{eqn_existence_uniq_WE_main}: 
$$(\lamol, \lamtl) = \left(0, \Lambda \right),$$
and show that the selection ensures the required continuity properties. Towards this, consider any~$\beta_n \to 0$ and consider corresponding~$\lim_{\beta_n \to 0} \lambda_1(\VPhi; \beta_n)$. It suffices to show that~$\lambda_1(\VPhi; \beta_n) \to 0$. 
Note that
$$\PB_1(\phi_1; 0, 0) = 1 - f(\phi_1) > 1 - f(\phi_2) = \PB_2(\phi_2; \Lambda, 0).$$
Consider any~$\epsilon > 0$ with~$2\epsilon < f(\phi_2) - f(\phi_1).$ By Lemma~\ref{lem_PB_joint_continuity_from_scratch}, there exists a~$\beta_\epsilon$ such that for all~$\beta \leq \beta_\epsilon,$ we have~$|\PB_1(\phi_1; 0, \beta) - \PB_1(\phi_1; 0, 0)| < \epsilon,$ and~$|\PB_2(\phi_2; \Lambda, \beta) - \PB_2(\phi_2; \Lambda, 0)| < \epsilon.$ Thus for~$\beta \leq \beta_\epsilon,$ we have
$$
\PB_1(\phi_1; 0, \beta) > \PB_2(\phi_2; \Lambda, \beta).
$$
Note that for~$\beta > 0,$~$\PB_1$ is a strictly increasing function of~$\lambda_1$; while~$\PB_2$ is a strictly decreasing function of~$\lambda_1.$ Therefore~$(0, \Lambda)$ is the unique WE for~$0 < \beta \leq \beta_{\epsilon}.$
\end{enumerate}

\ignore{and show that the selection ensures the required continuity properties. Towards this, consider any~$\beta_n \to 0$ and consider corresponding~$\lim_{\beta_n \to 0} \lambda_1(\VPhi; \beta_n)$. It suffices to show that~$\lambda_1(\VPhi; \beta_n) \to 0$. 
Note that
$$\PB_1(\phi_1; 0, 0) = 1 - f(\phi_1) > 1 - f(\phi_2) = \PB_2(\phi_2; \Lambda, 0).$$
Consider any~$\epsilon > 0$ with~$2\epsilon < f(\phi_2) - f(\phi_1).$ As $\PB_1$ is jointly continuous w.r.t. $(\lambda_1, \beta)$ when $\VPhi$ is fixed, there exists a~$\beta_\epsilon$ such that for all~$\beta \leq \beta_\epsilon,$ we have~$|\PB_1(\phi_1; 0, \beta) - \PB_1(\phi_1; 0, 0)| < \epsilon,$ and~$|\PB_2(\phi_2; \Lambda, \beta) - \PB_2(\phi_2; \Lambda, 0)| < \epsilon.$ Thus for~$\beta \leq \beta_\epsilon,$ we have
$$
\PB_1(\phi_1; 0, \beta) > \PB_2(\phi_2; \Lambda, \beta).
$$
Observe that for~$\beta > 0,$~$\PB_1$ is a strictly increasing function of~$\lambda_1$; while~$\PB_2$ is a strictly decreasing function of~$\lambda_1.$ Therefore~$(0, \Lambda)$ is the unique WE for~$0 < \beta \leq \beta_{\epsilon}.$}

This completes the proof of the claim that the WE is continuous at~$\beta = 0.$

The continuity of~$\MR_i$ with respect to~$\beta$ for~$\beta > 0$ follows from the continuity of the WE, while, the continuity at~$\beta = 0$ follows from the continuity of WE and using Lemma~\ref{lem_approx_MR}. 
\end{proof}

\begin{proof}[Proof of Theorem~\ref{thm_B_sys_limit_sys_main}, and Theorem~\ref{thm_EC}]
Observe that when~$\db(\phi_h) <0$, the existence of unique zero~$\phi_b$ can be shown using IVT since~$\db(0) > 0$, and~$\db$ is a continuous and strictly decreasing function (see Lemma~\ref{lem_monotonicity_of_h_1_function}). The proof heavily depends on~$m(\phi) = \Lambda (f(\phi) - \rho)\phi$ (also used in Table~\ref{tab:WE_MR_PB}); its derivative equals
\begin{equation}
\label{eqn_m_derivative}
    m'(\phi) = \Lambda\left(\db(\phi) + \frac{f(\phi)}{2} - \rho\right).
\end{equation} 
\begin{enumerate}[$(1)$]
    \item {\bf Case~1a:} Consider the case when~$\rho \le f(\phi_h)/2.$ Then from definition~\eqref{eqn_f_inverse},~$f^{-1}\left(2\rho \right) = \phi_h$. W.l.g.,  consider the best response (BR) of player~$i$ against any~$\phi_{-i}$, represented by~$\BR (\phi_{-i})$. Using Table~\ref{tab:WE_MR_PB},~$\MR_i(\phi_i, \phi_{-i}; 0)  = e\phi_i$ and hence 
~$\BR (\phi_{-i}) = \{\phi_h\}$. Thus~$(\phi_h, \phi_h)$ is the unique NE.

   \noindent {\bf Case~1b:} Next, consider~$f(\phi_h)/2 < \rho \le f(\phi_b)/2,$ which implies~$\phi_b < \phi_h.$ 
    Fix~$\phi_{-i} = f^{-1}\left(2 \rho \right)$. 
    Using Table~\ref{tab:WE_MR_PB}, the revenue rate of platform~$i$ against~$f^{-1}\left(2 \rho \right)$ is,
    \begin{equation}\label{eqn_MR_B_phiunderbar_one}
    \MR_i(\phi_i, \phi_{-i}; 0) =
     \begin{cases}
     e\phi_i & \text{ if } \phi_i \in [0,  f^{-1}\left(2 \rho \right) ), \\
     \frac{\Lambda}{2}f(\phi_i)\phi_i & \text{ if } \phi_i = \phi_{-i}, \\
     m(\phi_i) & \text{ if }  \phi_i \in ( f^{-1}\left(2 \rho \right), f^{-1}\left( \rho \right)  ), \\
     0 & \text{ if }  \phi_i \in [ f^{-1}\left( \rho \right), \phi_h ].
    \end{cases}   
    \end{equation}
    Note that~$\MR_i$ is strictly increasing for~$\phi_i < f^{-1}\left(2 \rho \right).$ Moreover, this function is continuous at~$f^{-1}\left(2 \rho \right).$ We now show below that~$\MR_i$ is decreasing for~$\phi_i > f^{-1}\left(2 \rho \right),$ which implies that~$(f^{-1}\left(2 \rho \right), f^{-1}\left(2 \rho \right))$ is a NE.
    
    From the hypothesis,~$\rho \leq f(\phi_b)/2$, which implies~$\phi_b \leq f^{-1}\left(2 \rho \right)$. Since~$\db$ is strictly decreasing  and~$\db(\phi_b) = 0$, we have~$\db(f^{-1}\left(2 \rho \right)) \le 0$, and hence we obtain~$m'(f^{-1}(2\rho)) \le 0$ (see~\eqref{eqn_m_derivative}). Further, by Lemma~\ref{lem_monotonicity_of_h_1_function},~$m'$ is decreasing, and hence we have~$m'(\phi) < 0$ for any~$\phi > f^{-1}(2\rho)$; thus~$m$ is decreasing in that range. 
    
    Now, to prove the uniqueness, we use a contradiction based argument. Say~$(\tilde \phi, \tilde \phi)$ is also a NE. The proof follows in the following two cases.
\begin{enumerate}[$(a)$]
    \item If~$\tilde \phi < f^{-1}\left(2 \rho \right)$

\noindent    The pair~$(\tilde \phi, \tilde \phi)~$ is not a NE because from Table~\ref{tab:WE_MR_PB} and with~$0 < \epsilon < f^{-1}(2\rho) - \tilde \phi$, we have 
$$\MR_i(\tilde \phi, \tilde \phi;0) = e\tilde \phi < e(\tilde \phi+ \epsilon) = \MR_i(\tilde \phi + \epsilon, \tilde \phi;0).$$ 
   \item If~$\tilde \phi > f^{-1}\left(2 \rho \right)$
   
\noindent   Using Table~\ref{tab:WE_MR_PB}, the relevant components of revenue rate of platform~$i$ against~$\tilde \phi$ is,
    \begin{equation*}
    \MR_i(\phi_i, \tilde \phi; 0) =
     \begin{cases}
     e\phi_i & \text{ if } \tilde \phi \in \left(f^{-1}(2\rho),  f^{-1}(\rho)\right) \text{ and } \phi_i < \tilde \phi, \\
     \Lambda f(\phi_i)\phi_i & \text{ if } \tilde \phi >   f^{-1}(\rho) \text{ and } \phi_i < \tilde \phi, \\
     \frac{\Lambda}{2} f(\phi_i) \phi_i & \text{ if }  \phi_i = \tilde \phi.
    \end{cases}   
    \end{equation*}
Observe from Lemma~\ref{lemma_oscillatory} that~$e {\tilde \phi} >  \nicefrac{\Lambda}{2} f({\tilde \phi}) {\tilde \phi}~$ as well as~$\Lambda f(\tilde \phi) {\tilde \phi} >  \nicefrac{\Lambda}{2} f({\tilde \phi}) {\tilde \phi}~$, hence one can choose a~$\phi_i \in (f^{-1}\left(2 \rho \right), {\tilde \phi})~$ sufficiently close to~${\tilde \phi}$ such that 
$$\MR_i(\tilde \phi, \tilde \phi;0) = \nicefrac{\Lambda}{2} f({\tilde \phi}) {\tilde \phi}  < e \phi_i = \MR_i(  \phi_i , \tilde \phi;0),$$
 and hence a contradiction to~$(\tilde \phi, \tilde \phi)$ being a NE.
\end{enumerate}
    Thus~$(f^{-1}\left(2 \rho \right), f^{-1}\left(2 \rho \right))$ is the unique NE.
    
\item  \textit{Proof of Equilibrium Cycle (EC):} 

We will first prove that, 
\begin{equation}
    f^{-1}\left(2 \rho \right) < \eL < \eU \leq f^{-1}\left(\rho \right). \label{Eqn_claim1}
\end{equation}

By definition,~$\eL = m(\eU)/e,$ and~$\eU$ is the unique maximizer of~$m$ (by Lemma~\ref{lem_concavity_of_m},~$m$ is strictly concave).

We first prove that~$f^{-1}(2\rho) < \eU.$ If~$m'(\eU) > 0,$ then~$\eU = \phi_h$. From hypothesis, 
we have,~$f^{-1}(2\rho) < \phi_b \leq \phi_h,$  and thus~$f^{-1}(2\rho) < \eU.$
Else, if~$m'(\eU) = 0$, by Lemma~\ref{lem_monotonicity_of_h_1_function},~$\db$ is strictly decreasing and hence~$\db(f^{-1}(2\rho)) > \db (\phi_b) = 0$ and  thus~$m'(f^{-1}(2\rho)) > 0.$
Using strict concavity of~$m,$  
$f^{-1}(2\rho) < \eU.$

Next,~$e f^{-1}(2\rho) = m (f^{-1}(2\rho)) < m(\eU) = e \eL$ and so,~$f^{-1}(2\rho) < \eL.$ 

Further, by Lemma~\ref{lemma_oscillatory} (given that~$f^{-1}(2\rho) < \eU$), we have~$e\eL = m(\eU) < e\eU,$  and thus~$\eL < \eU.$

Finally, if~$f^{-1} (\rho) = \phi_h$, then~$\eU \leq f^{-1}\left(\rho \right).$ If not, it follows that $f(f^{-1}(\rho)) = \rho,$ which implies~$m(0) = 0 = m(f^{-1}(\rho)).$ Again by strict concavity of~$m$,  we  have,~$\eU \leq f^{-1}\left(\rho \right).$ 
Thus~\eqref{Eqn_claim1} is proved.
From~\eqref{Eqn_claim1} and Table~\ref{tab:WE_MR_PB}, for any~$\phi_{-i} \in \left[\eL, \eU\right],$ and for any~$\phi_i$, we have
\begin{equation} \label{eqn_MR_equilibrium_cycle_proof}
\MR_i(\phi_i, \phi_{-i}; 0) =
 \begin{cases}
 e \phi_i & \text{if } \phi_i < \phi_{-i} \\
 \frac{\Lambda}{2} f(\phi_{i}) (\phi_{i}) & \text{if }  \phi_i = \phi_{-i} \\
 m(\phi_i) \mathds{1}_{\{\phi_i < f^{-1}(\rho)\}} + 0\ \mathds{1}_{\{\phi_i \geq f^{-1}(\rho)\}} = \max \{ m(\phi_i),0\}
 & \text{if }  \phi_i > \phi_{-i}.
\end{cases}   
\end{equation}
We first analyse the revenue rate of platform~$i,$ when it plays outside the EC, and the opponent plays inside the EC. Fix~$\phi_{-i} \in [\eL,\eU].$ From~\eqref{eqn_MR_equilibrium_cycle_proof},  
\begin{itemize}
    \item For any~$\tilde \phi_i \in [0, \eL)$,~$\MR_i(\tilde \phi_i, \phi_{-i}; 0) = e\tilde \phi_i < e \eL = m(\eU).$
    \item For any~$\tilde \phi_i \in (\eU, \phi_h],$~$\MR_i(\tilde \phi_i, \phi_{-i}; 0) = \max\{ m(\tilde \phi_i), 0\} < m(\eU).$
\end{itemize}
Thus in all, for any~$\phi_{-i} \in [\eL, \eU]$ and~$\tilde \phi_i \in [0,\phi_h] \setminus [\eL,\eU]$,
\begin{equation} \label{eqn_MR_equilibrium_cycle_proof_bigger_phir}
\MR_i(\tilde \phi_i, \phi_{-i}; 0) < m(\eU) = e\eL.   
\end{equation}

Now, we prove the three conditions of the Equilibrium Cycle (EC).

\noindent \textit{Proof of condition~(i):} 
If~$\phi_{-i} = \eL$ then choose~$\phi_i = \eU$. For any outside~$\tilde \phi_i \in [0,\phi_h] \setminus [\eL,\eU]$, from~\eqref{eqn_MR_equilibrium_cycle_proof} and~\eqref{eqn_MR_equilibrium_cycle_proof_bigger_phir}, 
~$\MR_i(\tilde \phi_i, \phi_{-i}; 0) < m(\eU) = \MR_i(\phi_i, \phi_{-i}; 0)~$ .

On the other hand, if~$\phi_{-i} \in (\eL, \eU]$ then the same is true similarly by choosing a $\phi_i \in (\eL, \phi_{-i})$,
\begin{equation}\label{eqn_MR_equilibrium_cycle_proof_three}
\MR_i(\phi_i, \phi_{-i}; 0) = e \phi_i > e \eL  = m(\eU)> \MR_i(\tilde \phi_i, \phi_{-i}; 0) \text{ for all } \tilde \phi_i \in [0,\phi_h] \setminus [\eL, \eU]. \nonumber
\end{equation}

\noindent \textit{Proof of condition~(ii):} Consider any~$\VPhi \in [\eL, \eU]^2$. W.l.g., it suffices to show the result for following cases:
\begin{enumerate}[$(a)$]
    \item~$\phi_i < \phi_{-i}$: choose any
~$\phi_{i}' \in (\phi_{i}, \phi_{-i})$.  From~\eqref{eqn_MR_equilibrium_cycle_proof},
$$
  \MR_{i}(\phi_i', \phi_{-i} ;0 ) = e\phi_{i}' >  e\phi_{i} = \MR_{i}(\phi_i, \phi_{-i}; 0). 
$$
\item~$\phi_i = \phi_{-i} \neq \eL$: choose\footnote{\label{footnote_1} For any~$\phi_i > \eL,$ from~\eqref{Eqn_claim1},~$\eL > f^{-1} \left( 2\rho \right)~$,  and thus by Lemma~\ref{lemma_oscillatory},  
~$\Lambda f(\phi_{i}) \phi_{i}  /2  < e\phi_{i}~$.
   Thus,~$  \max \left \{\nicefrac{\Lambda f(\phi_{i})\phi_{i}}{2e}, \eL \right \}    < \phi_{i}$.
   } a~$\phi_{i}' \in \left( \max \left \{\nicefrac{\Lambda f(\phi_{i})\phi_{i}}{2e}, \eL \right \} , \phi_{i} \right) \subset [\eL, \phi_{i}).$ From~\eqref{eqn_MR_equilibrium_cycle_proof},
$$\MR_{i}(\phi_i', \phi_{-i}; 0) = e\phi_{i}' > \frac{\Lambda f(\phi_{i})\phi_{i}}{2} =  \MR_{i}(\phi_i, \phi_{-i}; 0).$$
    \item~$\phi_i = \phi_{-i} = \eL$: choose~$\phi_i' = \eU.$ From~\eqref{eqn_MR_equilibrium_cycle_proof},
 $$ \MR_{i}(\phi_i', \phi_{-i}; 0) = m(\eU) = e \eL
 >\frac{\Lambda f\left(\eL\right)\eL}{2}   = \MR_{i}(\phi_i, \phi_{-i}; 0).$$
\end{enumerate}
The second part of condition~$(ii)$ for each of the above three parts follows from~\eqref{eqn_MR_equilibrium_cycle_proof_bigger_phir} as~$\MR_i(\phi_i',\phi_{-i};0) \geq e\eL.$

\noindent \textit{Proof of condition~(iii):} To the contrary, assume that the interval~$[c,d] \subsetneq [\eL,\eU]$ is an EC. Note that EC is a closed interval, thus either one of~$c > \eL$ or~$c = \eL$ must hold.

\begin{itemize}
    \item If~$c > \eL$, then fix~$\phi_{-i} = c$. Choose\footref{footnote_1}~$\phi_{i}' \in \left( \max \left \{\nicefrac{\Lambda f(c)c}{2e}, \eL \right \}, c \right)$. Clearly, we have~$e\phi_i' > e\eL = m(\eU) \geq m(\phi_i)$ for any~$\phi_i$ and from Lemma~\ref{lemma_oscillatory}~$e\phi_i' > \nicefrac{\Lambda}{2}f(c)c$, and so from~\eqref{Eqn_claim1},~\eqref{eqn_MR_equilibrium_cycle_proof},  and by our choice of~$\phi_{i}'$, 
$$
\MR_i(\phi_i', c; 0) = e \phi_i' > \max \left\{ \frac{\Lambda}{2}f(c)c, \ m(\phi_i) \right \} \geq \MR_i(\phi_i, c; 0) \text{ for any~$\phi_i \in [c, d]$},
$$
which is a
contradiction to condition~(i) in the definition of the EC.
\item If~$c = \eL$ then~$d < \eU$. Fix~$\phi_{-i} = c.$ Choose~$\phi_i' = \eU,$ then from~\eqref{Eqn_claim1},~$\eU < f^{-1}(\rho)$, from~\eqref{eqn_MR_equilibrium_cycle_proof}, and Lemma~\ref{lemma_oscillatory}, for all~$\phi_i \in [c, d]$
$$\MR_i(\eU, c; 0) =  m(\eU) = e \eL > \max \left \{\frac{\Lambda}{2}f(c)c ,  \ m(\phi_i) \right \}  \geq  \MR_i(\phi_i, c; 0).$$
 Since~$d < \eU$, we know~$\eU \notin [c,d]$, a contradiction to condition~(i) of the definition of the EC.

\end{itemize}

\noindent \textit{Proof of Mixed NE:}

\noindent By definition,~$\sigma^*([\eL, \eL]) = 0,$~$\sigma^*([\eL, \eU]) = 1$, and~$\sigma^*([\eL, \phi])$ is a strictly increasing function of~$\phi$; thus~$\sigma^*$ is a probability measure with support in~$[\eL,\eU]$.  From~\eqref{eqn_MR_equilibrium_cycle_proof}, the revenue rate of platform~$i$ for any~$\phi_i \in [\eL,\eU]$, when opponent plays the mixed strategy~$\sigma^*$ is,
\begin{align}
\MR_i (\phi_i, \sigma^*; 0) &=  \int_{\eL}^{\phi_i} m(\phi_i) d\sigma^*(\phi_{-i}) + \int_{\phi_i}^\eU
 e\phi_i d\sigma^* (\phi_{-i}) \nonumber \\
 &=  m(\phi_i) \sigma^*([\eL, \phi_i])+ e\phi_i (1- \sigma^*([\eL, \phi_i])) = m(\eU), \label{eqn_mr_mixed_NE}
 \end{align}
where the last equality is achieved by substituting~$\sigma^*([\eL, \phi_i])$ from the hypothesis. Thus,~$\MR_i (\phi_i, \sigma^*; 0)$  is constant for all~$\phi_i \in [\eL,\eU]$. Further, from~\eqref{eqn_MR_equilibrium_cycle_proof}, and~\eqref{eqn_mr_mixed_NE},
\begin{itemize}
    \item for any~$\tilde \phi \in [0, \eL)$,~$\MR_i(\tilde \phi, \sigma^*; 0) = e\tilde \phi < e \eL = m(\eU),$ and
    \item for any~$\tilde \phi \in (\eU, \phi_h],$~$\MR_i(\tilde \phi, \sigma^*; 0) = m(\tilde \phi) <  m(\eU).$
\end{itemize}
Thus $(\sigma^*, \sigma^*)$ is a mixed strategy Nash Equilibrium.

\item Consider~$ \delta~$ as given in the hypothesis. With~$\rho \geq 1 \geq f(\phi)$ for any~$\phi$, and hence from~\eqref{eqn_f_inverse},~$f^{-1}(\rho) = 0$.
Thus from Table~\ref{tab:WE_MR_PB}, the revenue rate of platform~$i$ against such~$\delta~$ is, 
\begin{equation*}
\MR_i(\phi_i, \delta; 0) =
 \begin{cases}
 0 & \text{ if } \phi_i > \delta \\
 (\nicefrac{\Lambda}{2}) f(\delta) \delta & \text{ if }  \phi_i = \delta \\
 \Lambda f(\phi_i) \phi_i & \text{ if }  \phi_i < \delta
\end{cases}   
\end{equation*}
From the hypothesis of the theorem,~$\MR_i(\phi_i,\delta;0) < \epsilon$ and thus~$(\delta,\delta)$ is an~$\epsilon$-NE.
\end{enumerate}
\end{proof}

\begin{proof}[Proof of Lemma~\ref{lem_mixed_NE_payoff_and_sec_value}] 
The average payoff of each platform under the mixed NE $(\sigma^*,\sigma^*)$ follows easily from~\eqref{eqn_mr_mixed_NE}. 

Consider any platform~$i$. From Table\ref{tab:WE_MR_PB}, the revenue rate of platform~$i$ is
\begin{equation}
\label{eqn_security_MR}
    \MR_i(\VPhi) = \begin{cases}
    e\phi_i & \text{if } \phi_{i} \leq f^{-1} (2\rho) ,
    \\
    \mathds{1}_{\{\phi_{i} < \phi_{-i} \}} e \phi_i + \mathds{1}_{\{\phi_{i} = \phi_{-i}\}} \frac{\Lambda}{2}f(\phi_i) \phi_i  + \mathds{1}_{\{\phi_{i} > \phi_{-i} \}} m(\phi_i) & \text{if } \phi_{i}\in [f^{-1} (2\rho), f^{-1} (\rho)], 
    \\
    \mathds{1}_{\{\phi_{i} < \phi_{-i}\}} \min\{e, \Lambda f(\phi_i) \} \phi_i + \mathds{1}_{\{\phi_{i} = \phi_{-i}\}} \frac{\Lambda}{2}f(\phi_i) \phi_i  + \mathds{1}_{\{\phi_{i} > \phi_{-i} \}}  0 & \text{if } \phi_{i} \geq f^{-1} (2\rho).
    \end{cases}
\end{equation}
We calculate the minimum value that platform~$i$ can obtain when it plays strategy~$\phi_i$ in each of the three ranges considered in~\eqref{eqn_security_MR}. From the definition of the security value (see Footnote~\ref{footnote_Security_value_def}),
\begin{align*}
\underline{\MR_i} 
&= \max \left\{ \max_{\phi_i \leq f^{-1} (2\rho) }\min_{\phi_{-i}} \MR_i(\VPhi), 
\max_{\phi_i \in [f^{-1} (2\rho), f^{-1} (\rho)] }\min_{\phi_{-i}} \MR_i(\VPhi),
\max_{\phi_i \geq f^{-1} (\rho) }\min_{\phi_{-i}} \MR_i(\VPhi)
\right \} \\
&=
\max \{ e f^{-1}(2\rho), m(\eU), 0 \} .
\end{align*}
The second term in the above maximization follows because: (a)~from Lemma~\ref{lemma_oscillatory}, we have~$e\phi_i > \nicefrac{\Lambda}{2}f(\phi_i) \phi_i > m(\phi_i)$ for~$\phi_i \in [ f^{-1} (2\rho),  f^{-1} (\rho)],$ and (b)~$\eU$ is the maximizer of the function~$m$ with~$ f^{-1}(2\rho) < \eU < f^{-1}(\rho)$.
Further, observe that~$e f^{-1}(2\rho) = m (f^{-1}(2\rho)) < m(\eU),$ and thus we have 
$\underline{\MR_i} =  m(\eU)$ and the security strategy is $\eU.$
\end{proof}

\begin{proof}[Proof of Theorem~\ref{thm_B_system_epsilon_equilibrium}]
From Lemma~\ref{lem_mr_continuity_B_system}, the function~$\MR$ is continuous over~$\hat{\mathcal{F}} \setminus \mathcal{H}$, where~$\hat{\mathcal{F}}:= [0, \phi_h]^2 \times [0, \bar \beta],$~$\bar \beta$ is sufficiently large positive constant, and~$\mathcal{H} := \{ (\phi_1, \phi_2, 0) : \phi_1 = \phi_2, \phi_1 > f^{-1}(2\rho) \} .$ We use this continuity to prove the following cases; where in the first case is  when~$\rho \leq f(\phi_b)/2$
\begin{enumerate}[$(1)$]
\item
Using uniform continuity of~$\MR_i,$ on the compact set~$\{(\phi_1, \phi_2, \beta) \in \hat{\mathcal{F}} : \phi_2 = f^{-1}(2\rho) \},$
for any~$\epsilon > 0,$  there exists~$\bar \beta_\epsilon > 0$ such that whenever~$\beta < \bar \beta_\epsilon,$ for every~$\phi \in [0, \phi_h],$
\begin{align*}
\MR_i(\phi, f^{-1}(2\rho); \beta) \leq \MR_i(\phi, f^{-1}(2\rho); 0) + \frac{\epsilon}{2}
&\leq \MR_i(f^{-1}(2\rho), f^{-1}(2\rho); 0) + \frac{\epsilon}{2} \\
&\leq \MR_i(f^{-1}(2\rho), f^{-1}(2\rho); \beta) + \epsilon.
\end{align*}
Here the first, and the last inequality follows from the uniform continuity of~$\MR_i$, while the middle one follows as~$(f^{-1}(2\rho), f^{-1}(2\rho))$ is NE in the \idp\ regime by Theorem~\ref{thm_B_sys_limit_sys_main}, when~$\rho \leq f(\phi_b)/2$.

\item \textit{Proof of Mixed NE:}

Consider the probability measure~$\sigma^*$ as given in Theorem~\ref{thm_B_sys_limit_sys_main} over~$[0, \phi_h]$. Note that~$\bar M := \Lambda \phi_h~$ is a finite uniform upper bound on~$\MR_i$. Let~$B(\phi_i, \eta)$ denote an open ball of radius~$\eta$ centered at~$\phi_i.$ Choose~$\eta$ small enough such that~$\sigma^*(B(\phi_i, \eta))  \le \nicefrac{\epsilon}{8 {\bar M}}.$ Let~$\mathcal{F}_\eta := [0, \phi_h]^2 \times [0, \bar \beta] \setminus B(\mathcal{H}, \eta),$ where~$B(\mathcal{H}, \eta)$ is an open ball around the set~$\mathcal{H}$ with radius~$\eta$; thus~$\mathcal{F}_\eta~$ is a compact set.
From Lemma~\ref{lem_mr_continuity_B_system},  the function~$\MR_i$ is uniformly continuous on the set~$\mathcal{F}_\eta$.
Now using uniform continuity of function~$\MR_i$ on the set~$\mathcal{F}_{\eta}$, there exists~$ \bar \beta_{\epsilon} > 0$ such that,
\begin{equation}\label{eqn_mr_EC_close_by_beta_two}
| \MR_i(\phi_i, \phi_{-i}; 0) - \MR_i(\phi_i, \phi_{-i}; \beta) | < \frac{\epsilon}{4}, \text{ when }  (\VPhi, \beta), (\VPhi, 0) \in \mathcal{F}_{\eta}   \text{ and when } \beta \leq \bar \beta_\epsilon .
\end{equation}
Using~\eqref{eqn_mr_EC_close_by_beta_two}, for any~$\phi_i \in [0,\phi_h]$, and if~$ \beta \le {\bar \beta}_\epsilon$
\begin{align} \label{eqn_epsilon_mix_NE}
\left| \MR_i(\phi_i, \sigma^*; \beta) - \MR_i(\phi_i, \sigma^*; 0) \right| 
 & \leq \int_{[0, \phi_h] \setminus  B(\phi_i,  \eta )} \left| \MR_i(\phi_i, \phi_{-i}; \beta) - \MR_i(\phi_i, \phi_{-i}; 0) \right| d\sigma^*(\phi_{-i}) \nonumber \\ & \hspace{1cm} + 2\bar M \sigma^*(B(\phi_i, \eta)) \nonumber \\
 & \leq \frac{\epsilon}{4} + \frac{\epsilon}{4} =  \frac{\epsilon}{2} . 
\end{align}
It now follows that for all~$\beta \leq \bar \beta_\epsilon$,
\begin{align} 
    \left| \MR_i(\sigma^*, \sigma^*; \beta) - \MR_i(\sigma^*, \sigma^*; 0) \right|
    &\leq  \int_{0}^{\phi_h} \left| \MR_i(\phi_i, \sigma^*; \beta) - \MR_i(\phi_i, \sigma^*; 0) \right| d\sigma^*(\phi_{i}) \nonumber \\
    &\leq \frac{\epsilon}{2} \int_{0}^{\phi_h} d\sigma^*(\phi_{i}) = \frac{\epsilon}{2}. \label{eqn_epsilon_mix_NE_one}
\end{align}
From Theorem~\ref{thm_B_sys_limit_sys_main},~$\MR_i(\phi_i, \sigma^*; 0) \leq \MR_i(\sigma^*, \sigma^*; 0)$ for any~$\phi_i$, as~$\sigma^*$ is mixed NE when $\beta = 0.$ 
Thus
from~\eqref{eqn_epsilon_mix_NE},~\eqref{eqn_epsilon_mix_NE_one}, 
\begin{equation*}
    \MR_i(\phi_i, \sigma^*, \beta) \leq \MR_i(\phi_i, \sigma^*; 0)  + \frac{\epsilon}{2} \leq \MR_i(\sigma^*, \sigma^*; 0) + \frac{\epsilon}{2} \leq \MR_i(\sigma^*, \sigma^*, \beta) + \epsilon.
\end{equation*}
Therefore using Theorem~\ref{thm_B_sys_limit_sys_main}, and the definition~\ref{def_epsilon_equilibria}, the interval~$[\eL,\eU]$ is an~$\epsilon$-mixed NE. 

\textit{Proof~$\epsilon$-EC:}

From Table~\ref{tab:WE_MR_PB}, for any~$\epsilon > 0,$ define
\begin{align} 
 \delta_\epsilon &= \inf_{\phi_{i} \in [\eL +\epsilon, \eU]} \left(\MR_i(\phi_{i} - \epsilon, \phi_{i}; 0) - \MR_i(\phi_{i}, \phi_{i}; 0)\right) \nonumber \\
 &= \inf_{\phi_{i} \in [\eL +\epsilon, \eU]} e(\phi_i-\epsilon)-\frac{\Lambda}{2}f(\phi_i)\phi_i = e \eL - \frac{\Lambda}{2} f(\eL + \epsilon)(\eL + \epsilon) .  \label{eqn_positive_delta_EC}
\end{align}
The last equality follows as~$(e - \nicefrac{\Lambda}{2}f(\phi))\phi$ is increasing in~$\phi.$  From Lemma~\ref{lemma_oscillatory},~$\delta_\epsilon$ at $\epsilon = 0$,~$(e - \nicefrac{\Lambda}{2} f(\eL ))\eL > 0,$ and thus there exists~$\bar \epsilon > 0$ such that~$\delta_\epsilon > 0$ for all~$\epsilon \leq \bar \epsilon.$ Further choose smaller~$\bar \epsilon$ if required such that~$\eL + \epsilon < \eU - \epsilon$ for all~$\epsilon \leq \bar \epsilon.$ Consider one such~$\epsilon > 0,$ and choose~$\eta > 0$ such that

\begin{equation} \label{eqn_eta_definition}
\eta < \frac{\min\left \{e, (m(\eU) - m(\eU + \epsilon)) , (m(\eU) - m(\eU - \epsilon)), \delta_\epsilon, \left(\frac{e \eL}{2} - \frac{\Lambda f(\eL) \eL}{4}\right), \left(m(\eU) -\frac{\Lambda f(\eL)\eL}{2e} \right) \right\}  }{\epsilon},
\end{equation}
From Lemma~\ref{lem_eta_positive_epsilon_EC_proof}, the numerator of the above equation is positive; thus one can choose such~$\eta > 0.$
For such~$\epsilon$ and~$\eta$ (which depends on~$\epsilon$), from uniform continuity of~$\MR$ over~$\hat{\mathcal{F}} \setminus {B}(\mathcal{H},\epsilon/2),$ choose~$\bar \beta_{\epsilon} > 0$ such that for all~$ \beta \leq \bar \beta_{\epsilon}$,
\begin{align}\nonumber
\text{for any } (\phi_i, \phi_{-i})  \in [0, \phi_h]^2, \ \text{ with } |\phi_i - \phi_{-i}| \geq \frac{\epsilon}{2}, \text{ we have } \\
| \MR_i(\phi_i, \phi_{-i}; 0) - \MR_i(\phi_i, \phi_{-i}; \beta) | <  \frac{\eta \epsilon}{2} .\label{eqn_mr_EC_close_by_beta}
\end{align}
\noindent \textit{Proof of condition~(i):}

\noindent If~$\phi_{-i} \in [\eL + \epsilon, \eU - \epsilon]$ then choose~$\phi_i = \eL$.

\noindent If~$\tilde \phi_i < \eL - \epsilon$ then for all~$\beta \leq \bar \beta_\epsilon$, using~\eqref{eqn_mr_EC_close_by_beta}, Table~\ref{tab:WE_MR_PB}, and~\eqref{eqn_eta_definition},
$$
\MR_i(\eL, \phi_{-i}; \beta) - \MR_i(\tilde \phi_i, \phi_{-i}; \beta) \geq e \eL - e \tilde \phi_i - \eta \epsilon\geq e \eL- e (\eL - \epsilon) - \eta \epsilon = e\epsilon - \eta\epsilon > 0 .
$$
If~$\tilde \phi_i > \eU + \epsilon$ then for all~$\beta \leq \bar \beta_\epsilon$, using~\eqref{eqn_mr_EC_close_by_beta}, Table~\ref{tab:WE_MR_PB}, and~\eqref{eqn_eta_definition}, (recall~$m$ decreases beyond~$\eU$, and~$m(\eU) = e\eL$)
\begin{align*}
\MR_i(\eL, \phi_{-i}; \beta) - \MR_i(\tilde \phi_i, \phi_{-i}; \beta)
&\geq e\eL - m(\tilde \phi_i) - \eta \epsilon \\ 
&\geq e\eL - m(\eU + \epsilon) - \eta \epsilon \\
&= m(\eU) - m(\eU + \epsilon) - \eta \epsilon > 0. 
\end{align*}

\noindent \textit{Proof of condition~(ii):}
for any~$\VPhi \in [\eL + \epsilon, \eU - \epsilon]^2$

\noindent If~$\phi_i > \phi_{-i},$ then choose~$\phi_{i}' =  \eL$.
Thus, for all~$\beta \leq \bar \beta_\epsilon$, using~\eqref{eqn_mr_EC_close_by_beta}, Table~\ref{tab:WE_MR_PB}, and~\eqref{eqn_eta_definition}, ($m$ is increasing till~$\eU$)
$$
\MR_i(\phi_i', \phi_{-i}; \beta) - \MR_i(\phi_i, \phi_{-i}; \beta) \geq e \eL - m(\phi_{i}) - \eta \epsilon \geq m(\eU) - m(\eU - \epsilon) - \eta \epsilon > 0 .
$$
\noindent If~$\phi_i = \phi_{-i}$ then choose~$\phi_i' = \phi_{i} - \epsilon$.Thus, for all~$\beta \leq \bar \beta_\epsilon$, using~\eqref{eqn_mr_EC_close_by_beta},~\eqref{eqn_positive_delta_EC}, Table~\ref{tab:WE_MR_PB}, and~\eqref{eqn_eta_definition},
\begin{equation*}
\MR_i(\phi_i', \phi_{-i}; \beta) - \MR_i(\phi_{i}, \phi_{-i}; \beta) \geq  \delta_\epsilon - \eta \epsilon > 0.
\end{equation*}
The second part of the~$\epsilon$-EC condition (ii) can be deduced similarly to condition (i) presented earlier. Specifically, if~$\phi_i > \phi_{-i}$, the proof can be obtained by setting~$\phi_i' = \eL$ and considering~$\phi_{-i} \in [\eL + \epsilon, \eU - \epsilon]$. On the other hand, if~$\phi_i = \phi_{-i}$, a similar approach as that used for condition~(i) can be applied, and
\begin{itemize}
    \item if~$\tilde \phi_i < \eL - \epsilon,$ then~$\MR_i(\phi_i', \phi_{-i}; \beta) - \MR_i(\tilde \phi_i, \phi_{-i}; \beta) \geq e (\phi_{i} - \epsilon) - e \tilde \phi_i - \eta \epsilon \geq e \eL - e \tilde \phi_i - \eta \epsilon > 0$,
    \item if~$\tilde \phi_i > \eU + \epsilon,$ then~$\MR_i(\phi_i', \phi_{-i}; \beta) - \MR_i(\tilde \phi_i, \phi_{-i}; \beta) \geq e (\phi_{i} - \epsilon) - m( \tilde \phi_i) - \eta \epsilon \geq e \eL - m(\eU) - \eta \epsilon > 0$.
\end{itemize}

\noindent \textit{Proof of condition~(iii):}
Towards condition~(iii), we prove a stronger statement, to be precise we will show:  for any subset~$[c,d] \subset [a + \epsilon,b - \epsilon],$ there exists~$ i \in \mathcal{N}$, and~$\phi_{-i} \in [ c , d ],$  such that, 
\begin{equation}\label{eqn_epsilon_EC_con_3}
\sup_{\tilde \phi_i \in [0, \phi_h]\setminus [c- \epsilon, d+\epsilon]} \mathcal{M}_i(\tilde \phi_i,\phi_{-i}) > \mathcal{M}_i(\phi_i,\phi_{-i})  \text{ for all } \phi_i \in [c,d].    
\end{equation}

To the contrary, assume that the interval~$[c,d] \subsetneq [\eL + \epsilon,\eU - \epsilon]$ is an~$\epsilon$-EC. Note that~$\epsilon$-EC is a closed interval, thus either one of~$c > \eL + \epsilon$ or~$c =  \eL + \epsilon$ must hold.

\begin{enumerate}[$(1)$]
    \item If~$c > \eL + \epsilon$, choose~$\phi_{-i} = c$ and~$\phi_{i}' = \left( \max \left \{\nicefrac{\Lambda f(c)c}{2e} , \eL \right \} + c \right)/2$; observe\footref{footnote_1}~$\phi_i' \in [c,d]^c$. We have a contradiction to condition~(i) in the~$\epsilon$-EC definition because by~\eqref{eqn_mr_EC_close_by_beta}, Table~\ref{tab:WE_MR_PB},~\eqref{eqn_eta_definition},  and by our choice of~$\phi_{i}',$ we have
    \begin{equation} \label{eqn_epsilon_EC_eq_1}
    \begin{aligned}
     \MR_i(\phi_i', c; \beta) \geq e \phi_i' -\frac{ \eta \epsilon}{2} &\stackrel{(a)}{>} \left( \mathds{1}_{\{\phi_i = c\}} \frac{\Lambda}{2}f(c)c + \mathds{1}_{\{\phi_i \in (c,d]\}} m(\phi_i) \right) +\frac{ \eta \epsilon}{2} \\
     &\geq  \MR_i(\phi_i, c; \beta) \text{ for any~$\phi_i \in [c, d]$},  
     \end{aligned}
    \end{equation}
    because of the following reasons:
    \begin{enumerate}[(i)]
        \item firstly inequality~$(a)$ is true with~$\phi_{i}' = (\nicefrac{\Lambda f(c)c}{2e} + c)/2$
        and for   choice of~$\eta$  as in~\eqref{eqn_eta_definition}, because:
    \begin{align}
    &e \phi_i' - \frac{\Lambda}{2}f(c)c = \frac{c}{2} \left (e - \frac{\Lambda}{2}f(c)\right) > \frac{\eL}{2}\left(e - \frac{\Lambda}{2}f(\eL)\right ) > \eta \epsilon, \mbox{ and, } \nonumber \\
    &e \phi_i' - m(\phi_i) > e \eL - m(d) \geq m(\eU) - m(\eU - \epsilon) > \eta \epsilon, \mbox{ for any }  \phi_i \in (c,d], \label{eqn_epsilon_EC_con_3_2}
    \end{align}
        towards the above, observe~$e > \nicefrac{\Lambda}{2}f(c)$,~$\phi_i' > \eL,$ and  recall~$m$ is concave and it's maximizer is at~$\eU > d$,

    \item  clearly~\eqref{eqn_epsilon_EC_con_3_2} also holds for~$\phi_i' = (\eL + c)/2$ and any~$\phi_i \in (c,d]$ and thus
    inequality~$(a)$ is also true with~$\phi_i' = (\eL + c)/2,$ as further, 
\begin{align*}
    e \phi_i' - \frac{\Lambda}{2}f(c)c = e \frac{(c+ \eL)}{2} - \frac{\Lambda}{2}f(c)c &= \frac{c(e - \Lambda f(c)) + e \eL}{2} \\
    & >   \frac{\eL (2e - \Lambda f(\eL))}{2} > \eta \epsilon .    
\end{align*}

\end{enumerate}
    
\item If~$c = \eL + \epsilon
$ then~$d  < \eU - \epsilon$. Choose~$\phi_{-i} = c,$ and~$\phi_i' = \eU;$ then we have a contradiction to condition~(i) in the~$\epsilon$-EC definition because from~\eqref{eqn_mr_EC_close_by_beta}, Table~\ref{tab:WE_MR_PB},~\eqref{eqn_eta_definition}, for all~$\phi_i \in [c, d]$, by our choice of~$\phi_i'$, and~$\eta~$,

\begin{align*} 
\MR_i(\eU, c; \beta) \geq   m(\eU) - \frac{ \eta \epsilon}{2}  &= e \eL - \frac{ \eta \epsilon}{2} \\
& > \left( \mathds{1}_{\{\phi_i = c\}} \frac{\Lambda}{2}f(c)c + \mathds{1}_{\{\phi_i \in (c,d]\}} m(\phi_i) \right) + \frac{ \eta \epsilon}{2} \\
&\geq  \MR_i(\phi_i, c; \beta),     
\end{align*}  
and this is because:  a) we have~$e\eL - \nicefrac{\Lambda}{2}f(c)c = \delta_\epsilon > \eta \epsilon$ from~\eqref{eqn_eta_definition} (also recall~$e \eL = m(\eU)$), and b) for all~$\phi_i \in (c, d],$ from~\eqref{eqn_eta_definition}, 
    \begin{equation*}
     e \eL - m(\phi_i) =  m(\eU) - m(\phi_i) \geq m(\eU) - m(d) \geq m(\eU) - m(\eU - \epsilon) > \eta \epsilon.   
    \end{equation*}
\end{enumerate}
\end{enumerate}
\end{proof}

\section{Proofs related to Section~\ref{sec:delay}} \label{sec:appendix_wait_metric}

\begin{proof}[Proof of Theorem~\ref{thm_WE_MR_new}] 
The proof of the theorem follows  as in that of  Theorem~\ref{thm_WE_MR_PB}, but with extra arguments required to accommodate the parameter $\alpha$.

Using Lemma \ref{lem_PB_joint_continuity_from_scratch},~$\WB_i$ as defined in \eqref{eqn_wait_metric_defn}--\eqref{eqn_QoS_WB_zero_beta} is jointly continuous in $(\VPhi, \lambda_i, \beta)$ for each $i$ and continuity w.r.t. $\alpha$ is obvious. 
 Thus the  objective function $\Gamma_D = (\WB_1(\lambda) - \WB_2(\Lambda - \lambda))^2$ given in~\eqref{eqn_existence_uniq_WE_main} with $Q_i = \WB_i$ of~\eqref{eqn_wait_metric_defn}--\eqref{eqn_QoS_WB_zero_beta} is jointly continuous in $(\lambda_1, \beta, \alpha)$, which provides the following:  a) 
 the WE is a minimizer of $\Gamma_D$   over the  compact domain $[0, \Lambda]$, and thus a WE exists for each $(\VPhi, \beta, \alpha)$;
 b)  by Maximum theorem~\cite{maximum} for parametric minimization of $\Gamma_D$ with $(\beta, \alpha) $ as parameters, the correspondence of WEs,  $(\beta, \alpha) \mapsto \lambda_i(\VPhi; \beta, \alpha)$ is upper semi-continuous. 

The proof is completed in two more steps:
 a) we show the uniqueness of WE at $\beta > 0$;  and 
 b) we choose WE at $\beta = 0$ as in Table~\ref{tab:WE_MR_PB} for $\alpha = 0,$ and  as in Table~\ref{tab:WE_MR_WB} for $\alpha > 0$ such that the mapping $(\beta, \alpha) \mapsto \lambda_i(\VPhi; \beta, \alpha)$ is continuous.

(a) \textit{Uniqueness of WE at $\beta >0$}:  Towards this, we first  write $\WB_i$ of \eqref{eqn_wait_metric_defn} as below:
\begin{align*}
  \WB_i(\VPhi; \lambda_i, \beta , \alpha) 
    &=  (1 - f(\phi_i))  +  f(\phi_i) \E\left[ \alpha \tau(N_i)  \right] , \text{ where } \\
  \alpha  \tau(n) &=  
    \begin{cases}
      1  &  \mbox{if } n = 0 , \\
    \alpha \left(\frac{1}{n}\right) \indic{n \leq \N} &  \mbox{else.} 
    \end{cases}  
\end{align*}
By hypothesis, $\alpha < 1$ and so $\tau$ is strictly decreasing in $n$ for $n \leq \N$, for each $i$. Thus, by Lemma~\ref{lem_coupling_arguments}, $\WB_i$ is strictly increasing function of $\lambda_i$ for any given $(\VPhi, \beta, \alpha),$ with $\beta > 0.$ 

 (b) \textit{Choosing WE at $\beta = 0$}: Since the correspondence~$(\beta, \alpha) \mapsto \lambda_i(\VPhi; \beta, \alpha)$ is compact and upper semi-continuous, and since WE is unique for every $\beta > 0$, the correspondence becomes a function (see, \cite{maximum}), which is continuous on $\{ (\beta, \alpha): \beta > 0$\}. Thus $\lim_{\beta \to 0} (\lambda_1(\VPhi; \beta, \alpha), \lambda_2(\VPhi; \beta, \alpha))$ for any $\alpha$ exists, and the rest of the proof computes these limits and shows them equal to the ones reported in Tables \ref{tab:WE_MR_PB}, \ref{tab:WE_MR_WB}.
 
 By upper semi-continuity of the set of WEs, these limits are WEs at $\beta = 0$. Thus, we derive the limits of Tables \ref{tab:WE_MR_PB}, \ref{tab:WE_MR_WB}, either by establishing the uniqueness of WE at $\beta = 0$ (in some cases) or by finding the limits of  WEs as $\beta \to 0$ (in the rest of the cases).  
This part is established case-wise. Towards the end, we discuss $\MR_i$.

\noindent \textit{Case~$(1)$}:~$\phi_1 = \phi_2$. From Lemma~\ref{lem_exist_unique_WE},  Lemma~\ref{lem_coupling_arguments}
($\WB_i$ is strictly increasing) and~\eqref{eqn_wait_metric_defn}, it is easy to see that the unique WE for~$\beta > 0$ and any $\alpha$  is~$(\lambda_1(\VPhi; \beta, \alpha), \lambda_1(\VPhi; \beta, \alpha)) =(\nicefrac{\Lambda}{2}, \nicefrac{\Lambda}{2})$.  Thus, we set $\lambda_1(\VPhi; 0, \alpha) = \nicefrac{\Lambda}{2}$ in both the Tables~\ref{tab:WE_MR_PB}, \ref{tab:WE_MR_WB}.

\noindent \textit{Case~$(2)$}:~$\phi_1 > \phi_2$

\noindent \textit{Case $(2a)$}:~$\phi_1, \phi_2 \in [0, f^{-1}\left( \rho \right))$. Observe,

\vspace{-4mm}
{\small\begin{equation}\label{eqn_WB_wait_metric_IDP}
\WB_i(\VPhi; \lambda_i, 0, \alpha) = \begin{cases}
  1 -  \frac{e}{\lambda_i}   +  \alpha   f(\phi_i)   \left(1 - \frac{e}{\lambda_i f(\phi_i)}\right)  \sum_{n=1}^\N  \frac{1}{n}\left(\frac{e}{\lambda_i f(\phi_i)}\right)^n  \hspace{-4mm} & \mbox{if }   \frac{e}{\lambda_i f(\phi_i)} < 1 ,\\ 
 1 -  f(\phi_i)    & \mbox{else. } 
 \end{cases}     
\end{equation}}
Using simple calculus-based arguments, it is easy to observe that $\WB_1$ is strictly increasing when $\nicefrac{e}{(\lambda_1 f(\phi_1))} < 1$ or $\lambda_1 > \nicefrac{e}{f(\phi_1)}$ and $\WB_2$ is strictly decreasing w.r.t. $\lambda_1$  when $\nicefrac{e}{(\lambda_2 f(\phi_2))} < 1$ or when $\lambda_1 < \Lambda - \nicefrac{e}{f(\phi_2)}$, for any given $(\VPhi, \alpha)$. Thus, $\WB_1 - \WB_2$ is non-decreasing over $[0, \Lambda]$ (since $\phi_1 > \phi_2,$ we have $1 - f(\phi_1) > 1 - f(\phi_2)$).

In view of the above, to prove the uniqueness\footnote{Once it is established that $\WB_1 - \WB_2 < 0$ at $\lambda_1 = 0$ and $\WB_1 - \WB_2 > 0$ at $\lambda_1 = \Lambda$, we will further have three sub-cases based on values of  $\Lambda - \nicefrac{e}{f(\phi_2)}, \nicefrac{e}{f(\phi_1)} $; in each sub-case using simple algebra, it is easy to prove uniqueness.}  of WE for this sub-case, it suffices to show that $\WB_1(\VPhi; \lambda_1, 0, \alpha) - \WB_2(\VPhi; \Lambda-\lambda_1, 0, \alpha) < 0$ at $\lambda_1 = 0$ and $\WB_1(\VPhi;   \lambda_1, \alpha) - \WB_2(\VPhi; \Lambda-\lambda_1, 0, \alpha) > 0$ at $\lambda_1 = \Lambda$. Towards that consider the difference at $\lambda_1 = 0$,
\begin{align} \nonumber 
\WB_1(\VPhi; 0, 0, \alpha) - \WB_2(\VPhi; \Lambda, 0, \alpha) \hspace{-35mm} &\\
\nonumber 
&=  1 -  f(\phi_1) - \left(1 -  f(\phi_2)  + \alpha    f(\phi_2)   \left( 1 - \frac{e}{\Lambda f(\phi_2)} \right) \left( \frac{1}{ \alpha} + \sum_{n=1}^\N  \frac{\left( \frac{e}{\Lambda f(\phi_2)} \right)^n}{n}\right)\right) \\
\nonumber 
&=  f(\phi_2) - f(\phi_1)   - \alpha    \left(f(\phi_2) - \frac{e}{\Lambda}\right) \left( \frac{1}{ \alpha} + \sum_{n=1}^\N  \frac{\left( \frac{e}{\Lambda f(\phi_2)} \right)^n}{n}\right) \\
\nonumber
&=  f(\phi_2) - f(\phi_1)   -     \left(f(\phi_2) - \frac{e}{\Lambda}\right) - \left( \alpha    \left(f(\phi_2) - \frac{e}{\Lambda}\right) \sum_{n=1}^\N  \frac{\left( \frac{e}{\Lambda f(\phi_2)} \right)^n}{n}\right) \\
&=  \left( \frac{e}{\Lambda} - f(\phi_1) \right) - \left( \alpha    \left(f(\phi_2) - \frac{e}{\Lambda}\right) \sum_{n=1}^\N  \frac{\left( \frac{e}{\Lambda f(\phi_2)} \right)^n}{n}\right). \label{eqn_WB_1_minus_WB_2}
\end{align}
Since $\phi_1 > \phi_2$ and $\phi_1 < f^{-1}(\rho),$ we have $\nicefrac{e}{\Lambda} - f(\phi_1) < 0$ and $f(\phi_2) - \nicefrac{e}{\Lambda} > 0,$ which shows that $\WB_1(\VPhi; 0, 0, \alpha) - \WB_2(\VPhi; \Lambda, 0, \alpha) < 0$. Along exactly similar lines, we have $\WB_1(\VPhi; \Lambda, 0, \alpha) - \WB_2(\VPhi; 0, 0, \alpha) > 0$. 

By Maximum theorem, the above uniqueness implies continuity of WE for $\alpha \in (0,1)$ and $\beta = 0$. In the proof of Theorem~\ref{thm_WE_MR_PB} (see Case~(2a), Case~(2b)) it is proved that the WE is unique at $(\beta, \alpha) = (0,0)$ for this sub-case. Thus, we have  continuity of WE for $\alpha \in [0,1)$ and $\beta = 0$.

\noindent \textit{Case $(2b)$}:~When~$\phi_1 \in [f^{-1}\left( \rho \right), \phi_h]$, from~\eqref{eqn_WB_wait_metric_IDP},~$ \WB_1(\VPhi; \lambda_1, 0, \alpha) = 1 - f(\phi_1)$ irrespective of~$\lambda_1$.
Then the following difference (simplified as in \eqref{eqn_WB_1_minus_WB_2}),  equals, 
{\small \begin{align}
\nonumber
 \WB_1( &\VPhi; \lambda_1, 0, \alpha) - \WB_2(\VPhi; \Lambda - \lambda_1, 0, \alpha)  \\
&=  \left \{ \begin{array}{ll}
  \left( \frac{e}{\Lambda - \lambda_1} - f(\phi_1) \right) + \alpha    \left(\frac{e}{\Lambda - \lambda_1} - f(\phi_2)  \right) \sum_{n=1}^\N  \frac{\left( \frac{e}{(\Lambda - \lambda_1) f(\phi_2)} \right)^n}{n}   & \mbox{ if }  \lambda_1 \leq \Lambda - \frac{e}{f(\phi_2)} \\
 f(\phi_2) - f(\phi_1)    &  \mbox{ else.}
 \end{array}
  \right . \label{eqn_WB_1_minus_WB_2_next_case}
\end{align}}
\noindent Using derivatives, one can show that the above increases monotonically w.r.t. $\lambda_1$ when $ \lambda_1 \leq \Lambda - \nicefrac{e}{f(\phi_2)}$. We now consider further sub-cases based on $\phi_2$:
\begin{enumerate}
\item If $\phi_2 \leq f^{-1}(\rho)$, i.e., when $f(\phi_2) \geq \nicefrac{e}{\Lambda}$. 
In this sub-case,  from \eqref{eqn_WB_1_minus_WB_2_next_case},
$ \WB_1(\VPhi; \lambda_1, 0, \alpha) - \WB_2(\VPhi; \Lambda - \lambda_1, 0, \alpha) > 0$ at $\lambda_1 = 0$, $\alpha = 0$, and   there exists a $\bar \alpha > 0$ such that at  $\lambda_1 = 0 ,$
\begin{align*}
\WB_1(\VPhi; \lambda_1, 0, \alpha) - \WB_2(\VPhi; \Lambda - \lambda_1, 0, \alpha) & \geq 0 \mbox{ for all } \alpha \le \bar \alpha,  \mbox{ and } \\
\WB_1(\VPhi; \lambda_1, 0, \alpha) - \WB_2(\VPhi; \Lambda - \lambda_1, 0, \alpha) &< 0 \mbox{ if }  \alpha >  \bar \alpha  .
\end{align*}
Thus for all $0\le \alpha < \min\{1, \bar \alpha\} $, by monotonicity, the difference in \eqref{eqn_WB_1_minus_WB_2_next_case} remains positive for all $\lambda_1$   and hence the WE from definition \eqref{eqn_existence_uniq_WE_main} is $\lambda_1^* = 0$.

If $\bar \alpha < 1$, then  for all $\alpha \in [\bar \alpha, 1)$,  the WE as before is the unique zero of \eqref{eqn_WB_1_minus_WB_2_next_case} (observe the difference at $\lambda_1 = \nicefrac{e}{f(\phi_2)}$ equals $f(\phi_2)- f(\phi_1) >0$). Thus the WE is unique, as given in Tables \ref{tab:WE_MR_PB},~\ref{tab:WE_MR_WB} and is continuous for all $\alpha \ge 0$ and $\beta = 0$ even for this sub-case.



%
\item If $\phi_2 > f^{-1}(\rho)$, i.e., when $f(\phi_2) < \nicefrac{e}{\Lambda}$ then from~\eqref{eqn_WB_wait_metric_IDP},~$ \WB_2(\VPhi; \Lambda-\lambda_1, 0, \alpha) = 1 - f(\phi_2)$ for all~$\lambda_2$. Thus
$$
\WB_1(\VPhi; \lambda_1, 0, \alpha) - \WB_2(\VPhi; \Lambda - \lambda_1, 0, \alpha) = f(\phi_2) - f(\phi_1) \mbox{ for all } \lambda_1.
$$
Thus from \eqref{eqn_existence_uniq_WE_main}, any $\lambda_1 \in [0, \Lambda]$ is a WE, and for continuity, we chose the following as limit WE: 
%
%
%
$$(\lambda_1(\VPhi; 0, \alpha),  \lambda_2(\VPhi; 0, \alpha)) = \left(0, \Lambda \right)$$
Now, proof to show that the above selection ensures the required continuity properties follows exactly as given in the proof of Theorem~\ref{thm_WE_MR_PB}, Case~(2c) with $\phi_2 > f^{-1}(\rho)$ --- observe continuity of $\WB_i$ is also provided in Lemma~\ref{lem_PB_joint_continuity_from_scratch} and the term $\beta_\epsilon$ given in that proof also depends upon $\alpha$. Thus, the WE in this case is the same as that provided in both the Tables~\ref{tab:WE_MR_PB},~\ref{tab:WE_MR_WB}.
\end{enumerate}

The continuity of~$\MR_i$ (after defining it at $\beta = 0$ using the limit WE) follows by Lemma~\ref{lem_approx_MR} -- observe from above proof that as $(\beta_n, \alpha_n) \to (0, \alpha), $ we have $\lambda_i(\VPhi; \beta_n, \alpha_n) \to \lambda_i(\VPhi; 0, \alpha) $.

\end{proof}

{\begin{proof}[Proof of Theorem~\ref{thm_new_system_epsilon_equilibrium}]
Using similar techniques as  Lemma~\ref{lem_mr_continuity_B_system}, the function~$\MR$ is continuous over~$\hat{\mathcal{F}} \setminus \mathcal{H}$, where~$\hat{\mathcal{F}}:= [0, \phi_h]^2 \times [0, \bar \beta]\times [0, \hat \alpha],$~$\bar \beta$ is positive constant (finite), $\hat \alpha \in (0, \bar \alpha)$ and~$\mathcal{H} := \{ (\phi_1, \phi_2, 0, \alpha) : \phi_1 = \phi_2, \phi_1 > f^{-1}(2\rho) \} .$ Using this continuity, the remaining part of the proof follows arguments analogous to those in the proof of Theorem~\ref{thm_B_system_epsilon_equilibrium}.
\ignore{\begin{enumerate}[$(1)$]
\item
Using uniform continuity of~$\MR_i,$ on the compact set~$\{(\phi_1, \phi_2, \beta, \alpha) \in \hat{\mathcal{F}} : \phi_2 = f^{-1}(2\rho) \},$
for any~$\epsilon > 0,$  there exists~$\bar \beta_\epsilon > 0$, $\bar \alpha_\epsilon \in (0, \bar \alpha]$ such that whenever~$\beta < \bar \beta_\epsilon$ and $ \alpha < \bar \alpha_\epsilon$, for every~$\phi \in [0, \phi_h],$
\begin{align*}
\MR_i(\phi, f^{-1}(2\rho); \beta, \alpha) \leq \MR_i(\phi, f^{-1}(2\rho); 0, 0) + \frac{\epsilon}{2}
&\leq \MR_i(f^{-1}(2\rho), f^{-1}(2\rho); 0, 0) + \frac{\epsilon}{2} \\
&\leq \MR_i(f^{-1}(2\rho), f^{-1}(2\rho); \beta, \alpha) + \epsilon.
\end{align*}
Here, the first and the last inequality follows from the uniform continuity of~$\MR_i$, while the middle one follows as~$(f^{-1}(2\rho), f^{-1}(2\rho))$ is NE in the \idp\ regime by Theorem~\ref{thm_B_sys_limit_sys_main}, when~$\rho \leq f(\phi_b)/2$.

\item \textit{Proof of Mixed NE:}

\noindent Consider the probability measure~$\sigma^*$ as given in Theorem~\ref{thm_B_sys_limit_sys_main} over~$[0, \phi_h]$. Note that~$\bar M := \Lambda \phi_h~$ is a finite uniform upper bound on~$\MR_i$. Let~$B(\phi_i, \eta)$ denote an open ball of radius~$\eta$ centered at~$\phi_i.$ Choose~$\eta$ small enough such that~$\sigma^*(B(\phi_i, \eta))  \le \nicefrac{\epsilon}{8 {\bar M}}.$ Let~$\mathcal{F}_\eta := [0, \phi_h]^2 \times [0, \bar \beta] \times [0, \bar \alpha] \setminus B(\mathcal{H}, \eta),$ where~$B(\mathcal{H}, \eta)$ is an open ball around the set~$\mathcal{H}$ with radius~$\eta$; thus~$\mathcal{F}_\eta~$ is a compact set.
From Lemma~\ref{lem_mr_continuity_Wait_system},  the function~$\MR_i$ is uniformly continuous on the set~$\mathcal{F}_\eta$.
Now using uniform continuity of function~$\MR_i$ on the set~$\mathcal{F}_{\eta}$, there exists~$\bar \beta_\epsilon > 0$, $\bar \alpha_\epsilon \in (0, \bar \alpha]$ such that,
\begin{equation}\label{eqn_mr_EC_close_by_beta_two}
| \MR_i(\phi_i, \phi_{-i}; 0, 0) - \MR_i(\phi_i, \phi_{-i}; \beta, \alpha) | < \frac{\epsilon}{4}, \text{ when }  (\VPhi, \beta, \alpha), (\VPhi, 0, 0) \in \mathcal{F}_{\eta}   \text{ and } \beta \leq \bar \beta_\epsilon, \alpha \leq \bar \alpha_\epsilon .
\end{equation}
Using~\eqref{eqn_mr_EC_close_by_beta_two}, for any~$\phi_i \in [0,\phi_h]$, and if~$ \beta \le {\bar \beta}_\epsilon, \alpha \leq \bar \alpha_\epsilon,$
\begin{align} \label{eqn_epsilon_mix_NE}
\left| \MR_i(\phi_i, \sigma^*; \beta, \alpha) - \MR_i(\phi_i, \sigma^*; 0, 0) \right| 
 & \leq \int_{[0, \phi_h] \setminus  B(\phi_i,  \eta )} \left| \MR_i(\phi_i, \phi_{-i}; \beta, \alpha) - \MR_i(\phi_i, \phi_{-i}; 0, 0) \right| d\sigma^*(\phi_{-i}) \nonumber \\ & \hspace{1cm} + 2\bar M \sigma^*(B(\phi_i, \eta)) \nonumber \\
 & \leq \frac{\epsilon}{4} + \frac{\epsilon}{4} =  \frac{\epsilon}{2} . 
\end{align}
It now follows that for all~$\beta \leq \bar \beta_\epsilon, \alpha \leq \bar \alpha_\epsilon,$
\begin{align} 
    \left| \MR_i(\sigma^*, \sigma^*; \beta, \alpha) - \MR_i(\sigma^*, \sigma^*; 0, 0) \right|
    &\leq  \int_{0}^{\phi_h} \left| \MR_i(\phi_i, \sigma^*; \beta, \alpha) - \MR_i(\phi_i, \sigma^*; 0, 0) \right| d\sigma^*(\phi_{i}) \nonumber \\
    &\leq \frac{\epsilon}{2} \int_{0}^{\phi_h} d\sigma^*(\phi_{i}) = \frac{\epsilon}{2}. \label{eqn_epsilon_mix_NE_one}
\end{align}
From Theorem~\ref{thm_B_sys_limit_sys_main},~$\MR_i(\phi_i, \sigma^*; 0, 0) \leq \MR_i(\sigma^*, \sigma^*; 0, 0)$ for any~$\phi_i$, as~$\sigma^*$ is mixed NE when $\beta = 0, \ \alpha = 0.$ 
Thus
from~\eqref{eqn_epsilon_mix_NE},~\eqref{eqn_epsilon_mix_NE_one}, 
\begin{equation*}
    \MR_i(\phi_i, \sigma^*, \beta, \alpha) \leq \MR_i(\phi_i, \sigma^*; 0, 0)  + \frac{\epsilon}{2} \leq \MR_i(\sigma^*, \sigma^*; 0, 0) + \frac{\epsilon}{2} \leq \MR_i(\sigma^*, \sigma^*, \beta, \alpha) + \epsilon.
\end{equation*}
Therefore using Theorem~\ref{thm_B_sys_limit_sys_main}, and the definition~\ref{def_epsilon_equilibria}, the interval~$[\eL,\eU]$ is an~$\epsilon$-mixed NE.

\textit{Proof~$\epsilon$-EC:}

From Table~\ref{tab:WE_MR_PB}, for any~$\epsilon > 0,$ define
\begin{align} 
 \delta_\epsilon &= \inf_{\phi_{i} \in [\eL +\epsilon, \eU]} \left(\MR_i(\phi_{i} - \epsilon, \phi_{i}; 0, 0) - \MR_i(\phi_{i}, \phi_{i}; 0, 0)\right) \nonumber \\
 &= \inf_{\phi_{i} \in [\eL +\epsilon, \eU]} e(\phi_i-\epsilon)-\frac{\Lambda}{2}f(\phi_i)\phi_i = e \eL - \frac{\Lambda}{2} f(\eL + \epsilon)(\eL + \epsilon) .  \label{eqn_positive_delta_EC}
\end{align}
The last equality follows as~$(e - \nicefrac{\Lambda}{2}f(\phi))\phi$ is increasing in~$\phi.$  From Lemma~\ref{lemma_oscillatory},~$\delta_\epsilon$ at $\epsilon = 0$,~$(e - \nicefrac{\Lambda}{2} f(\eL ))\eL > 0,$ and thus there exists~$\bar \epsilon > 0$ such that~$\delta_\epsilon > 0$ for all~$\epsilon \leq \bar \epsilon.$ Further, choose smaller~$\bar \epsilon$ if required such that ~$\eL + \epsilon < \eU - \epsilon$ for all~$\epsilon \leq \bar \epsilon.$ Consider one such~$\epsilon > 0,$ and choose~$\eta > 0$ such that
{\small 
\begin{equation} \label{eqn_eta_definition}
\eta < \frac{\min\left \{e, (m(\eU) - m(\eU + \epsilon)) , (m(\eU) - m(\eU - \epsilon)), \delta_\epsilon, \left(\frac{e \eL}{2} - \frac{\Lambda f(\eL) \eL}{4}\right), \left(m(\eU) -\frac{\Lambda f(\eL)\eL}{2e} \right) \right\}  }{\epsilon},
\end{equation}} %
From Lemma~\ref{lem_eta_positive_epsilon_EC_proof}, the numerator of the above equation is positive; thus one can choose such~$\eta > 0.$
For such~$\epsilon$ and~$\eta$ (which depends on~$\epsilon$), from uniform continuity of ~$\MR$ over~$\hat{\mathcal{F}} \setminus {B}(\mathcal{H},\epsilon/2),$ choose~$\bar \beta_{\epsilon} > 0$, $\bar \alpha_\epsilon \in (0, \bar \alpha]$ such that for all~$ \beta \leq \bar \beta_{\epsilon}$, $\alpha \leq \bar \alpha_{\epsilon},$
\begin{equation}\label{eqn_mr_EC_close_by_beta}
\text{for any } (\phi_i, \phi_{-i})  \in [0, \phi_h]^2, \ \text{ with } |\phi_i - \phi_{-i}| \geq \frac{\epsilon}{2}, \text{ we have } | \MR_i(\phi_i, \phi_{-i}; 0, 0) - \MR_i(\phi_i, \phi_{-i}; \beta, \alpha) | <  \frac{\eta \epsilon}{2} .
\end{equation}
\noindent \textit{Proof of condition~(i):}

\noindent If~$\phi_{-i} \in [\eL + \epsilon, \eU - \epsilon]$ then choose~$\phi_i = \eL$.

\noindent If~$\tilde \phi_i < \eL - \epsilon$ then for all~$\beta \leq \bar \beta_\epsilon$, $\alpha \leq \bar \alpha_{\epsilon},$ using~\eqref{eqn_mr_EC_close_by_beta}, Table~\ref{tab:WE_MR_PB}, and~\eqref{eqn_eta_definition},
$$
\MR_i(\eL, \phi_{-i}; \beta, \alpha) - \MR_i(\tilde \phi_i, \phi_{-i}; \beta, \alpha) \geq e \eL - e \tilde \phi_i - \eta \epsilon\geq e \eL- e (\eL - \epsilon) - \eta \epsilon = e\epsilon - \eta\epsilon > 0 .
$$
If~$\tilde \phi_i > \eU + \epsilon$ then for all~$\beta \leq \bar \beta_\epsilon$, using~\eqref{eqn_mr_EC_close_by_beta}, Table~\ref{tab:WE_MR_PB}, and~\eqref{eqn_eta_definition}, (recall~$m$ decreases beyond~$\eU$, and~$m(\eU) = e\eL$)
\begin{align*}
\MR_i(\eL, \phi_{-i}; \beta) - \MR_i(\tilde \phi_i, \phi_{-i}; \beta) \geq e\eL - m(\tilde \phi_i) - \eta \epsilon  &\geq e\eL - m(\eU + \epsilon) - \eta \epsilon \\
&= m(\eU) - m(\eU + \epsilon) - \eta \epsilon > 0.    
\end{align*}

\noindent \textit{Proof of condition~(ii):}
for any~$\VPhi \in [\eL + \epsilon, \eU - \epsilon]^2$

\noindent If~$\phi_i > \phi_{-i},$ then choose~$\phi_{i}' =  \eL$.
Thus, for all~$\beta \leq \bar \beta_\epsilon$, $\alpha \leq \bar \alpha_{\epsilon},$ using~\eqref{eqn_mr_EC_close_by_beta}, Table~\ref{tab:WE_MR_PB}, and~\eqref{eqn_eta_definition}, ($m$ is increasing till~$\eU$)
$$
\MR_i(\phi_i', \phi_{-i}; \beta, \alpha) - \MR_i(\phi_i, \phi_{-i}; \beta, \alpha) \geq e \eL - m(\phi_{i}) - \eta \epsilon \geq m(\eU) - m(\eU - \epsilon) - \eta \epsilon > 0 .
$$
\noindent If~$\phi_i = \phi_{-i}$ then choose~$\phi_i' = \phi_{i} - \epsilon$.Thus, for all~$\beta \leq \bar \beta_\epsilon$, $\alpha \leq \bar \alpha_{\epsilon},$ using~\eqref{eqn_mr_EC_close_by_beta},~\eqref{eqn_positive_delta_EC}, Table~\ref{tab:WE_MR_PB}, and~\eqref{eqn_eta_definition},
\begin{equation*}
\MR_i(\phi_i', \phi_{-i}; \beta, \alpha) - \MR_i(\phi_{i}, \phi_{-i}; \beta, \alpha) \geq  \delta_\epsilon - \eta \epsilon > 0.
\end{equation*}
The second part of the~$\epsilon$-EC condition (ii) can be deduced similarly to condition (i) presented earlier. Specifically, if~$\phi_i > \phi_{-i}$, the proof can be obtained by setting~$\phi_i' = \eL$ and considering~$\phi_{-i} \in [\eL + \epsilon, \eU - \epsilon]$. On the other hand, if~$\phi_i = \phi_{-i}$, a similar approach as that used for condition~(i) can be applied, and
\begin{itemize}
    \item if ~$\tilde \phi_i < \eL - \epsilon,$ then~$\MR_i(\phi_i', \phi_{-i}; \beta, \alpha) - \MR_i(\tilde \phi_i, \phi_{-i}; \beta, \alpha) \geq e (\phi_{i} - \epsilon) - e \tilde \phi_i - \eta \epsilon \geq e \eL - e \tilde \phi_i - \eta \epsilon > 0$,
    \item if ~$\tilde \phi_i > \eU + \epsilon,$ then~$\MR_i(\phi_i', \phi_{-i}; \beta, \alpha) - \MR_i(\tilde \phi_i, \phi_{-i}; \beta, \alpha) \geq e (\phi_{i} - \epsilon) - m( \tilde \phi_i) - \eta \epsilon \geq e \eL - m(\eU) - \eta \epsilon > 0$.
\end{itemize}

\noindent \textit{Proof of condition~(iii):}
Towards condition~(iii), we prove a stronger statement, to be precise we will show:  for any subset~$[c,d] \subset [a + \epsilon,b - \epsilon],$ there exists~$ i \in \mathcal{N}$, and~$\phi_{-i} \in [ c , d ],$  such that, 
\begin{equation}\label{eqn_epsilon_EC_con_3}
\sup_{\tilde \phi_i \in [0, \phi_h]\setminus [c- \epsilon, d+\epsilon]} \mathcal{M}_i(\tilde \phi_i,\phi_{-i}) > \mathcal{M}_i(\phi_i,\phi_{-i})  \text{ for all } \phi_i \in [c,d].    
\end{equation}

To the contrary, assume that the interval~$[c,d] \subsetneq [\eL + \epsilon,\eU - \epsilon]$ is an~$\epsilon$-EC. Note that~$\epsilon$-EC is a closed interval, thus either one of~$c > \eL + \epsilon$ or~$c =  \eL + \epsilon$ must hold.

\begin{enumerate}[$(1)$]
    \item If~$c > \eL + \epsilon$, choose~$\phi_{-i} = c$ and ~$\phi_{i}' = \left( \max \left \{\nicefrac{\Lambda f(c)c}{2e} , \eL \right \} + c \right)/2$; observe\footref{footnote_1} ~$\phi_i' \in [c,d]^c$. We have a contradiction to condition~(i) in the~$\epsilon$-EC definition because by~\eqref{eqn_mr_EC_close_by_beta}, Table~\ref{tab:WE_MR_PB},~\eqref{eqn_eta_definition},  and by our choice of~$\phi_{i}',$ we have
    \begin{equation} \label{eqn_epsilon_EC_eq_1}
    \begin{aligned}
     \MR_i(\phi_i', c; \beta, \alpha) \geq e \phi_i' -\frac{ \eta \epsilon}{2} &\stackrel{(a)}{>} \left( \mathds{1}_{\{\phi_i = c\}} \frac{\Lambda}{2}f(c)c + \mathds{1}_{\{\phi_i \in (c,d]\}} m(\phi_i) \right) +\frac{ \eta \epsilon}{2} \\
     &\geq  \MR_i(\phi_i, c; \beta, \alpha) \text{ for any~$\phi_i \in [c, d]$},  
     \end{aligned}
    \end{equation}
    because of the following reasons:
    \begin{enumerate}[(i)]
        \item firstly inequality~$(a)$ is true with  ~$\phi_{i}' = (\nicefrac{\Lambda f(c)c}{2e} + c)/2$
        and for   choice of~$\eta$  as in ~\eqref{eqn_eta_definition}, because:
    \begin{align}
    &e \phi_i' - \frac{\Lambda}{2}f(c)c = \frac{c}{2} \left (e - \frac{\Lambda}{2}f(c)\right) > \frac{\eL}{2}\left(e - \frac{\Lambda}{2}f(\eL)\right ) > \eta \epsilon, \mbox{ and, } \nonumber \\
    &e \phi_i' - m(\phi_i) > e \eL - m(d) \geq m(\eU) - m(\eU - \epsilon) > \eta \epsilon, \mbox{ for any }  \phi_i \in (c,d], \label{eqn_epsilon_EC_con_3_2}
    \end{align}
        towards the above, observe~$e > \nicefrac{\Lambda}{2}f(c)$, ~$\phi_i' > \eL,$ and  recall ~$m$ is concave and it's maximizer is at~$\eU > d$,

    \item  clearly~\eqref{eqn_epsilon_EC_con_3_2} also holds for ~$\phi_i' = (\eL + c)/2$ and any~$\phi_i \in (c,d]$ and thus
    inequality~$(a)$ is also true with ~$\phi_i' = (\eL + c)/2,$ as further, 
   ~$$
    e \phi_i' - \frac{\Lambda}{2}f(c)c = e \frac{(c+ \eL)}{2} - \frac{\Lambda}{2}f(c)c = \frac{c(e - \Lambda f(c)) + e \eL}{2}>   \frac{\eL (2e - \Lambda f(\eL))}{2} > \eta \epsilon .
   ~$$
    \end{enumerate}
    
\item If~$c = \eL + \epsilon
$ then~$d  < \eU - \epsilon$. Choose~$\phi_{-i} = c,$ and~$\phi_i' = \eU;$ then we have a contradiction to condition~(i) in the~$\epsilon$-EC definition because from~\eqref{eqn_mr_EC_close_by_beta}, Table~\ref{tab:WE_MR_PB},~\eqref{eqn_eta_definition}, for all~$\phi_i \in [c, d]$, by our choice of~$\phi_i'$, and~$\eta~$,

\begin{align*} 
\MR_i(\eU, c; \beta, \alpha) &\geq   m(\eU) - \frac{ \eta \epsilon}{2} \\
& = e \eL - \frac{ \eta \epsilon}{2} > \left( \mathds{1}_{\{\phi_i = c\}} \frac{\Lambda}{2}f(c)c + \mathds{1}_{\{\phi_i \in (c,d]\}} m(\phi_i) \right) + \frac{ \eta \epsilon}{2} \geq  \MR_i(\phi_i, c; \beta, \alpha),     
\end{align*}  
and this is because:  a) we have~$e\eL - \nicefrac{\Lambda}{2}f(c)c = \delta_\epsilon > \eta \epsilon$ from~\eqref{eqn_eta_definition} (also recall  ~$e \eL = m(\eU)$), and b) for all~$\phi_i \in (c, d],$ from~\eqref{eqn_eta_definition}, 
    \begin{equation*}
     e \eL - m(\phi_i) =  m(\eU) - m(\phi_i) \geq m(\eU) - m(d) \geq m(\eU) - m(\eU - \epsilon) > \eta \epsilon.   
    \end{equation*}
\end{enumerate}
\end{enumerate}
}
\end{proof}
}

\ignore{
\begin{lemma}\label{lem_wait_metric_joint_continuity}
The function $\WB_i:[0, \phi_h] \times [0,\Lambda] \times [0, \infty) \times [0, \gamma_i) \to \mathbb{R}$ is jointly continuous for each $i$.
\end{lemma}
\begin{proof}[Proof of Lemma~\ref{lem_wait_metric_joint_continuity}]
From \eqref{eqn_exp_delay_expression}, the continuity of $\WB_i$ w.r.t $\phi_i$ is immediate, once $f$ is continuous; we consider continuity w.r.t. the remaining three parameters $(\lambda_i, \beta, \alpha)$.
The required  joint continuity for the first term in \eqref{eqn_QoS_WB_positive_beta}-\eqref{eqn_QoS_WB_zero_beta} is already established in Lemma~\ref{lem_PB_joint_continuity_from_scratch}; continuity of the second term at all points with $\beta > 0$ follows as in Lemma~\ref{lem_PB_joint_continuity_from_scratch}, thus, it suffices to show the continuity of the following term of \eqref{eqn_QoS_WB_positive_beta}-\eqref{eqn_QoS_WB_zero_beta} at any $(\lambda_i, \beta, \alpha)$ with $\beta = 0$,
\begin{align}\label{eqn_d_continuity_eqn}
  \alpha \frac{ \sum_{n=1}^{\infty} \varrho^n \mu_{i,n}}{\sum_{n=1}^{\infty} \mu_{i,n}}  \mathds{1}_{\{\beta > 0\}} - \mathds{1}_{\{\beta = 0\}} \mathds{1}_{\{\rho_i < 1\}}      \alpha \frac{ \varrho \rho_i (1 - \rho_i) }{(1 - \varrho \rho_i)}  , \ \text{ where } \rho_i = \frac{e}{\lambda_i f(\phi_i)}. 
\end{align}
Observe from \eqref{eqn_QoS_WB_positive_beta}-\eqref{eqn_QoS_WB_zero_beta} that for $\rho_i < 1,$ the denominator term $\sum_{n = 0}^{\infty} \mu_{i,n}$
converges to $\nicefrac{1}{(1 - \rho_i)}$ and the numerator term $ \sum_{n=1}^{\infty} \varrho^n \mu_{i,n}$
converges to $\nicefrac{ \varrho \rho_i }{(1 - \varrho \rho_i)}  $, both as $\beta \to 0$ and with the ordered pair $(\lambda_i, \phi_i)$ fixed; and thus the continuity at any point $(\lambda_i, \beta, \alpha)$ with $\beta = 0$ and $\rho_i < 1$ again follows as in Lemma~\ref{lem_PB_joint_continuity_from_scratch}.
From \eqref{eqn_d_continuity_eqn}, the proof is complete if we show  for $\rho_i \geq 1$ the following is true\footnote{
For any $\beta > 0,$  $\lim_{k \to \infty} \sum_{n=1}^\infty \mu_n$ exists (finite), and thus $\lim_{k \to \infty} \sum_{n=1}^\infty \nicefrac{\mu_n}{n}$ also exists, because:
$$
\mu_n = \prod_{a =1}^n \frac{e}{c + a \beta}  \le \left(\frac{e}{\beta}\right)^n 
\frac{1}{n!}.
$$}:
\begin{align}\label{eqn_series_conv_wait_metric}
\frac{ \lim_{k \to \infty}\sum_{n=1}^{k} \varrho^n \mu_{i,n} }{ \lim_{k \to \infty} \sum_{n=0}^{k} \mu_{i,n}} \to 0, \text{ as } \beta \to 0 \text{ with } (\lambda_i, \phi_i) \text{ fixed} .
\end{align}
For ease of notation, lets represent $c = \lambda_i f(\phi_i)$ only for the following paragraphs of the rest of the proof. 
Consider a $\Delta$ such that $\rho_i = \nicefrac{e}{c} >  \Delta > 1$ and  $\varrho \Delta < 1$ and then 
define 
\begin{equation} \label{eqn_delta_bound_series_conv_proof}
N_\beta := \min_n \left\{ n : \frac{e}{c + n \beta} <   \Delta \right \};   
\end{equation}
Now consider dividing the numerator and denominator of \eqref{eqn_series_conv_wait_metric}
with $\mu_{N_\beta}$ and observe that  uniformly across all $\beta $ we have:
$$
\frac{ \sum_n  \mu_n  } {\mu_{N_\beta}} \ge 1 \mbox{ for any } \beta . 
$$
Thus suffices to prove that 
$$
\sum_{n=1}^\infty \varrho^n \frac{ \mu_n  } {\mu_{N_\beta}} \to 0 \mbox{ as } \beta  \to 0
$$
Observe $\nicefrac{e}{c+a\beta} < \Delta$ for $a \ge  N_\beta$ and $\nicefrac{e}{c+a\beta} > \Delta$ for $a <  N_\beta$ and hence:
\begin{eqnarray*}
\sum_{n=1}^\infty \varrho^n \frac{ \mu_n  } {\mu_{N_\beta}}  &=& \sum_{n=1}^{N_\beta} \varrho^n \frac{ \mu_n  } {\mu_{N_\beta}}  + \sum_{n=N_\beta }^\infty \varrho^n \frac{ \mu_n  } {\mu_{N_\beta}} \\
&=&  \sum_{n=1}^{N_\beta} \varrho^n \frac{ 1   } {\prod_{a=n+1}^{N_\beta} \frac{e}{c + a \beta}}  + \sum_{n=N_\beta }^\infty \varrho^n  \prod_{a= N_\beta+1}^n \frac{e}{c + a \beta} \\
&\le &  \sum_{n=1}^{N_\beta} \varrho^n \frac{ 1   } {\Delta^{N_\beta-n} }  + \sum_{n=N_\beta }^\infty \varrho^n  \Delta^{n -N_\beta+1}  \\
&\le&  \Delta^{- N_\beta }   \frac{ 1 }{1 - \nicefrac{\varrho}{\Delta}}  + \Delta^{- N_\beta }  \frac{1} {1 - \varrho \Delta } \to 0
\end{eqnarray*}
as $\Delta^{- N_\beta }  \to 0$ once $\beta \to 0.$

\ignore{

To prove \eqref{eqn_series_conv_wait_metric}, the idea is to scale the numerator and denominator appropriately and show that the scaled numerator converges to zero as $\beta \to 0$, while the scaled denominator is uniformly bounded away from 0 (for all small enough $\beta $). 
Towards this: a) consider a $\eta < 1$ with $1 < \nicefrac{1}{\eta} < \nicefrac{e}{c}$ (possible when $e > c$);  b) define $N_\beta$ as below for each $\beta > 0$,
\begin{equation} \label{eqn_delta_bound_series_conv_proof}
N_\beta := \min_n \left\{ n : \frac{e}{c + n \beta} <  \frac{1}{\eta} \right \};   
\end{equation}
and  c) consider the limit of the following `scaled terms' ($\nicefrac{b_\beta}{a_\beta}$ equals the ratio in \eqref{eqn_series_conv_wait_metric} for each $\beta > 0$):
\begin{align*}
\lim_{\beta > 0, \beta \to 0}  
\frac{b_\beta}{a_\beta}  \mbox{ with }  \ \ 
a_\beta =  \lim_{k \to \infty} \sum_{n=1}^k \mu_n  \ \eta^{N_\beta}  \  \text{ and } \  b_\beta =  \lim_{k \to \infty}\sum_{n=1}^k \frac{\mu_n}{n} \eta^{N_\beta}. 
\end{align*}

{\color{red}
Now, if we prove the following two claims, then the proof follows:
\begin{enumerate}[\text{Claim}~(1):]
    \item There exists $\bar \beta > 0$ and $L_\beta > 0$ such that $a_b \leq L_\beta $ fir all $\beta \leq \bar \beta.$
    \item $b_\beta \to 0$ as $\beta \to 0$
\end{enumerate}
  }

\textbf{Step (1):} 
With $M_\beta := \nicefrac{e}{ (c+N_\beta \beta ) }$, 
\begin{align*}
a_\beta = \lim_{k \to \infty} \sum_{n=1}^k \mu_n  \  \eta^{N_\beta} 
\geq &\eta^{N_\beta}  \sum_{n=1}^{N_\beta }
\mu_{n} \  
\ge \eta^{N_\beta}  \sum_{n=1}^{N_\beta } \left ( \frac{e}{c + N_\beta \beta} \right )^n \\
 = & \eta^{N_\beta} \left ( \frac{ 1 - (M_\beta )^{N_\beta +1}  } {1 - M_\beta} - 1 \right )   = \eta^{N_\beta}   \frac{ M_\beta - (M_\beta )^{N_\beta +1}  } {1 - M_\beta} 
\end{align*}
By \eqref{eqn_delta_bound_series_conv_proof},  $M_\beta  \to  \nicefrac{1}{\eta}$  (idea is to match the rates so that denominator has a non-zero lower bound), and so the RHS in the above,  
$$
\eta^{N_\beta}   \frac{ M_\beta - (M_\beta )^{N_\beta +1}  } {1 - M_\beta}  \to  \frac{ 1  } {  1 - \eta }   > 0.
$$
Thus there exists a $\bar \beta > 0$ such that 
$a_\beta > L > 0$ for all $\beta \le \bar \beta.$
Observe here that 
$$
N_\beta \beta \to  \eta e - c  > 0
$$so that $N_\beta \to \infty$ and hence $\eta^{N_\beta} \to 0$, as $\beta \to 0.$

{\bf Step (2):} Say $\eta = \delta \nicefrac{c}{e}$  Now consider 
\begin{eqnarray*}
    b_\beta &=& \lim_{k \to \infty}  \sum_{n=1}^k \mu_n \varrho^n \eta^{N_\beta} = \sum_{n \le N_\beta}   \eta^{N_\beta - n}    \varrho^n \prod_{a=1}^n \frac{ c \delta } { c + a\beta} +  \prod_{a=1}^{N_\beta} \frac{ c \delta } { c + a\beta}  \lim_{k \to \infty}  \sum_{n > N_\beta}^k    \varrho^n \prod_{a=N_\beta + 1}^n \frac{ e } { c + a\beta}  \\
    &\le &     \sum_{n \le N_\beta}  \eta^{ N_\beta}  \left ( \frac{ \varrho \delta  } { \eta} \right )^n+   \sum_{n > N_\beta}   \prod_{a=N_\beta + 1}^n \frac{ e } { c + a\beta}  
\end{eqnarray*} 

\newpage

\begin{eqnarray*}
    b_\beta &=& \lim_{k \to \infty}  \sum_{n=1}^k\frac{ \mu_n} {n} \eta^{N_\beta} = \sum_{n \le N_\beta}   \eta^{N_\beta - n}  \frac{1}{n} \prod_{a=1}^n \frac{ c \delta } { c + a\beta} + \prod_{a=1}^{N_\beta} \frac{ c \delta } { c + a\beta}  \lim_{k \to \infty}  \sum_{n > N_\beta}^k     \frac{1}{n} \prod_{a=N_\beta + 1}^n \frac{ e } { c + a\beta}  \\
    &\le &  \prod_{a=1}^{N_\beta} \frac{ c \delta } { c + a\beta} \frac{1}{N_\beta}  \left (  \sum_{n \le N_\beta}  \eta^{ N_\beta - n}  \left (1+ \frac{N_\beta -n}{ n}\right ) +   \sum_{n > N_\beta}   \prod_{a=N_\beta + 1}^n \frac{ e } { c + a\beta}  \right )
\end{eqnarray*} 

\newpage

Now, there exists a $\bar \beta> 0$ and $L > 0$ such that,
\begin{align*}
\lim_{k \to \infty} \sum_{n=1}^k \mu_n \left( \frac{c}{e} \right)^{N_\beta} 
\geq 
\mu_{N_\beta} \left( \frac{c}{e} \right)^{N_\beta}  
= 
\prod_{a=1}^{N_\beta} \left( \frac{e}{c+a\beta}. \frac{c}{e} \right)
=
\prod_{a=1}^{N_\beta} \left( \frac{c}{c+a\beta} \right)
  =: L_\beta > L > 0,
\end{align*}
for all $\beta \le {\bar \beta}$, because of  the following observation (as the summation can be viewed as Riemann sum and because   $N_\beta \beta \to \nicefrac{ (e- \delta c) }{\delta} $, as $\beta \to 0$, see \eqref{eqn_delta_bound_series_conv_proof}):
\begin{align*}
\beta \log (L_\beta) &= 
\beta \log \left( \prod_{a=1}^{N_\beta}  \frac{c}{c+a\beta} \right) =  \sum_{a=1}^{N_\beta} \log \left(   
   \frac{c}{c+a\beta} \right) \beta \\
\stackrel{\text{as } \beta \to 0}{\longrightarrow}
&\int_0^{ \frac{e-\delta c } {\delta} }
 \log\left(\frac{c}{c+x} \right) dx = (e-c) (\log(c) + 1) - e \log (e) - c \log (c) > -\infty.
\end{align*}
Thus,
$$
\prod_{a=1}^{N_\beta} \left( \frac{c}{c+a\beta} \right) \geq \exp((e-c) (\log(c) + 1) - e \log (e) - c \log (c)) 
 = L_\beta > 0.
$$
This proves claim (1).

\textbf{Claim (2)}:

\newpage 
Towards this, define the sequence $S_n$ as in the following, and it suffices to show that $\lim_{n \to \infty} S_n = 0$ for any $\rho_i \geq 0$:
\begin{align*}
S_n = \frac{\sum_{k=1}^{n} \frac{\rho_i^k}{k}}{\sum_{k=0}^{n} \rho_i^k}
\end{align*}
Observe that for $\rho_i = 1,$ the above expression can be simplifies as
$$
 S_n = \frac{\sum_{k=1}^n \frac{1}{k}}{n} \ \text{ and thus, } \lim_{n \to \infty} S_n = \lim_{n \to \infty} \frac{\sum_{k=1}^n \frac{1}{k}}{n} \leq \lim_{n \to \infty} \frac{1 + \ln(n)}{n} = 0.
$$
For $\rho_i > 1,$ dividing both sides by $\rho_i^n,$ and changing the variables such that $y = \nicefrac{1}{\rho_i},$ we have,
\begin{align*}
\sum_{k=1}^{n} \frac{\rho_i^k}{k}
&= \sum_{k=0}^{n-1} \frac{y^k}{n-k}
= \frac{1}{n}\sum_{k=0}^{n-1} \frac{n y^k}{n-k}
\leq \frac{1}{n}\sum_{k=0}^{n-1} \left( 1 + \frac{k}{n-k}\right) y^k 
< \frac{1}{n}\sum_{k=0}^{n-1} \left( 1 + k \right) y^k\\
&= \frac{1}{n} \frac{d}{dy} \sum_{k=1}^{n} y^k
= \frac{1}{n} \frac{d}{dy} \left( \frac{1 - y^{n+1}}{1-y} - 1 \right) = \frac{1}{n} \frac{(1 - y^{n+1}) - (n+1)y^n(1-y)}{(1-y)^2}.
\end{align*}
Therefore,
\begin{align*}
 S_n =  \frac{\frac{1}{n} \frac{(1 - y^{n+1}) - (n+1)y^n(1-y)}{(1-y)^2}}{\frac{1 - y^{n+1}}{1-y}}, \text{ and thus }  \lim_{n \to \infty}  S_n = 0.
\end{align*}
}
Thus, we have that~$\WB_i(\bphi; \lambda_i,\beta, \alpha)$ is jointly continuous over~$[0, \Lambda]\times [0, \infty) \times [0, \gamma_i)$.
\end{proof}
}

\section{Proofs related to Section~\ref{sec:comparison}}
\label{sec:appendix_comparison}
\begin{proof}[Proof of Lemma~\ref{lem_mono_price_less_than_D}]
Observe that the optimal/equilibrium price is decreasing function of~$\rho,$ thus from Theorems~\ref{thm_monopoly_opt} and~\ref{thm_B_sys_limit_sys_main}, it suffices to show~$ \phi_m \geq \phi_b.$ Recall that~$f$ is strictly concave and decreasing function (see Assumption {\bf A}.1). From the definitions of~$d_m$~and~$ d_b$, 
$$
\text{for any } \phi \in [0, \phi_h], \text{ we have } d_m(\phi) < d_b(\phi).
$$
As~$\phi_m,$ and~$\phi_b$ are respective zeros (if the zero exists in~$[0, \phi_h),$ else~$\phi_h$) of~$d_m$ and~$ d_b$, the first statement of the lemma follow.

Next, note that 
$$ \phi_m \geq \phi_b \Rightarrow \frac{f(\phi_b)}{2} \geq \frac{f(\phi_m)}{2}.$$ 
We now prove the following claim which will be used repeatedly:
\begin{equation} \label{eqn_lem_price_comparison_claim}
\frac{\Lambda}{2} f(\phi_m) \phi_m \geq \frac{\Lambda}{2} f(\eU) \eU, \mbox{ and } \frac{\Lambda}{2} f(\phi_m) \phi_m > m(\eU) \mbox{ when } \eU > f^{-1}(2\rho).
\end{equation}
Firstly, observe from Theorem~\ref{thm_monopoly_opt} that $\phi_m$ is the maximizer of the function $\phi \mapsto f(\phi) \phi,$ thus the first inequality of~\eqref{eqn_lem_price_comparison_claim} follows. From Lemma~\ref{lemma_oscillatory}, $\nicefrac{\Lambda}{2} f(\eU) \eU  > m(\eU)$ when $\eU > f^{-1}(2\rho).$ Together, these inequalities imply the second inequality in~\eqref{eqn_lem_price_comparison_claim}.

Observe that the equilibria under consideration are symmetric, and thus, both the platforms get half of the market share $\nicefrac{\Lambda}{2}$.
From Theorem~\ref{thm_B_sys_limit_sys_main} and Table~\ref{tab:WE_MR_PB}, the equilibrium payoff for $\rho \leq \nicefrac{f(\phi_m)}{2}$ is $e f^{-1}(2\rho),$ which also equals the monopoly payoff. For $\rho \in (\nicefrac{f(\phi_m)}{2}, \nicefrac{f(\phi_b)}{2}],$ the equilibrium payoff equals $$e f^{-1}(2\rho) = \frac{\Lambda f(f^{-1}(2\rho)) f^{-1}(2\rho)}{2} < \frac{\Lambda f(\phi_m) \phi_m}{2},$$ which is the monopoly payoff (the last inequality above follows since~$\phi_m$ is the maximizer of~$f(\phi)\phi$).
    For $\rho > \nicefrac{f(\phi_b)}{2},$ the payoff QoS metric $\PB$ is $\max\{m(\eU), 0\},$ which is again less than or equal to monopoly payoff given that $\eU > f^{-1}(2\rho)$ (see~\eqref{eqn_lem_price_comparison_claim}, and~\eqref{Eqn_claim1}). 
\end{proof}

\section{Additional proofs}

\begin{lemma} \label{lem_concavity_of_m}
Assuming \textbf{A.1} and \textbf{A.2}, for all~$\phi \in [0, \phi_h]$, the function~$m(\phi) :=  \Lambda  f(\phi)\phi - e\phi$ is strictly concave function of~$\phi$.
\end{lemma}

\begin{proof}[Proof of Lemma~\ref{lem_concavity_of_m}]
From \textbf{A.1},~$f$ is a strictly concave and decreasing function. Thus,~$f'$ is negative and decreasing, and~$m'$ is a strictly decreasing function of~$\phi$. Therefore,~$m$ is a strictly concave function of~$\phi$.
\end{proof}

\ignore{\begin{lemma}[Joint continuity of~$\MR_i$ in $\DA$ system] \label{lem_MR_joint_cont}
Define, $\mu_n^i(\VPhi_m;\beta_m) = \prod_{a=1}^n \frac{e}{h(\VPhi_m) + a\beta_{m}}.$
For any~$i\in \mathcal{N}, \ \phi_{i,m} \in [0, \phi_h]$ and~$ \beta_{m} \in [0, \infty)$, 
$$\text{if } (\phi_{1,m}, \phi_{2,m}, \beta_m) \rightarrow (\phi_1, \phi_2, \beta) \text{ then } \sum_{n=0}^{\infty} \mu_n^i(\VPhi_m;\beta_m) \to \sum_{n=0}^{\infty} \mu_n^i(\VPhi;\beta), \text{ and } \MR_i(\phi_{1,m}, \phi_{2,m}; \beta_m) \rightarrow \MR_i(\phi_1, \phi_2; \beta) .$$
\end{lemma}}

\ignore{\begin{lemma}[Joint continuity of~$\PB_i$]\label{lem_PBi_continuity}
For any~$i\in \mathcal{N}, \ \phi_{i,m} \in [0, \phi_h],$~$\lambda_{i,m} \in [0, \Lambda],$ and~$ \beta_{m} \in [0, \infty)$, 
$$\text{if } (\phi_{i,m}, \lambda_{i,m}, \beta_m) \rightarrow (\phi_i, \lambda_{i}, \beta) \text{ then } \PB_i(\phi_{i,m}, \lambda_{i,m}; \beta_m) \rightarrow \PB_i(\phi_i, \lambda_{i}; \beta) .$$
\end{lemma}

\begin{proof}[Proof of Lemma~\ref{lem_PBi_continuity}]
For any~$(\phi_i, \lambda_{i}, \beta) \in [0, \phi_h] \times [0, \Lambda] \times [0, \infty),$ consider a sequence $(\phi_{i,m}, \lambda_{i,m}, \beta_m) \rightarrow (\phi_i, \lambda_{i}, \beta).$ From~\eqref{Eq_overall_exact},
\begin{align*}
\PB_i(\phi_{i,m}, \lambda_{i,m}; \beta_m) &= \DA_i(\phi_{i,m}, \lambda_{i,m}; \beta_m) f(\phi_{i, m}) +(1-f(\phi_{i, m})) \\
&= 1 - f(\phi_{i, m}) \left( 1 - \frac{1}{\sum_{n=0}^{\infty} \prod_{a=1}^n\frac{e}{\lambda_{i,m} f(\phi_{i, m}) + a \beta_m}}
\right).
\end{align*}
From the above equation, it suffices to show the continuity for the infinite series. 
We prove this continuity in the two cases similar to the proof of the Lemma~\ref{lem_MR_joint_cont}. The only difference is that here we have~$\lambda_{i, m} f(\phi_{i,m}) \to \lambda_{i} f(\phi_{i})$, and in Lemma~\ref{lem_MR_joint_cont}, we had~$h(\VPhi_m) \to h(\VPhi).$ As~$h$ and~$f$ are both continuous and bounded functions, the proof follows via exactly the same arguments as in the proof of Lemma~\ref{lem_MR_joint_cont}; we omit the details.
\end{proof}
}

\begin{lemma}[$\MR$ continuity in~$\PB$ system]\label{lem_mr_continuity_B_system}
Let~$\mathcal{F}:= [0, \phi_h]^2 \times [0, \infty),$ and~$\mathcal{H}:= \{ (\phi_1, \phi_2, 0) : \phi_1 = \phi_2, \phi_1 > f^{-1}(2\rho) \} .$ Then the function~$\MR_i$ is continuous on the set~$\mathcal{F} \setminus \mathcal{H}$. 
\end{lemma}
\begin{proof}[Proof of Lemma~\ref{lem_mr_continuity_B_system}]
For any~$(\VPhi, \beta) \in \mathcal{F} \setminus \mathcal{H}$, consider a sequence~$(\VPhi_m,\beta_m) \to (\VPhi, \beta)$; from Lemma~\ref{lem_mr_derivation}, the revenue rate of platform~$i$ is,
\begin{align}
\MR_i(\VPhi_m; \beta_m) 
&= \lambda_i(\VPhi_m; \beta_m) f(\phi_{i,m})\phi_{i,m} \left( 1 - \DA_i(\VPhi_m; \beta_m) \right) \nonumber \\
&= \lambda_i(\VPhi_m; \beta_m) f(\phi_{i,m}) \phi_{i,m}
\left( 1 - \frac{1}{\sum_{n=0}^{\infty} \mu_{i,n}(\VPhi_m, \lambda_{i,m}(\VPhi_m;\beta_m);\beta_m)} \right) . \label{eqn_mr_PB_cont_proof}
\end{align}
We prove this Lemma in the following two cases:~$\beta > 0$ and~$\beta = 0$.

\noindent \textbf{Case~1:} If~$\beta > 0$: From Lemma~\ref{lem_PB_joint_continuity_from_scratch}, the function~$\PB_i$ is jointly continuous on~$[0, \phi_h] \times [0, \Lambda] \times [0, \infty).$  For~$\beta > 0,$ by Lemma~\ref{lem_exist_unique_WE}, there exists a unique WE (optimizer)~$\lambda_i(\VPhi; \beta) \in [0, \Lambda].$ Using the Maximum Theorem~\cite{maximum}, we deduce the continuity of~$\lambda_i$ with respect to~$(\VPhi, \beta).$
Finally, using Lemma~\ref{lem_PB_joint_continuity_from_scratch},~$\sum_{n=0}^{\infty} \mu_{i,n}$ given in \eqref{eqn_mr_PB_cont_proof} converges uniformly for any~$(\VPhi, \beta) \in \mathcal{F}$ with~$\beta > 0$; which implies the joint continuity of the function~$\MR_i$ with respect to~$(\VPhi, \beta).$

\noindent \textbf{Case~2:} If~$\beta = 0$ (\idp\ regime): In this case, from the proof of Theorem~\ref{thm_WE_MR_PB}, there exist unique WE~$\lambda_i(\VPhi; 0)$ in the \idp\ regime in some cases; while in others there are more than one.
When the WE is unique in the \idp\ regime, again similar to the previous case, using the Maximum Theorem~\cite{maximum}, the proof follows. From the proof of Theorem~\ref{thm_WE_MR_PB}, there are only two cases when the WE is not unique: Case~(a)~if~$(\VPhi, \beta) \in \mathcal{H}$ (here we do not claim continuity), and Case~(b) if~$\phi_1 \geq f^{-1}(\rho),$ and~$\phi_1 > \phi_2.$ We now prove the continuity of~$\MR_i$ in Case~(b). 

\noindent Case~(b): We initially show the joint continuity of~$\lambda_1$ and using it we will show the joint continuity of~$\MR_1.$ Since $\lambda_2 = \Lambda - \lambda_1,$ it suffices to show the joint continuity of $\lambda_1,$ and $\MR_1.$ Consider the sequence~$(\VPhi_m,\beta_m) \to (\VPhi, 0)$ and corresponding $\lim_{(\VPhi_m,\beta_m) \to (\VPhi, 0)} \lambda_1(\VPhi_m; \beta_m)$. We will show~$\lambda_1(\VPhi_m; \beta_m) \to \lambda_1(\VPhi; 0)$.
From Theorem~\ref{thm_WE_MR_PB}, equation \eqref{Eq_overall_aproxx},
$$(\lamol, \lamtl) = \left(0, \Lambda \right), \text{ and } \PB_1(\VPhi; 0, 0) = 1 - f(\phi_1) > 1 - f(\phi_2) = \PB_2(\VPhi; \Lambda, 0).$$
Consider any~$\epsilon > 0$ with~$2\epsilon < f(\phi_2) - f(\phi_1).$ By Lemma~\ref{lem_PB_joint_continuity_from_scratch}, there exists a~$\delta_\epsilon > 0$ such that: 
\begin{itemize}
    \item~$\Lambda - \nicefrac{e}{f(f^{-1}(\rho) - \delta_\epsilon)} < \epsilon$ (see Theorem~\ref{thm_WE_MR_PB}), and
    \item for~$| (\tilde \VPhi, \tilde \beta) - (\VPhi, 0)|< \delta_\epsilon,$ we have~$|\PB_1(\tilde \VPhi; 0, \tilde \beta) - \PB_1(\VPhi; 0, 0)| < \epsilon,$ and~$|\PB_2(\tilde \VPhi; \Lambda, \tilde \beta) - \PB_2(\VPhi; \Lambda, 0)| < \epsilon.$ 
\end{itemize}
Thus for~$| (\tilde \VPhi, \tilde \beta) - (\VPhi, 0)|< \delta_\epsilon,$ 
$$
\PB_1(\tilde \VPhi; 0, \tilde \beta) > \PB_2(\tilde \VPhi; \Lambda, \tilde \beta).
$$
For~$\tilde \beta > 0$ and~$| (\tilde \VPhi, \tilde \beta) - (\VPhi, 0)|< \delta_\epsilon,$ the function~$\PB_1$ is a strictly increasing in~$\lambda_1$, while~$\PB_2$ is a strictly decreasing in~$\lambda_1$; thus we have~$\lambda_1(\tilde \VPhi; \tilde \beta) = 0.$
\noindent For~$\tilde \beta  = 0$ and~$| (\tilde \VPhi, \tilde \beta) - (\VPhi, 0)|< \delta_\epsilon,$
\begin{itemize}
    \item If~$\tilde \phi_1 \geq f^{-1}(\rho),$ it follows from Table~\ref{tab:WE_MR_PB} that~$\lambda_1(\tilde \VPhi; 0) = 0.$
    \item If~$\tilde \phi_1 < f^{-1}(\rho),$ it follows from Table~\ref{tab:WE_MR_PB} that~$\lambda_1(\tilde \VPhi; 0) =  \Lambda - \nicefrac{e}{f(\tilde \phi_1)}.$ Note that~$\tilde \phi_1 > f^{-1}(\rho) - \delta_\epsilon,$ that is,~$f(\tilde \phi_1 ) < f(f^{-1}(\rho) - \delta_\epsilon).~$ Thus,
    $$\lambda_1(\tilde \VPhi; 0) = \Lambda - \frac{e}{f(\tilde \phi_1)} <  \Lambda - \frac{e}{f(f^{-1}(\rho) - \delta_\epsilon)} < \epsilon.$$
\end{itemize}
The above arguments imply that~$\lambda_1(\VPhi_m; \beta_m) \to \lambda_1(\VPhi; 0)$. Therefore, using Lemma~\ref{lem_PB_joint_continuity_from_scratch}, it can be proved that~$\sum_{n=0}^\infty \mu_{i,n}(\VPhi_m, \lambda_{i,m}(\VPhi_m;\beta_m), \beta_m) \rightarrow \infty.$
Thus we have $\MR_i(\VPhi_m,\beta_m) \to \lambda_i(\VPhi; 0) f(\phi_{i}) \phi_{i} = \MR_i(\VPhi;0).$
\end{proof}

\ignore{
\begin{lemma}[$\MR$ continuity in~$\mathcal{\WB}$ system]\label{lem_mr_continuity_Wait_system}
Let~$\mathcal{F}:= [0, \phi_h]^2 \times [0, \infty) \times [0, \gamma),$ where $\gamma = \min\{ \gamma_1, \gamma_2\},$ and~$\mathcal{H}:= \{ (\phi_1, \phi_2, 0, \alpha) : \phi_1 = \phi_2, \phi_1 > f^{-1}(2\rho) \} .$ Then the function~$\MR_i$ is continuous on the set $\mathcal{F} \setminus \mathcal{H}$. 
\end{lemma}
\begin{proof}[Proof of Lemma~\ref{lem_mr_continuity_Wait_system}]

For any~$(\VPhi, \beta, \alpha) \in \mathcal{F} \setminus \mathcal{H}$, consider a sequence~$(\VPhi_m,\beta_m, \alpha_m) \to (\VPhi, \beta, \alpha)$; from Lemma~\ref{lem_mr_derivation}, the matching revenue of platform~$i$ is,
\begin{align*}
\MR_i(\VPhi_m; \beta_m, \alpha_m) 
&= \lambda_i(\VPhi_m; \beta_m, \alpha_m) f(\phi_{i,m})\phi_{i,m} \left( 1 - \DA_i(\VPhi_m; \beta_m, \alpha_m) \right)  \\
&= \lambda_i(\VPhi_m; \beta_m, \alpha_m) f(\phi_{i,m}) \phi_{i,m}
\left( 1 - \frac{1}{\sum_{n=0}^{\infty} \mu_{i,n}(\VPhi_m(\lambda_{i,m}(\VPhi_m;\beta_m, \alpha_m);\beta_m, \alpha_m)} \right) .
\end{align*}
Towards this, we first consider the joint continuity of WE $\lambda_i$.
From Lemma~\ref{lem_PB_joint_continuity_from_scratch}, the function~$\PB_i$ is jointly continuous on~$[0, \phi_h] \times [0, \Lambda] \times [0, \infty)$,  using similar arguments,  the function~$\WB_i$ is jointly continuous on~$[0, \phi_h] \times [0, \Lambda] \times [0, \infty) \times [0, \gamma)$  (see \eqref{eqn_QoS_WB_positive_beta}-\eqref{eqn_QoS_WB_zero_beta}). Because of joint continuity in $\WB_i,$ using Maximum theorem based arguments similar to the proof of Theorem~\ref{thm_WE_MR_PB}, we have upper semi-continuity of the set of WEs w.r.t. $(\VPhi, \beta, \alpha)$. By Theorem \ref{thm_WE_MR_new}, the WE is unique when $\beta > 0$, thus we have continuity of WE $(\lambda_i)$ at all $(\VPhi, \beta, \alpha)$ with $\beta > 0$ (see proof of Theorem~\ref{thm_WE_MR_PB}).  
Further, using Lemma~\ref{lem_PB_joint_continuity_from_scratch},~$\sum_{n=0}^{\infty} \mu_{i,n}$ converges uniformly for any~$(\VPhi, \beta, \alpha) \in \mathcal{F}$ with~$\beta > 0$; which implies the joint continuity of the function~$\MR_i$ with respect to~$(\VPhi, \beta, \alpha)$ when $\beta  > 0.$




Next, we consider a case when~$\beta = 0$. Recall from the proof of Theorem~\ref{thm_WE_MR_PB}, at $\beta = 0,$ there exist unique WE~$\lambda_i(\VPhi; 0, 0)$ in some cases; while in others there are more than one.
When the WE is unique, again similar to the previous case, using the Maximum Theorem~\cite{maximum}, the proof follows. From the proof of Theorem~\ref{thm_WE_MR_PB}, there are only two cases when the WE is not unique: Case~(a)~if~$(\VPhi, \beta, \alpha) \in \mathcal{H}$ (here we do not claim continuity), and Case~(b) if~$(\VPhi, \beta, \alpha) \in \mathcal{F}$ such that $\phi_1 \geq f^{-1}(\rho),$ and~$\phi_1 > \phi_2.$ We now prove the continuity of~$\MR_i$ in Case~(b). 

\noindent Case~(b): We initially show the joint continuity of~$\lambda_i$ and using it we will show the joint continuity of~$\MR_i.$ Consider the sequence~$(\VPhi_m,\beta_m, \alpha_m) \to (\VPhi, 0, \alpha)$ and corresponding~$\lim_{(\VPhi_m,\beta_m, \alpha_m) \to (\VPhi, 0, \alpha)} \lambda_1(\VPhi_m; \beta_m, \alpha_m)$. We will show~$\lambda_1(\VPhi_m; \beta_m, \alpha_m) \to \lambda_1(\VPhi; 0, \alpha)$.
From the proof of Theorem~\ref{thm_WE_MR_PB} and~\eqref{eqn_QoS_WB_zero_beta},
$$(\lambda_1(\VPhi; 0, \alpha), \lambda_2(\VPhi; 0, \alpha)) = \left(0, \Lambda \right), \text{ and } \PB_1(\VPhi; 0, 0, \alpha) = 1 - f(\phi_1) > 1 - f(\phi_2) = \PB_2(\VPhi; \Lambda, 0, \alpha).$$
Consider any~$\epsilon > 0$ with~$2\epsilon < f(\phi_2) - f(\phi_1).$ By Lemma~\ref{lem_PB_joint_continuity_from_scratch}, there exists a~$\delta_\epsilon > 0$ such that: 
\begin{itemize}
    \item~$\Lambda - \nicefrac{e}{f(f^{-1}(\rho) - \delta_\epsilon)} < \epsilon,$ and
    \item for~$| (\tilde \VPhi, \tilde \beta, \tilde \alpha) - (\VPhi, 0, \alpha)|< \delta_\epsilon,$ we have~$|\PB_1(\tilde \VPhi; 0, \tilde \beta, \tilde \alpha) - \PB_1(\VPhi; 0, 0, \alpha)| < \epsilon,$ and~$|\PB_2(\tilde \VPhi; \Lambda, \tilde \beta, \tilde \alpha) - \PB_2(\VPhi; \Lambda, 0, \alpha)| < \epsilon.$ 
\end{itemize}
Thus for~$| (\tilde \VPhi, \tilde \beta, \tilde \alpha) - (\VPhi, 0, \alpha)|< \delta_\epsilon,$ 
$$
\PB_1(\tilde \VPhi; 0, \tilde \beta, \tilde \alpha) > \PB_2(\tilde \VPhi; \Lambda, \tilde \beta, \tilde \alpha).
$$
For~$\tilde \beta > 0$ and~$| (\tilde \VPhi, \tilde \beta, \tilde \alpha) - (\VPhi, 0, \alpha)|< \delta_\epsilon,$ the function~$\PB_1$ is a strictly increasing in~$\lambda_1$, while~$\PB_2$ is a strictly decreasing in~$\lambda_1$; thus we have~$\lambda_1(\tilde \VPhi; \tilde \beta, \tilde \alpha) = 0.$
\noindent For~$\tilde \beta  = 0$ and~$| (\tilde \VPhi, \tilde \beta, \tilde \alpha) - (\VPhi, 0, \alpha)|< \delta_\epsilon,$ 
\begin{itemize}
    \item If~$\tilde \phi_1 \geq f^{-1}(\rho),$ it follows from Table~\ref{tab:WE_MR_PB} that~$\lambda_1(\tilde \VPhi; 0, \tilde \alpha) = 0.$
    \item If~$\tilde \phi_1 < f^{-1}(\rho),$ it follows from Table~\ref{tab:WE_MR_PB} that~$\lambda_1(\tilde \VPhi; 0, \tilde \alpha) =  \Lambda - \nicefrac{e}{f(\tilde \phi_1)}.$ Note that~$\tilde \phi_1 > f^{-1}(\rho) - \delta_\epsilon,$ that is,~$f(\tilde \phi_1 ) < f(f^{-1}(\rho) - \delta_\epsilon).~$ Thus,
    $$\lambda_1(\tilde \VPhi; 0) = \Lambda - \frac{e}{f(\tilde \phi_1)} <  \Lambda - \frac{e}{f(f^{-1}(\rho) - \delta_\epsilon)} < \epsilon.$$
\end{itemize}
The above arguments imply that~$\lambda_1(\VPhi_m; \beta_m, \alpha_m) \to \lambda_1(\VPhi; 0, \alpha)$. 

Finally, using Lemma~\ref{lem_PB_joint_continuity_from_scratch}, it can be proved that $\sum_{n=0}^\infty \mu_{i,n}(\VPhi_m, \lambda_{i,m}(\VPhi_m;\beta_m, \alpha_m);\beta_m, \alpha_m) \rightarrow \infty.$
Thus, we have $\MR_i(\VPhi_m;\beta_m, \alpha_m) \to \lambda_i(\VPhi; 0, \alpha) f(\phi_{i}) \phi_{i} = \MR_i(\VPhi;0,\alpha).$
\end{proof}
}

\begin{lemma}\label{lem_eta_positive_epsilon_EC_proof}
If~$m$,~$\eL,$~$\eU$ are as defined in Theorem~\ref{thm_B_sys_limit_sys_main}, and~$\delta_\epsilon$ is as defined in~\eqref{eqn_positive_delta_EC}, then
{
\begin{align*}
\min\Bigg \{e, (m(\eU) - m(\eU + \epsilon)) , (m(\eU) - m(\eU - \epsilon)), \delta_\epsilon, \left(\frac{e \eL}{2} - \frac{\Lambda f(\eL) \eL}{4}\right) &, \\
 \left(m(\eU) -\frac{\Lambda f(\eL)\eL}{2e} \right) &\Bigg\}  > 0.    
\end{align*}}
\end{lemma}
\begin{proof}[Proof of Lemma~\ref{lem_eta_positive_epsilon_EC_proof}]
Observe that~$e>0$, and from~\eqref{eqn_positive_delta_EC},~$ \delta_\epsilon > 0.$ As~$m$ is strictly concave and~$\eU$ is the maximizer of~$m,$ we have~$m(\eU) - m(\eU + \epsilon) > 0,$ and~$m(\eU) - m(\eU - \epsilon) > 0.$ From Lemma~\ref{lemma_oscillatory},~$e\eL - \nicefrac{\Lambda}{2}f(\eL)\eL > 0 ,$ and by definition~$\eL = m(\eU),$ thus~$m(\eU) -\nicefrac{\Lambda f(\eL)\eL}{2e} > 0.$
\end{proof}

We state the following lemma, whose proof is immediate.
\begin{lemma}\label{lemma_oscillatory}
If~$\phi \in \left(f^{-1}\left(2 \rho\right), \phi_h \right)$ then~$e\phi > \frac{\Lambda}{2}f(\phi)\phi > \Lambda f(\phi)\phi - e\phi = m(\phi)$.
\end{lemma}

\begin{lemma}\label{lem_coupling_arguments}
Consider two identical systems $S$ and $\tilde S$ with the only difference being in the passenger arrival rates $\lambda$ and $\tilde{\lambda}$ respectively, such that $\lambda > \tilde{\lambda}$. Then (i)
\begin{equation} \label{eqn_coup_claim}
n(t) \leq \tilde{n}(t) \text{ and } Z(t) \leq \tilde{Z}(t) \text{ at all time $t$, where } Z(t) = n(t) + r(t).
\end{equation}
(ii) Let $\pi^\lambda$ be the stationary distribution of the system $( n(t), r(t) )$  with passenger arrival rate $\lambda$ and let $f$ be a strict monotone increasing function. Then the following is a strictly monotone increasing function:
$$
\lambda \mapsto \sum_{n=0}^{\infty} f(n) \pi^\lambda (n) .
$$
Part~(ii)
 is also true if $f(0) > f(1).$
 
 \end{lemma}
\begin{proof}[Proof of Lemma~\ref{lem_coupling_arguments}]
We use coupling arguments to compare the sample paths of the two systems. The $S$ system sees all passenger arrivals (realization of Poisson process at rate $ \lambda$), while the $\tilde S$ system sees only a fraction of them sampled independently with probability $\nicefrac{\tilde \lambda}{ \lambda}$; thus one can have a passenger arrival in both the systems or in $S$ system and not in $\tilde S$ system. Both systems see the driver arrivals at the same time.  Let $n(t), r(t)$ represent the number of drivers waiting and riding in the $S$ system at time $t$; the pair $\tilde n(t), \tilde r(t)$ represents similar quantities in $\tilde S$ system. The residual quantities like residual inter-arrival times, residual impatient times, etc., (at any $t$ or after any event) are again exponentially distributed with the same parameters as the corresponding initial quantities. Recall also that all these parameters are the same in both the systems except for $\lambda$ and $\tilde \lambda$. Based on these properties, we further use the following coupling construction arguments:

\begin{itemize}
    \item For example, if $n(t) > \tilde n(t),$ at some $t$ (one can also take $t = \tau^+,$ where $\tau$ is an event, e.g., passenger arrival), we consider $n(t)$ independent copies of exponentially distributed random variables with parameter $\beta$ to represent (residual) impatient times of the waiting drivers in $S$ system; without loss of generality, we consider the first $\tilde n(t)$ of them to represent the impatient times in $\tilde S$ system; the vice versa arguments are used if $\tilde n(t) > n(t) $ and the same quantities are used when $\tilde n(t) = n(t) .$ 
    \item Similar construction arguments are used for riding times, based on $r(t), \tilde r(t).$
    \item If the drivers complete the rides simultaneously in both the systems, the random joining back flags (the probability of which equals  $p$) are also coupled. 
    
\end{itemize}
%
%
%
Based on the above coupling arguments, one can observe only one of the following joint events in the combined (coupled) $(S, \tilde S)$ system. 
\begin{enumerate}
    \item Passenger arrival to both the systems
    \item Passenger arrival to only $S$ system
    \item Driver arrival to both the systems 
    \item Driver impatience event either in one of the systems or in both
    \item Ride completion of a driver either in one of the systems or in both
\end{enumerate}
Let $t^- \geq 0$, in general, represent the time epoch just before any of the above events; to have uniform notation, it also represents the initial time epoch $0.$  \textit{By hypothesis, the system state satisfies the condition~\eqref{eqn_coup_claim} at $t^- = 0$. We assume it to satisfy \eqref{eqn_coup_claim} at any~$t^- $ and  prove that the system at corresponding $t^+$ (time instance just after the event) also satisfies the same, which completes the proof.} Observe that the state at~$t^+$ equal that at time $t$, when one considers right continuous sample paths. 
\ignore{\begin{enumerate}[(a)]
    \item $n(t^-) <  \tilde{n}(t^-)$ and $Z(t^-) = \tilde{Z}(t^-)$
    \item $n(t^-) = \tilde{n}(t^-)$ and $Z(t^-) = \tilde{Z}(t^-)$
    \item $n(t^-) <  \tilde{n}(t^-)$ and $Z(t^-) <  \tilde{Z}(t^-)$
    \item $n(t^-) = \tilde{n}(t^-)$ and $Z(t^-) <  \tilde{Z}(t^-)$
\end{enumerate}}
To begin with, observe that at any change epoch time $t$, only one event (e.g., passenger arrival) occurs.\footnote{two events occurring simultaneously have probability~$0$ in continuous random variables.} 
We obtain the required proof for each of the above mentioned five events in the following. Without loss of generality, we only consider the case when passenger arrivals accept the price (by replacing $\lambda f(\phi)$ with $\lambda$, and $\tilde \lambda f(\phi)$ with $\tilde \lambda$).
\begin{enumerate}
    \item Passenger arrival to both the systems:

    If $n(t^-) > 0$ then using \eqref{eqn_coup_claim}, $\tilde n(t^-) > 0$; at such $t,$ both the systems accommodate the arriving passenger and hence the waiting drivers reduce by one, and thus $n(t) = n(t^-) - 1$ and $\tilde n(t) = \tilde n(t^-) - 1$.
    
    On the other hand, if $n(t^-) = 0 $ and $ \tilde n( t^-) > 0$ then   only $\tilde S$ system accommodates the passenger  and hence $n(t) = n(t^-)$ and  $\tilde n(t) = \tilde n(t^-) - 1$. 

    Finally, if $n(t^-) = \tilde n(t^-) = 0$, there is no change in the system state.
    
    In all the cases, $Z(t) = Z(t^-)$ and  $\tilde Z(t) = \tilde Z(t^-)$, thus~\eqref{eqn_coup_claim} is satisfied at $t$; observe that $Z$ or $\tilde Z$ changes only when a new driver arrives to the system or when a driver leaves the system after completing the ride.

 \ignore{\begin{enumerate}
        \item If $n(t^-) = \tilde n(t^-) = 0$:
        This is a trivial case as none of the systems can accept the passenger, and the system state does not change.
        \item If $n(t^-) = 0 $ and $ \tilde n( t^-) > 0$:
        In this case, at time $t,$ we have $n(t) = n(t^-)$, $\tilde n(t) = \tilde n(t^-) - 1$ and $Z(t) = Z(t^-)$ and  $\tilde Z(t) = \tilde Z(t^-),$ thus the claim in \eqref{eqn_coup_claim} is satisfied.
        \item If $n(t^-) > 0$: In this case, at time $t,$ we have $n(t) = n(t^-) - 1$, $\tilde n(t) = \tilde n(t^-) - 1$ and $Z(t) = Z(t^-)$ and  $\tilde Z(t) = \tilde Z(t^-),$ thus the claim in \eqref{eqn_coup_claim} is satisfied.
    \end{enumerate}}
    %
    %
    \item Passenger arrival to only $S$ system:
    
    If $n(t^-) > 0$ then at time $t,$ the $S$ system accommodate the arriving passenger and thus we have $n(t) = n(t^-) - 1$, else $n(t) = n(t^-) = 0$,  and the system state does not change.    
    In both the cases, $Z(t) = Z(t^-)$ and  $\tilde Z(t) = \tilde Z(t^-),$ thus~\eqref{eqn_coup_claim} is satisfied at $t$.
\ignore{
    \begin{enumerate}
        \item If $n(t^-) = 0$:
        This is a trivial case as the $S$ system cannot accept the passenger, and the system state does not change.
        \item If $n(t^-) > 0$: In this case, at time $t,$ we have $n(t) = n(t^-) - 1$, and $Z(t) = Z(t^-)$ and  $\tilde Z(t) = \tilde Z(t^-),$ thus the claim in \eqref{eqn_coup_claim} is satisfied.
    \end{enumerate}
    }
    \item Driver arrival to both the systems:

    Recall by coupling, drivers arrive simultaneously to both the systems. 
    At such~$t$, we have $n(t) = n(t^-) + 1$, $\tilde n(t) = \tilde n(t^-) + 1$ and $Z(t) = Z(t^-) + 1$ and  $\tilde Z(t) = \tilde Z(t^-) + 1.$ Thus, if \eqref{eqn_coup_claim} is satisfied at $t^-$ then it is also satisfied at $t$.
    \item Driver impatience event either in one of the systems or in both:
    
    Firstly, observe that such an event is possible only when the corresponding $n(t^-)$ or $ \tilde n(t^-)$ is strictly positive.
    
    If $n(t^-) = \tilde n(t^-)$, once again, by coupling, a driver gets impatient in both systems.
    Therefore, at such $t$, we have $n(t) = n(t^-) - 1$, $\tilde n(t) = \tilde n(t^-) - 1$. As before, there is no change in $Z$, $\tilde Z$.

    If $n(t^-) < \tilde n(t^-)$ then the $\tilde S$ system has $ \tilde n(t^-) - n(t^-)$ more drivers at time~$t^-.$ Hence, a driver can get impatient in both the systems, or a driver in $\tilde S$ gets impatient and not in the $S$ system. 
    In the latter case, at such $t$, we have $n(t) = n(t^-) $, $Z(t^-) = Z(t)$  and $\tilde n(t) = \tilde n(t^-) - 1$, $\tilde Z(t^-) = \tilde Z(t)$. In both cases, it is easy to observe~\eqref{eqn_coup_claim} is satisfied at $t$. 
 
\ignore{\begin{enumerate}
        \item If $n(t^-) = \tilde n(t^-)$:
        Since the number of waiting drivers are same and all residual times are coupled, we have $n(t) = n(t^-) - 1$, $\tilde n(t) = \tilde n(t^-) - 1$, $Z(t) = Z(t^-)$ and  $\tilde Z(t) = \tilde Z(t^-),$ thus the claim in \eqref{eqn_coup_claim} is satisfied.
        \item If $n(t^-) < \tilde n(t^-)$: In this case, the $tilde S$ system has $ \tilde n(t^-) - n(t^-)$ more drivers at time $t^-.$ If the drivers from $S$ system looses patience then driver from $\tilde S$ system also looses patience and  we have $n(t) = n(t^-) - 1$, and thus the claim in \eqref{eqn_coup_claim} is satisfied. If the driver from only $\tilde S$ system, which is not coupled with $S$ system looses patience then also, $\tilde n(t) = \tilde n(t^-) - 1$, which implies that $n(t) \leq \tilde n(t)$,       
        $Z(t) = Z(t^-)$ and  $\tilde Z(t) = \tilde Z(t^-),$ thus the claim in \eqref{eqn_coup_claim} is satisfied.
    \end{enumerate}
    }
    \item Ride completion of a driver either in one of the systems or in both: 

This event is not possible when  $r(t^-) = \tilde r(t^-)  = 0$, i.e., without riding drivers.


\ul{If $r(t^-) = \tilde r(t^-) > 0$}, then by coupling arguments again, a driver completes the ride simultaneously in both the systems; the drivers either leave the systems or join back, this flag is also coupled. If the drivers rejoin then, 
\begin{eqnarray*}
    \big  (n(t),  \ r(t)) \ = \ (n(t^-) + 1, \  r(t^-) - 1  \big ) \text{ and } \big ( \tilde n(t), \ \tilde r(t) \big )  =  \big ( \tilde n(t^-) + 1,  \  \tilde r(t^-) - 1  \big  ) \\
\mbox{ and hence, }      Z(t) \ = \ Z(t^-),  \mbox{ and }   \tilde Z(t) = \tilde Z(t^-). \hspace{50mm}
\end{eqnarray*}
If they leave the system,   we have 
\begin{eqnarray*}
    \big  (n(t),  \ r(t)) \ = \ (n(t^-)  \  r(t^-) - 1  \big ) \text{ and } \big ( \tilde n(t), \ \tilde r(t) \big )  =  \big ( \tilde n(t^-) ,  \  \tilde r(t^-) - 1  \big  ) \\
\mbox{ and then, }      Z(t) \ = \ Z(t^-) - 1,  \   \tilde Z(t) = \tilde Z(t^-) -1.  \hspace{40mm}
\end{eqnarray*}
In both the cases, \eqref{eqn_coup_claim} is satisfied at $t$.

\ul{If $r(t^-) < \tilde r(t^-)$} then the system $\tilde S$ has $\tilde r(t^-) - r(t^-)$ more drivers in ride than the system $S$ at time $t^-$. Therefore, a driver can finish the ride in both the systems, or a driver in the~$\tilde S$ system can finish the ride but not in the system $S$.  In the latter case, at such $t,$

\begin{itemize}
    \item if the driver rejoins the system $\tilde S$, then we have 
$$
\left( n(t), r(t) \right)  = \left( n(t^-),  \tilde r(t^-) - 1 \right )
\text{ and }
\left( \tilde n(t), \tilde r(t) \right)  = \left( \tilde n(t^-) + 1, \tilde r(t^-) - 1 \right), 
 $$
    leading to $Z(t) = Z(t^-)$ and $\tilde Z(t) = \tilde Z(t^-)$. 
    
    \item on the other hand if the driver leaves $\tilde S$, we will have $Z(t) = Z(t^-1)$, but  $\tilde Z(t) = \tilde Z(t^-) - 1$; nonetheless \eqref{eqn_coup_claim} is satisfied at $t$ because:
    $$
    \tilde Z(t^-) - Z(t^-) =  \tilde n(t^-) - n(t^-) + \tilde r(t^-) - r(t^-) \ge \tilde r(t^-) - r(t^-) \ge 1,
    $$in this sub-case.

    %
\end{itemize}

\ul{If $r(t^-) > \tilde r(t^-)$} then system $S$ has $r(t^-) - \tilde r(t^-)$ more drivers  at $t^-$.  We will either have 
$$
\tilde n(t) - n(t)  = \tilde n(t^-) - n(t^-) \mbox{ with } \tilde Z(t) -  Z(t) = \tilde Z(t^-) - Z(t^-).
$$
or  
$$
\tilde n(t) - n(t) < \tilde n(t^-) - n(t^-) \mbox{ with } \tilde Z(t) -  Z(t) = \tilde Z(t^-) - Z(t^-).
$$
Nonetheless, the inequality~\eqref{eqn_coup_claim}  is  again satisfied at $t$ (using similar logic as above) because we now have: 
    $$
    \tilde n(t^-) - n(t^-) = \tilde Z(t^-) - Z(t^-) - \left (\tilde r(t^-) - r(t^-) \right ) \ge r(t^-) - \tilde  r(t^-) \ge 1.
    $$

\ignore{If the driver from only system $S$ completes the ride, then
\begin{itemize}
    \item if the driver rejoins the system $ S$, then we have $n(t) = n(t^-) + 1$, $\tilde n(t) = \tilde n(t^-)$ and $r(t) = r(t^-) -1 $, $\tilde r(t) = \tilde r(t^-)$ which implies that at $t,$ we have $\tilde  r(t) - r(t) \geq  n(t) - \tilde  n(t)$ and thus, we have $Z(t) \leq \tilde  Z(t)$ and thus \eqref{eqn_coup_claim} is satisfied at $t$.
    \item if the driver leaves the system $S$, then we have $n(t) = n(t^-)$, $\tilde n(t) = \tilde n(t^-)$ and $r(t) = r(t^-) - 1$, $\tilde r(t) = \tilde r(t^-)$, which implies that $Z(t) \leq \tilde Z(t),$ and thus \eqref{eqn_coup_claim} is satisfied at $t$.
\end{itemize}
     }

\ignore{\begin{enumerate}
     \item If $r(t^-) = \tilde r(t^-):$ Since all riding the drivers are coupled, the ride-completing drivers complete the ride simultaneously in both the systems; either both rejoin or both leave the systems simultaneously. Thus, again, the claim in \eqref{eqn_coup_claim} is satisfied at $t$.
     
     \item If $r(t^-) < \tilde r(t^-):$ In this case, system $\tilde S$ has $\tilde r(t^-) - r(t^-)$ more drivers in ride. If drivers from system $S$ complete the ride, then the driver from system $\tilde S$ also completes the ride as the riding drivers in both systems are coupled. If the driver from $\tilde S$ system only completes the ride (out of $\tilde r(t^-) - r(t^-)$ drivers), then again the claim in \eqref{eqn_coup_claim} is satisfied at $t$.
     
     \item If $r(t^-) > \tilde r(t^-):$
     Again, if the drivers from both systems complete the ride, then the claim in \eqref{eqn_coup_claim} is satisfied as both drivers rejoin back or leave the respective systems simultaneously. If the driver from only $S$ system looses patience then $r(t) = r(t^-) -1 ,$ $n(t) = n(t^-) + 1,$ and $\tilde r(t) = \tilde r(t^-) ,$ $\tilde n(t) = \tilde n(t^-),$ thus using \eqref{eqn_claim_equivalence}, the claim in \eqref{eqn_coup_claim} is satisfied at $t$.
     
\end{enumerate}
}
\end{enumerate}

\ignore{Now observe that at any time t, 
\begin{equation} \label{eqn_claim_equivalence}
 Z(t) \leq \tilde Z(t) \implies
 r(t) - \tilde r(t) \leq \tilde n(t) - n(t) \text{ OR } 
\tilde  r(t) - r(t) \geq  n(t) - \tilde  n(t).
\end{equation}
}

Thus, \eqref{eqn_coup_claim} is satisfied for all $t$ in almost all the sample paths. This completes the proof of part~(i).

Towards part (ii), 
\ignore{using part (i) and standard text book arguments (almost surely):
$$
\sum_{n=0}^{\infty} f(n) \pi^\lambda (n) = \lim_{T \to \infty}  \frac{1}{T} \int_0^T f( n(t) ) dt \le  \lim_{T \to \infty}  \frac{1}{T} \int_0^T f( \tilde n(t) ) dt   = \sum_{n=0} f(n) \pi^{\tilde \lambda } (n).
$$
Towards strict monontonicity, one can use RRT as in the proof of Lemma \ref{lem_mr_derivation} and obtain that 
$$
\lim_{T \to \infty}  \frac{1}{T} \int_0^T f( n(t) ) dt =  \frac{ E \left [ \int_0^{\tau (s)} f( n(s) ) ds \right ] } {E [\tau (s)]}
$$
where, the random variable $\tau$ denotes the length of one renewal cycle, $s = (0,0,0,0)$ is any arbitrary state that 

Thus, we established stochastic dominance with respect to $\lambda$.
}
consider a joint system which combines  $S$ and $\tilde S$ systems coupled as in part~(i), and where the state of the system  is given by ${\bf s}  = (n, r, \tilde n, \tilde r).$ Using part~(i) and standard textbook RRT based arguments (almost surely),
$$
\sum_{n=0}^{\infty} f(n) \pi^\lambda (n) = \lim_{T \to \infty}  \frac{1}{T} \int_0^T f( n(t) ) dt \le  \lim_{T \to \infty}  \frac{1}{T} \int_0^T f( \tilde n(t) ) dt   = \sum_{n=0} f(n) \pi^{\tilde \lambda } (n).
$$
{\bf Strict Monotonicity:}
Consider the special state ${\bf s}^0 := (0,0,0,0)$ and define the stopping time for the combined system:
$\tau({\bf s}^0) =  \inf_t \{ t > \tau_1 : {\bf s} (t) = {\bf s}^0 \}$, where $\tau_1$ is the time of the first transition epoch after the combined system 
starts in ${\bf s}^0$ state at time $0$. Observe that $\tau({\bf s}^0)  < \infty$ with probability one when $\tilde S$ system is stationary -- in fact it is easy to observe that  $\tau({\bf s}^0) = \inf_t \{ t > \tau_1 : \tilde Z (t) = 0   \} . $
 One can use RRT as in the proof of Lemma \ref{lem_mr_derivation} and obtain that almost surely
\begin{align*}
\lim_{T \to \infty}  \frac{1}{T} \int_0^T f( n(t) ) dt &=  \frac{ E_{{\bf s}^0} \left [ \int_0^{\tau ({\bf s}^0)} f( n(t) ) dt \right ] } {E_{{\bf s}^0} \left[ \tau ({\bf s}^0) \right]}
\mbox{ and} \\
\lim_{T \to \infty}  \frac{1}{T} \int_0^T f( \tilde n(t) ) dt & =  \frac{ E_{{\bf s}^0} \left [ \int_0^{\tau ({\bf s}^0)} f( \tilde n(t) ) dt \right ] } {E_{{\bf s}^0} \left[ \tau ({\bf s}^0) \right]}.    
\end{align*}
Towards strict monotonicity,
again as in Lemma \ref{lem_mr_derivation}, %
we prove that  $n (t) \le \tilde n (t) - 1$  for $t \in T$, where $T$ is some time interval  whose expected length is at least $\nicefrac{1}{(\nu + \tilde \lambda)}$, 
with strictly positive probability; this we prove for the case when $n(t) = 0$,  $\tilde n(t) >0 $ and hence the hypothesis is sufficient.  Basically, there exists a strictly positive probability at least equal to $ c (1 - \nicefrac{\tilde \lambda}{\lambda}) $ (with an appropriate $c >0$) such that: a) 
a driver enters both the systems at $\tau_1$; 
b) at the next event epoch, 
a passenger enters $S$ system but not  $\tilde S$ system; and c) and so the system remains with $n = 0 < \tilde n = 1 $  at least till the riding driver in $S$ system completes the ride and rejoins or till a passenger arrives to~$\tilde S$ system.

\begin{figure}[ht]
    \centering
    \includegraphics[trim = {0.2cm 2.2cm 0cm 3.2cm}, clip,
    scale = 0.38]{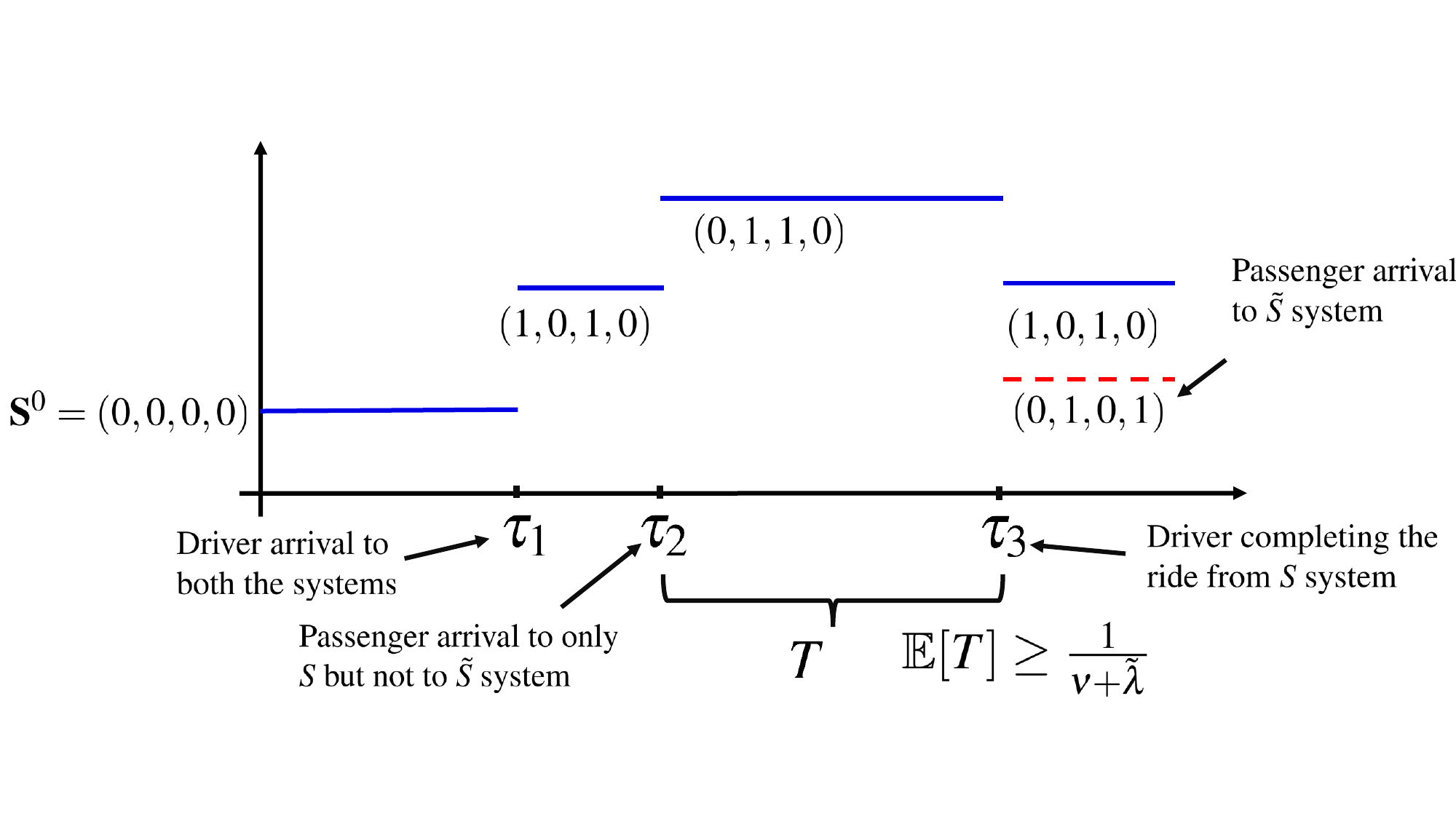}
{\caption{A sample path in the coupled system with initial state  ${\bf s}^0 = (0,0,0,0)$ and with $n (t) < \tilde n (t)$ for $t \in [\tau_2, \tau_3]$.}  
\label{fig:coupling_arguments}}
\end{figure}
\end{proof}
\ignore{\begin{enumerate}
    \item $n(t^-) > \tilde{n}(t^-)$ and $N(t^-) = \tilde{N}(t^-) \implies r(t^-) < \tilde{r}(t^-)$
\begin{itemize}
    \item Passenger arrival to both  the systems or only to $\tilde S$ system satisfies the claim.
    \item Loss of patience of a driver in $S$ system also satisfies the hypothesis (at maximum, the strict inequality can change into equality i.e., $n(t) = \tilde{n}(t)$).
    \item Ride completion (both rejoining as well as leaving the system) of the driver in $\tilde S $ system also satisfies the hypothesis (at maximum, the strict inequality can change into inequality $n(t) = \tilde{n}(t)$).
\end{itemize}

    \item $n(t^-) = \tilde{n}(t^-)$ and $N(t^-) = \tilde{N}(t^-) \implies r(t^-) = \tilde{r}(t^-)$
\begin{itemize}
    \item Passenger arrival to $\tilde S$ system can only change the equality into strict inequality ($n(t) > \tilde{n}(t)$), which still satisfies the hypothesis.
    \item Loss of patience of a driver or the ride completion (both rejoining as well as leaving the system) of the driver in both the systems satisfies the claim as all the drivers are coupled.
\end{itemize}

    \item $n(t^-) > \tilde{n}(t^-)$ and $N(t^-) > \tilde{N}(t^-) \implies r(t^-) > \tilde{r}(t^-)$ or $r(t^-) < \tilde{r}(t^-)$ or $r(t^-) = \tilde{r}(t^-)$
\begin{itemize}
    \item Passenger arrival to both systems or only to $\tilde S$ system satisfies the claim.
    \item Loss of patience of a driver in $S$ system also satisfies the hypothesis (at maximum, the strict inequality can change into equality i.e., $n(t) = \tilde{n}(t)$).
    \begin{itemize}
        \item If $r(t^-) > \tilde{r}(t^-)$: Ride completion (both rejoining as well as leaving the system) of the driver in $ S $ system also satisfies the hypothesis.
        \item $r(t^-) < \tilde{r}(t^-)$: Ride completion (both rejoining as well as leaving the system) of the driver in $ \tilde S $ system also satisfies the hypothesis (at maximum, the strict inequality can change into inequality $n(t) = \tilde{n}(t)$).
        \item $r(t^-) = \tilde{r}(t^-)$: Ride completion of a driver from both the platforms satisfy the claim as all the drivers in ride are coupled.
    \end{itemize}
\end{itemize}

    \item $n(t^-) = \tilde{n}(t^-)$ and $N(t^-) > \tilde{N}(t^-) \implies r(t^-) > \tilde{r}(t^-)$
\begin{itemize}
    \item Passenger arrival in $\tilde S $ system satisfies the claim.
    \item Loss of patience of driver satisfies the claim as all the waiting drivers are coupled.
    \item Ride completion (rejoining as well as leaving the system) of the driver in the~$ S $ system satisfies the hypothesis.
\end{itemize}
\end{enumerate}}

\end{document}